# Some classical formulæ for curves and surfaces
## by Thomas Dedieu

The goal of this text is to present the computation by Salmon, in the second half of the XIXth century (now a little less than two-hundred years ago), of various numbers enumerating planes with a prescribed tangency pattern with a smooth, and sufficiently general, surface $S$ in $\mathbf{P}^3$ (or, equivalently, of hyperplane sections of $S$ with prescribed singularities). Emblematic among these are the number $\check{t}$ of tritangent planes, given in (3.17), and the number $\check{\beta}$ of planes cutting out a curve with a tacnode, given in (3.15).

The aforementioned numbers are the degrees of various loci in the dual projective space $\check{\mathbf{P}}^3$ (the points of which correspond to planes in $\mathbf{P}^3$), among which the closures of the various equisingular strata of the dual surface $S^\vee$ (the points of the latter correspond to planes tangent to $S$). The results we present feature not only the computation of these degrees, but also the description of the intricate geometry of these loci and their counterparts on $S$. Among these loci, one finds the parabolic curve of $S$ (the locus of points at which the tangent plane cuts out a cuspidal curve; it corresponds to the cuspidal double curve of $S^\vee$), the node-couple curve (the locus of points at which the tangent plane is also tangent at another point; it corresponds to the ordinary double curve of $S^\vee$), and the flecnodal curve (the locus of points $p$ at which the tangent plane cuts out a curve with one local branch with a flex at $p$; its counterpart on $S^\vee$ is not visible in terms of equisingular strata). A fundamental aspect of this geometry is the duality relation between $S$ and its dual, and more generally between the loci in each of them.

Salmon's arguments are specific to hypersurfaces, and based on the manipulation of polynomials and elimination theory. For instance, he determines relatively explicit polynomials cutting out the node-couple and flecnodal curves on $S$ (the parabolic curve is cut out by the Hessian determinant). He also makes an essential use of polars, which for hypersurfaces may be defined in terms of polynomial calculus, as we recall in Appendix A.1. His method to compute the number $\check{t}$ of tritangent planes, which we present in Sections 3.2 and 3.3, is as follows. Let $\varpi$ be a general point of $\check{\mathbf{P}}^3$, and consider the intersection of $S^\vee$ with its first and second polars with respect to $\varpi$. It is supported on the locus of points $\varkappa$ such that the line $\langle \varpi, \varkappa \rangle$ has intersection multiplicity at least three with $S^\vee$ at $\varkappa$, which of course comprises the triple points of $S^\vee$, and among those the points corresponding to tritangent planes. Then, a fine analysis of this intersection provides a set of relations involving the enumerative numbers we are interested in (but not only them), see Theorem (3.7), and it turns out we are able to compute sufficiently many of the quantities in these relations in some other way (for instance, the aforementioned number $\check{\beta}$ is obtained by intersecting the parabolic and flecnodal curves on $S$) so that we can finally find $\check{t}$.

We also present, still following Salmon, the computation of enumerative numbers related to the already mentioned loci, but not directly to the enumeration of hyperplane sections with prescribed singularities. For instance, we find the number of apparent double points of the cuspidal and ordinary double curve of $S^\vee$, namely $g_{\text{He}} = \check{h}$ and $g_{\text{Kn}} = \check{k}$ given in (4.17) and (4.26), respectively. In fact, we study (in Section 4) the two developable surfaces in $\mathbf{P}^3$ dual to the cuspidal and ordinary double curves of $S^\vee$, and compute their numerical characters. This involves a fine study of the local geometry of $S$ along its parabolic curve (see Section 4.2), which complements the aforementioned description of the geometry of $S$ and its duality with $S^\vee$. This also comprises an analysis analogous to that described above for the enumeration of tritangent planes, in which one considers this time the intersection of $S^\vee$ and its first polar with respect to $\varpi$ with another hypersurface derived from $S^\vee$ and relative to $\varpi$, such that the intersection is supported on the locus of points $\varkappa$ such that the line $\langle \varpi, \varkappa \rangle$ has intersection multiplicity with $S^\vee$ greater than one both at $\varkappa$ and at some other point.

The last Section 5 is dedicated to relatively more recent formulas (from the first half of the XXth century, now a little less than one-hundred years ago), following the presentation given by Enriques in his landmark book [15]. The main novelty is that they involve new invariants from a different nature (Salmon's formulas are all in term of the degree $n$ of the surface only), such as the arithmetic genus, and are applicable to surfaces in $\mathbf{P}^3$ with ordinary singularities, so that in practice they apply to any general three-dimensional linear subsystem of a very ample linear system on a smooth surface. This study includes in particular proofs of the Noether formula for the arithmetic genus of a smooth surface, and the classical postulation formula.

Appendix B contains a presentation of developable surfaces in $\mathbf{P}^3$ in classical terms, following Edge's beautiful book [14]. It features Plücker formulas for curves in $\mathbf{P}^3$, that are necessary for the computations in Section 4. Perhaps more importantly, it is the occasion to describe the duality relations between a curve in $\mathbf{P}^3$ and its dual, following the guiding thread of this text. This all has been generalized, in modern times, to curves in a projective space of arbitrary dimension by Piene [36], and her results are included as well.

Appendix C proposes an application of a basic topological formula counting singular curves in a pencil, see [XI, Section **], to the proof of a duality statement, namely that the orthogonality correspondence takes asymptotic tangent lines of a surface in $\mathbf{P}^3$ to asymptotic tangent lines of its dual.

For the two central Sections 3 and 4, I have followed rather closely Salmon's magnificent book [45], in its fifth edition of 1915 (I have not dug in earlier text of his, or prior versions of this book). I have done my best to add explanations and justifications wherever I felt this was needed (and this is a lot of places, if I may be so bold). Yet a blind spot of the present text is that it leaves the question of the validity of Salmon's formulas aside: they apply to a smooth surface under suitable transversality conditions, and we make all along the blanket assumption that a sufficiently general surface indeed verifies the latter transversality conditions, but do not prove it; see Paragraph (3.8) for some references on this question.

Salmon's formulas have been revisited in several modern texts already, in particular by Piene [38], using her theory of polars [37]. Besides, the title of this text is a nod at the title of [38], itself a nod at a series of papers by Roth. There is also a proof by Ronga [42], based on the use of multiple-point formulas: these ideas are discussed in [XIII]. Witaszek's thesis [53] is very complete and valuable, concentrated on the study of very general quartic surfaces. I have done my best to quote all relevant works on the topic. While I certainly have quoted all those that I took inspiration from, I claim by no means to have been encyclopedic. There are undoubtly many references missing (in particular, I have intentionally left aside some of the older material). May their authors and, more importantly, the interested readers forgive me.

Let me quote G.C. Rota to conclude this introduction: « Le teorie vanno e vengono ma le formule restano. »[1]

**Thanks.** I am grateful to Ragni Piene for her observations on a preliminary version of this text.

**General notation.** The exponent (and sometimes the superscript) $\vee$ denotes the dual, in the sense of vector spaces and projective spaces, or of projective geometry (see Section 1.1). The symbol $\perp$ denotes the orthogonal with respect to the pairing between linear functionals and vectors as well as, most often, the induced notion between points in two dual projective spaces (we will usually dodge the difference between vector spaces and projective spaces, thus trading off slight abuses for lighter notations and, hopefully, clarity). Thus, if $p$ is a point in $\mathbf{P}^N$, then $p^\perp$ is the hyperplane in the dual projective space $\check{\mathbf{P}}^N$ whose points correspond to all hyperplanes of $\mathbf{P}^N$ that contain $p$; and if $H$ is a hyperplane of $\mathbf{P}^N$, $H^\perp$ is the corresponding point of $\check{\mathbf{P}}^N$. For any subset $A \subseteq \mathbf{P}^N$, $\langle A \rangle$ denotes the span of $A$; thus, if $p$ and $q$ are two distinct points, $\langle p, q \rangle$ is the line spanned by $p$ and $q$.

For a vector space $V$, we denote by $\mathbf{P}(V)$ the projective space of all vector lines of $V$. If $W$ is a linear subspace of $V$, we use the shorthand $\mathbf{P}(V)/\mathbf{P}(W)$ for the projective space $\mathbf{P}(V/W)$, so that if $\Lambda$ is a linear subspace of some projective space $\mathbf{P}^N$, then $\mathbf{P}^N/\Lambda$ denotes the corresponding quotient projective space. We will then usually denote by $\pi_\Lambda : \mathbf{P}^N \dashrightarrow \mathbf{P}^N/\Lambda \cong \mathbf{P}^m$ the projection from $\Lambda$. We may abuse notation and denote, for any linear subspace $F$ of $\mathbf{P}^N/\Lambda$, by $\pi_\Lambda^{-1}(F)$ the largest linear subspace of $\mathbf{P}^N$ surjecting onto $F$.

---

[1] « The theories may come and go but the formulas remain. » Quote and translation stolen from Ottaviani's Five lectures on projective invariants.

If $p$ is a point of a scheme $X$, we denote by $\mathrm{TC}_p X$ the tangent cone of $X$ at $p$; if $p$ is a smooth point of $X$, we denote by $T_p X$ its tangent space at $p$; if moreover $X$ is seen as lying in some projective space, then $\mathbf{T}_p X$ denotes its projective tangent space at $p$.

If $X$ is a degree $n$ hypersurface in some projective space $\mathbf{P}^N$, we denote by $\mathrm{D}^{p^k} X$ its $k$-th polar with respect to the point $p$; it is a degree $n-k$ hypersurface.

The notation $O(d)$ denotes a linear combination of monomials which all have degree at least $d$.

The letter $n$ is reserved throughout the text to the degree of a hypersurface, with the (hopefully only) exception of Appendix A.

# 1 – Background material

## 1.1 – Projective duality

The theme of this text is the enumeration of hyperplanes with a prescribed tangency pattern with a smooth surface in $\mathbf{P}^3$. To get started, let us recall the basics of projective duality, that we will use repeatedly in the course of the text. For a more complete introduction, see [49, Chap. 1].

**(1.1) Definition.** *Let $X$ be a variety in $\mathbf{P}^N$. The* dual variety *of $X$ is the Zariski closure $X^\vee$ in the dual projective space $\check{\mathbf{P}}^N$ of the set of points $\varpi \in \check{\mathbf{P}}^N$ such that the hyperplane $\varpi^\perp \subseteq \mathbf{P}^N$ is tangent to $X$ at one of its smooth points.*

A central result is the following *reflexivity theorem.* For an algebraic proof, see [26] (which also features a detailed history of this result), or [21, Example 16.20]. For a proof using symplectic geometry, see [19, Chap. 1].

**(1.2) Theorem.** *Let $X$ be a variety in $\mathbf{P}^N$. One has $(X^\vee)^\vee = X$. More precisely, let $p$ and $\varpi$ be smooth points of $X$ and $X^\vee$ respectively: The hyperplane $\varpi^\perp \subseteq \mathbf{P}^N$ is tangent to $X$ at $p$ if and only if the hyperplane $p^\perp \subseteq \check{\mathbf{P}}^N$ is tangent to $X^\vee$ at $\varpi$.*

This may be refined in the following way. For a proof, see [25], or [52, §14.1].

**(1.3) Theorem.** *Suppose $X$ is smooth and its dual is a hypersurface. A point $\varpi \in X^\vee$ is smooth if and only if the hyperplane $\varpi \subseteq \mathbf{P}^N$ is tangent to $X$ in a unique point $p$, and $\varpi^\perp \cap X$ has an ordinary double point at $p$.*

The general philosophy is that the singularity of the dual in a point $\varpi$ reflects the tangency scheme of the hyperplane $\varpi^\perp$ with $X$. I shall loosely refer to any result incarnating this philosophy as a *biduality statement.* For instance, when $X$ itself is a smooth hypersurface one has the following result, see [12, Prop. 11.24].

**(1.4) Theorem.** *Let $X$ be a smooth hypersurface in $\mathbf{P}^N$, and let $\varpi \in X^\vee$ be such that the singular locus of $\varpi^\perp \cap X$ consists of finitely many points $p_1, \ldots, p_r$. Then the tangent cone to $X^\vee$ at $\varpi$ is the union of all hyperplanes $p_1^\perp, \ldots, p_r^\perp$, each counted with multiplicity $\mu(X, \varpi, p_i)$, the Milnor number of the singularity of $\varpi^\perp \cap X$ at $p_i$.*

The following result (see [49, Tm. 1.21]) says that linear projections and taking linear sections are dual operations. Recall from elementary linear algebra that if $W$ is a linear subspace of a vector space $V$, then the transpose of the projection $V \to V/W$ canonically identifies with the injection $W^\perp \hookrightarrow V^\vee$.

**(1.5) Proposition.** *Let $X$ be a variety in $\mathbf{P}^N$, and $\Lambda$ be a general linear subspace of codimension $m + 1$ in $\mathbf{P}^N$. Consider the projection from $\Lambda$,*

$$\pi_\Lambda : \mathbf{P}^N \dashrightarrow \mathbf{P}^N/\Lambda \cong \mathbf{P}^m.$$

i) *If $m > \dim(X)$ and $\pi_\Lambda|_X$ is an isomorphism to its image, then $\pi_\Lambda(X)^\vee = \Lambda^\perp \cap X^\vee$.*
ii) *If $m = \dim(X)$, let $B \subseteq \mathbf{P}^m$ be the branch divisor of $\pi_\Lambda|_X$. One has $B^\vee = \Lambda^\perp \cap X^\vee$.*
iii) *If $m < \dim(X)$, let $\tilde{X} \to X$ be the blow-up of $X$ at $\Lambda \cap X$, and $\tilde{\pi}_\Lambda : \tilde{X} \to \mathbf{P}^m$ be the morphism induced by $\pi_\Lambda$. Let $\Delta \subseteq \mathbf{P}^m$ be the discriminant of $\tilde{\pi}_\Lambda$. One has $\Delta^\vee = \Lambda^\perp \cap X^\vee$.*

## 1.2 – Projective duality for hypersurfaces and polarity

When $X$ is a hypersurface, much information on $X^\vee$ may be obtained by studying the defining equation of $X$ using polynomial calculus. We gather here a few facts relevant for us, and refer to [13, Chap. 1] for a much more complete treatment. A key technical tool is polarity, the main results of which are included in Appendix A.1.

In this section, we always assume that $X$ is a hypersurface in $\mathbf{P}^N$, defined by a degree $n$ homogeneous polynomial $F$ in a system of homogeneous coordinates $(x_0 : \ldots : x_N)$.

**(1.6)** Let $q = (q_0 : \ldots : q_N)$ be a smooth point of $X$. By Euler's formula, the equation of the tangent hyperplane to $X$ at $q$ is

$$\partial_0 F(q) \cdot x_0 + \cdots + \partial_N F(q) \cdot x_N = \mathrm{D}_q F(x_0, \ldots, x_N) = \mathrm{D}^{(x_0,\ldots,x_N)} F(q) = 0,$$

in the notation of (A.1). In particular, the point $p = (x_0 : \ldots : x_N)$ sits on $\mathbf{T}_q X$ if and only if $q$ sits in the degree $n - 1$ hypersurface $\mathrm{D}^p X$.

**(1.7)** Suppose for simplicity that $X$ is smooth. It may be seen from (1.6) that for all $p \in \mathbf{P}^N$, the intersection $\mathrm{D}^p X \cap X$ is the locus of those points $q \in X$ such that the line $\langle p, q \rangle$ is tangent to $X$ at $q$. Thus, $\mathrm{D}^p X \cap X$ is the ramification divisor of the projection $\pi_p : X \to \mathbf{P}^{N-1}$. It is called the *apparent boundary* of $X$ as seen from the point $p$; I will denote it by $\mathcal{B}^p X$. The cone projecting the apparent boundary from $p$ is called the *circumscribed cone*; I will denote it by $\mathcal{C}^p X$. When $X$ is singular, the circumscribed cone is the closure of the set of all lines passing through $p$ and tangent to $X$ at a smooth point; we will use this definition in Paragraph (3.5).

**(1.8) Gauss map and degree of the dual.** Since there is a unique tangent hyperplane at each smooth point of $X$, there is a dominant rational map

$$\gamma_X : X \dashrightarrow X^\vee$$

called the *Gauss map*. If $X^\vee$ is a hypersurface, it follows from (1.3) that the Gauss map is birational.

It follows from (1.6) that the Gauss map of $X$ is the restriction to $X$ of the map defined on $\mathbf{P}^N$ by the $N$-dimensional linear system of all first polars of $X$, i.e., of all $\mathrm{D}^p X$, $p \in \mathbf{P}^N$. This is a linear subsystem of $|\mathcal{O}_{\mathbf{P}^N}(n - 1)|$, and it is base point free if and only if $X$ is smooth.

This gives a way to compute the degree of the dual of $X$, which is also called the *class* of $X$. If $X$ is smooth, one finds that the degree of $X^\vee$ equals $n(n-1)^{N-1}$. If $X$ is singular, this has to be corrected by taking into account the base locus of the linear system of first polars, which is supported on the singular locus of $X$. One may prove the Plücker formula $(P_1)$ of (1.12) in this way. More generally, Teissier [48, Appendice II] has obtained a formula for the class of any hypersurface with isolated singularities; the reader may also consult [13, Section 1.2.3]. The same formula has been proved by Laumon [33], with a method close in spirit to the topological method described in [XI, Section **]. This formula has been generalized to hypersurfaces with arbitrary singularities by Piene [37].

**(1.9)** Assume that $X$ is smooth. The Hessian hypersurface $\mathrm{Hess}(X)$ (see Appendix A.2) cuts out on $X$ the locus of those points $p \in X$ such that the tangent section $\mathbf{T}_p X \cap X$ has a singularity at $p$ worse than an ordinary double point. The intersection $X \cap \mathrm{Hess}(X)$ is the ramification divisor of the Gauss map $\gamma_X : X \to X^\vee$.

Before I close this section, let me mention that there exists a notion of polars for varieties in a projective space that are not necessarily hypersurfaces. See, for instance, [37] and [18]. Using these one may compute the class of varieties in projective space other than hypersurfaces, see, e.g., [37] and [51].

### 1.3 – Plücker formulae for plane curves

As a first and fundamental application of the previous results, let us explain the classical Plücker formulae for plane curves.

**(1.10)** Let $C \subseteq \mathbf{P}^2$ be a plane curve. It has a finite number of bitangent lines (i.e., lines tangent to $C$ in two distinct points) and inflection points, which give rise respectively to nodes and cusp of the dual curve $C^\vee \subseteq \check{\mathbf{P}}^2$: this can be seen for instance with Theorem (1.4). In general there are no lines with a larger tangent scheme with $C$, so that nodes and cusps are the only singularities of $C^\vee$.

The reflexivity theorem (1.2) suggests to look for formulae linking numerical characters of $C$ and $C^\vee$ in a symmetric fashion. In this perspective one should allow $C$ to have nodes and cusps.

**(1.11) Plückerian characters of a plane curve.** We consider a plane curve $C \subseteq \mathbf{P}^2$ such that both $C$ and its dual $C^\vee \subseteq \check{\mathbf{P}}^2$ have only ordinary double points (nodes) and ordinary cusps as singularities. We call:
— $n$ the degree of $C$;
— $\check{n}$ the degree of $C^\vee$;
— $\beta$ the number of bitangents of $C$, which is also the number of nodes of $C^\vee$;
— $\iota$ the number of flexes of $C$, which is also the number of cusps of $C^\vee$;
— $\delta$ the number of nodes of $C$, which is also the number of bitangents of $C^\vee$;
— $\kappa$ the number of cusps of $C$, which is also the number of flexes of $C^\vee$.

**(1.12) Plücker formulae for plane curves.** One has the following relations for a curve $C$ as in (1.11).

$$(P_1) \quad \check{n} = n(n-1) - 2\delta - 3\kappa; \qquad (\check{P}_1) \quad n = \check{n}(\check{n}-1) - 2\beta - 3\iota;$$
$$(P_2) \quad \iota = 3n(n-2) - 6\delta - 8\kappa; \qquad (\check{P}_2) \quad \kappa = 3\check{n}(\check{n}-2) - 6\beta - 8\iota,$$

For a more symmetric version one may prefer to write $\check{\delta}$ and $\check{\kappa}$ for $\beta$ and $\iota$, respectively.

We sketch the proof, as it is emblematic of our approach in this text, and refer to [16] for a detailed implementation.

One first computes the class of a plane curve with nodes and cusps, following (1.8): this will give the two relations on the first line by applying the computation to both $C$ and $C^\vee$. Let $p \in \mathbf{P}^2$ be a general point. It follows from Theorem (A.10) that the first polar $\mathrm{D}^p C$ passes through all nodes and cusps of $C$, is smooth at these points, and has general tangent line at the nodes and tangent line equal to the reduced tangent cone of $C$ at the cusps. Hence $\mathrm{D}^p C$ intersects $C$ with multiplicity $\geqslant 2$ at the nodes and $\geqslant 3$ at the cusps, and in fact these are the actual intersection multiplicities. The relation follows.

Then one computes the number of flexes of $C$ by considering its intersection with its Hessian (see (1.9)): this will give the two relations on the second row. Again, one has to substract the contribution of $C \cap \mathrm{Hess}(C)$ corresponding to the singularities of $C$ and not to actual inflection points. If $q$ is a node of $C$, then the Hessian curve also has an ordinary double point at $q$, with tangent cone equal to that of $C$. If $q$ is a cusp of $C$, then the Hessian curve has a triple point with two local branches at $q$, one smooth with general tangent, and one cuspidal with tangent cone equal to that of $C$. It follows that $C$ intersects its Hessian with multiplicity $\geqslant 6$ at the nodes and $\geqslant 8$ at the cusps, and again these are the actual intersection multiplicities, hence the formula.

**(1.13) Generators for the Plücker relations.** The four relations in (1.12) generate a codimension three complete intersection ideal, i.e., the Plücker formulae reduce to three formulae only, with the Koszul type relations as the only relations among them. Indeed, the following relation holds:
$$P_1 + \check{P}_1 = \frac{1}{3}(P_2 + \check{P}_2),$$
so that any three among $P_1, P_2, \check{P}_1, \check{P}_2$ generate all Plücker relations. It is arguably an endless game to search for a preferred set of three generators for the Plücker relations; the three relations below form an appealing generating set:

$$\begin{align}
(P_1') \quad & 3n - \kappa = 3\check{n} - \iota \\
(P_2') \quad & n(n-2) + \check{n}(\check{n}-2) = 2\delta + 3\kappa + 2\beta + 3\iota \\
(P_3') \quad & 18(\delta - \beta) = (\kappa - \iota)(\kappa - \iota + 6\check{n} - 27).
\end{align}$$

They are indeed relations as

$$\begin{align}
P_1' &= 3P_1 - P_2 = -3\check{P}_1 + \check{P}_2 \\
P_2' &= P_1 + \check{P}_1 = \frac{1}{3}(P_2 + \check{P}_2) \\
P_3' &= 3\left[3(n+\check{n}) + \kappa - \iota - 3\right] P_1 - \left[3(n+\check{n}) + \kappa - \iota\right] P_2 + 9\check{P}_1,
\end{align}$$

and they indeed generate the Plücker relations as

$$\begin{pmatrix} P_1 \\ P_2 \\ \check{P}_1 \\ \check{P}_2 \end{pmatrix} = \begin{pmatrix} \frac{1}{6}(n+\check{n}) + \frac{1}{18}(\kappa-\iota) & \frac{1}{2} & -\frac{1}{18} \\ \frac{1}{2}(n+\check{n}) + \frac{1}{6}(\kappa-\iota) - 1 & \frac{3}{2} & -\frac{1}{6} \\ -\frac{1}{6}(n+\check{n}) - \frac{1}{18}(\kappa-\iota) & \frac{1}{2} & \frac{1}{18} \\ -\frac{1}{2}(n+\check{n}) - \frac{1}{6}(\kappa-\iota) + 1 & \frac{3}{2} & \frac{1}{6} \end{pmatrix} \begin{pmatrix} P_1' \\ P_2' \\ P_3' \end{pmatrix}.$$

**(1.14) Smooth plane curves.** When the plane curve $C$ is smooth one has $\delta = \kappa = 0$, so one is left with four Plückerian characters and three independent Plücker relations. The two characters $\check{n}$ and $\iota$ are expressed in terms of $n$ by the relations $P_1$ and $P_2$, and one finds in addition:
$$\beta = \tfrac{1}{2}n(n-2)(n-3)(n+3).$$

**(1.15) The genus formula.** The Plücker formulae usually come adjoined with the following relation, gotten from the fact that the two curves $C$ and $\check{C}$, being birational, have the same geometric genus. Computing their respective geometric genera by means of the adjunction formula, one gets:

$$\tfrac{1}{2}(n-1)(n-2) - \delta - \kappa = \tfrac{1}{2}(\check{n}-1)(\check{n}-2) - \beta - \iota.$$

The latter relation may in fact be obtained from the Plücker relations: it is

$$\tfrac{1}{2}(-P_1 - 2\check{P}_1 + \check{P}_2) = -\tfrac{1}{2}\left[\tfrac{1}{3}(n+\check{n}) + \tfrac{1}{9}(\kappa - \iota) - 1\right] P_1' + \tfrac{1}{18} P_3'.$$

The condition that a plane curve and its dual both have nodes and cusps as their only singularities poses some delicate questions. We refer to [32] for recent results. Finally, let me mention [31], in which the author introduces a concept of virtual numbers of nodes and cusps, which allows him to give symmetric Plücker formulas valid for curves with arbitrary singularities.

## 2 – Double curves of the dual of a surface in $\mathbf{P}^3$

We now start our study of surfaces in $\mathbf{P}^3$. Here we give a first description of the dual of a smooth surface in $\mathbf{P}^3$, that will be refined later.

### 2.1 – Local geometry of a surface in $\mathbf{P}^3$ and its dual

We first take the occasion to introduce some necessary notions on the local geometry of a surface in $\mathbf{P}^3$. Let $S$ be a surface in $\mathbf{P}^3$.

**(2.1) Parabolic points.** A *parabolic point* of $S$ is a point $p \in S$ such that the Hessian of $S$ at $p$ vanishes (see (A.11)).

Equivalently, the polar quadric of $S$ at $p$, viz. $\mathrm{D}_{p^2}S$, is singular, in the sense that it is either a singular quadric indeed, or the whole $\mathbf{P}^3$. If $p$ is a smooth point of $S$, then the first (resp. second) possibility is equivalent to the tangent hyperplane section $S \cap \mathbf{T}_p S$ having a double point at $p$ with tangent cone a double line (resp. multiplicity strictly larger than 2 at $p$). Thus, in general, the tangent hyperplane section of $S$ at a parabolic point $p$ will have a cusp at $p$.

If the Hessian of $S$ does not vanish identically on $S$ (this is the case if and only if the dual of $S$ is a surface, see (B.5)), we shall call the curve $S \cap \mathrm{Hess}(S)$ the *parabolic curve* of $S$, or indifferently the *Hessian curve*, and denote it by $\mathrm{He}(S)$.

**(2.2) Asymptotic tangent lines.** A line $\Lambda \subseteq \mathbf{P}^3$ is an *asymptotic tangent* of $S$ at the point $p \in S$ if $\Lambda$ and $S$ have intersection multiplicity at least 3 at $p$.

By Theorem (A.8), the asymptotic tangent lines at $p$ are the lines contained in the intersection $\mathbf{T}_p S \cap \mathrm{D}_{p^2}S$ of the tangent hyperplane and the polar quadric of $S$ at $p$ (note that $\mathbf{T}_p S = \mathrm{D}_p S$ and $\mathrm{D}_{p^2}S$ are tangent at $p$ by Theorem (A.10)). Thus, there are two distinct asymptotic tangent lines at $p$ if $p$ is not a parabolic point, and one double asymptotic tangent lines at $p$ if $p$ is a parabolic point and the polar quadric $\mathrm{D}_{p^2}S$ is not the whole $\mathbf{P}^3$. If $\mathrm{D}_{p^2}S$ is the whole $\mathbf{P}^3$, then all lines in $\mathbf{T}_p S$ passing through $p$ are asymptotic tangents.

If the polar quadric $\mathrm{D}_{p^2}S$ is not the whole $\mathbf{P}^3$, then the sum of the asymptotic tangent lines at $p$ is the tangent cone at $p$ of the tangent hyperplane section $\mathbf{T}_p S \cap S$.

*(2.2.1) Flex tangent lines.* The asymptotic tangent lines may also be called *flex tangent lines*. Let $H \subseteq \mathbf{P}^3$ be a hyperplane containing $p \in S$. The curve $H \cap S$ has a flex at $p$ if and only if the plane $H$ contains an asymptotic tangent line of $S$ at $p$, as an elementary local computation shows.

*(2.2.2) A duality statement for asymptotic tangents.* For all smooth $p \in S$, if $\Lambda \subseteq \mathbf{P}^3$ is an asymptotic tangent line of $S$ at $p$, then the line $\Lambda^\perp \subseteq \check{\mathbf{P}}^3$ is an asymptotic tangent line of $S^\vee$ at $(\mathbf{T}_p S)^\perp$.

This is proved in [C, **], and in Appendix C; see also Proposition (4.6).

**(2.3) Codimension one singular loci of the dual.** Let $S$ be a general surface in $\mathbf{P}^3$. Then its dual $S^\vee$ is singular along the two curves parametrizing, respectively, bitangent hyperplanes of $S$, and hyperplanes tangent to $S$ at a parabolic point.

I will call the former the *ordinary double curve* of $S^\vee$. The surface $S^\vee$ indeed has an ordinary double point at its generic point by Theorem (1.4). The two local sheets of $S^\vee$ along the ordinary double curve correspond to deformations of bitangent hyperplanes preserving either one of the two tangency points. The *node-couple curve* $\mathrm{Kn}(S) \subseteq S$ is the closure of the locus of points $p \in S$ such that $\mathbf{T}_p S$ is a bitangent plane. Thus, the Gauss map realizes $\mathrm{Kn}(S)$ as a double cover of the ordinary double curve of $S^\vee$.

I will call *cuspidal double curve* of $S^\vee$ the image of the parabolic curve of $S$ by the Gauss map. The surface $S^\vee$ indeed has an $A_2$ double point at its generic point, as follows from the fact that the ramification divisor of the birational map $S \dashrightarrow S^\vee$ is the parabolic curve, see [13, Proposition 1.2.1].

We will make throughout the text the blanket assumption that the codimension one singular loci of the dual surface $S^\vee$ amount eaxctly to its ordinary and cuspidal double curves, to which we will refer as the *double curves* of $S^\vee$. We refer to Paragraph (3.8) for a discussion of this assumption.

## 2.2 – Degrees of the double curves

This section is dedicated to formulae that give the degrees of the ordinary and cuspidal double curves of the dual of a surface $S$ in $\mathbf{P}^3$; the main result is Corollary (2.7). These formulae have been given a modern proof in [24] (see, in particular, Section VII.1 therein).

The underlying idea is that the double curves of the dual $S^\vee$ amount to behaviours occuring in codimension two in the linear system of hyperplane sections of $S$, hence these behaviours are visible in a general codimension one linear subsystem; therefore, one should be able to find all the relevant information after projecting $S$ from a general point.

**(2.4) Setup.** Let $S$ be a smooth, degree $n$, surface in $\mathbf{P}^3$. For all $p \in \mathbf{P}^3$, we let $\mathcal{B}^p S$ be the apparent boundary of $S$ from $p$ (since $S$ is smooth, this is the curve $S \cap \mathrm{D}^p S$), and consider the plane curve $B_p \subseteq \mathbf{P}^2$ image of $\mathcal{B}^p S$ by the projection $\pi_p : \mathbf{P}^3 \dashrightarrow \mathbf{P}^2$ from the point $p$; it is the branch divisor of the projection $\pi_p|_S : S \to \mathbf{P}^2$.

The main result of [10] (see also [30]) asserts that for general $p$, the curve $B_p$ has only nodes and cusps as singularities (see [20] for a general result valid in all dimensions), and it is irreducible. In addition, we shall assume that the dual plane curve $B_p^\vee$ has only nodes and cusps as singularities as well: since $B_p^\vee$ is isomorphic to $p^\perp \cap S^\vee$, this is equivalent to the assumption that the codimension one singular loci of $S^\vee$ are as described in (2.3).

**(2.5) Plückerian characters of the branch curve.** Let $\nu$ be the degree of the branch curve $B_p \subseteq \mathbf{P}^2$, $\check{\nu}$ be its class, and $\delta, \kappa, \beta, \iota$ be its other Plückerian characters as in (1.12). These numbers may all be expressed as polynomials in the degree of $S$:

$$\nu = n(n-1) \qquad\qquad \check{\nu} = n(n-1)^2$$
$$\delta = \tfrac{1}{2}n(n-1)(n-2)(n-3) \qquad\qquad \beta = \tfrac{1}{2}n(n-1)(n-2)(n^3 - n^2 + n - 12)$$
$$\kappa = n(n-1)(n-2) \qquad\qquad \iota = 4n(n-1)(n-2).$$

**(2.6)** Let me first outline a quick proof of these formulae, based on the two facts that, on the one hand, $B_p$ is the birational projection from $p$ of the apparent boundary $\mathcal{B}^p S$, and on the other hand, it is the plane dual of the section of $S^\vee$ by the plane $p^\perp \subseteq \check{\mathbf{P}}^3$.

This implies that: i) the degree and class of $B_p$ equal the degrees of $\mathcal{B}^p S$ and $S^\vee$ respectively, which gives the formulae for $\nu$ and $\check{\nu}$ (see (2.8) below): ii) the geometric genus of $B_p$ equals that of $\mathcal{B}^p S$, which is easy to compute since the latter is a complete intersection. We can then derive all remaining characters using the Plücker formulae for plane curves. This is carried out in a more general context in Section 5.3.

Since $B_p$ is a general hyperplane section of $S^\vee$, the above formulae give the following.

**(2.7) Corollary.** *The numbers of binodal and cuspidal curves in a general net of hyperplane sections of $S$ are respectively*

$$\tfrac{1}{2}n(n-1)(n-2)(n^3 - n^2 + n - 12) \quad \text{and} \quad 4n(n-1)(n-2).$$

*These numbers are the degrees of the ordinary and cuspidal double curves of $S^\vee$ respectively.*

**(2.8)** Let me now give a description of the Plückerian characters in (2.5) in terms of the tangent hyperplanes to $S$. It gives a way to compute these characters directly (although for $\beta$ I don't know of an elementary way to do it), regardless of the quick proof given in (2.6).

*(2.8.1) Degree of $B_p$.* The curve $B_p$ is the projection of the circumscribed curve $\mathcal{B}^p S = S \cap \mathrm{D}^p S$ from the outer point $p$. Therefore, $\deg(B_p) = \deg(\mathcal{B}^p S) = \deg(S) \cdot \deg(\mathrm{D}^p S)$.

*(2.8.2) Class of $B_p$.* The curve $B_p^\vee \subseteq \check{\mathbf{P}}^2$ is the hyperplane section of $S^\vee$ by $p^\perp$. Indeed, for all $q \in B_p$, the hyperplane $\pi_p^{-1}(\mathbf{T}_q B_p)$ is tangent to $S$ at some point $q' \in \mathcal{B}^p S$ and contains $p$; conversely, if $H \subseteq \mathbf{P}^3$ contains $p$ and is tangent to $S$ at some point $q'$, then $q' \in \mathcal{B}^p S$ and $\pi_p(H)$ is the tangent line to $B_p$ at $\pi_p(q')$.
This implies that $\deg(B_p^\vee) = \deg(S^\vee) = \deg(S) \cdot \deg(\mathrm{D}^p S)^2$.

*(2.8.3) Nodes of $B_p$.* The nodes of $B_p$ are the points of $B_p$ with two distinct preimages in $\mathcal{B}^p S$, hence they correspond to unordered pairs of distinct points $q, q' \in \mathcal{B}^p S$ such that $p, q, q'$ are aligned. Since $\mathcal{B}^p S$ is the curve of points $q \in S$ such that the line $\langle p, q \rangle$ is tangent to $S$ at $q$, the nodes of $B_p$ correspond to the bitangent lines to $S$ passing by $p$. The number of such lines may be computed in a natural way using discriminants; this is explained in [C, Section **].

*(2.8.4) Flexes of $B_p$.* Since $B_p$ is the branch curve of the projection of $S$ from $p$, it has a flex at a point $q$ if and only if the hyperplane $\pi_p^{-1}(\mathbf{T}_q B_p)$ cuts out a curve with a cuspidal singularity on $S$.

This may be seen with the following local computation. Let $q$ be an inflection point of $B_p$, and $q' \in \mathcal{B}^p S$ be such that $\pi_p(q') = q$. Locally around $q'$, $S$ may be described as a double cover of a neighbourhood of $q$ in $\mathbf{P}^2$, hence as a hypersurface $z^2 = f(x, y)$ in affine coordinates $(x, y, z)$, such that the projection from $p$ is $(x, y, z) \mapsto (x, y)$ and $f(x, y) = 0$ is an equation of $B_q$ locally at $q$. Assume that the tangent line $\mathbf{T}_q B_p$ is $y = 0$. Then the section of $S$ by its tangent plane at $q'$ is the hypersurface $z^2 = f(x, 0)$ in affine coordinates $(x, z)$, and since $q$ is a flex of $B$, $f(x,0)$ vanishes to the third order at 0, hence the curve $z^2 = f(x,0)$ has a cusp at the origin, as we wanted to show. The converse is proved in the same fashion.

It follows that the number of flexes of $B_p$ is the number of hyperplanes passing through $p$ and tangent to $S$ at a parabolic point. This is the number of intersection points of $\mathcal{B}^p S$ with the parabolic curve, i.e., $\deg(S) \cdot \deg(\mathrm{D}^p S) \cdot \deg(\mathrm{Hess}(S))$.

*(2.8.5) Cusps of $B_p$.* The cusps of $B_p$ correspond to the points $q \in \mathcal{B}^p S$ such that the line $\langle p, q \rangle$ is tangent to $\mathcal{B}^p S$ at $q$, see Lemma (B.9). By Lemma (2.9) below, this happens if and only if the line $\langle p, q \rangle$ is an asymptotic tangent of $S$ at $q$. Therefore, the number of cusps of $B_p$ is $\deg(S) \cdot \deg(\mathrm{D}^p S) \cdot \deg(\mathrm{D}^{p^2} S)$.

*(2.8.6) Bitangents of $B_p$.* A line $\Lambda \subseteq \mathbf{P}^2$ is bitangent to $B_p$ if and only if the plane $\pi_p^{-1}(\Lambda) \subseteq \mathbf{P}^3$ is bitangent to $S$. In [45, §605–607] Salmon constructs from the equation of $S$, by rather elaborate elimination theory, a polynomial which cuts out on $S$ the locus of tangency points of all bitangent planes of $S$. This polynomial has degree $(n-2)(n^3 - n^2 + n - 12)$, and from this one deduces the number of bitangent planes of $S$ that pass through a given general point, and therefore the number of bitangent lines of $B_p$, see Theorem ** and Corollary ** in [C]. This construction by Salmon is the main topic of [C].

**(2.9) Lemma.** *Let $q \in \mathcal{B}^p S$. The line $\langle p, q \rangle$ is tangent to $\mathcal{B}^p S$ at $q$ if and only if it is an asymptotic tangent of $S$ at $q$.*

This may be proved in several ways that I deem interesting. One may first of all do it by a local computation: this is done in [10, Prop. 3.9], and I will leave it as an exercise. I shall now give two other proofs.

*Proof of (2.9) by polarity.* By assumption, $q \in \mathrm{D}^p S$, and the line $\langle p, q \rangle$ is tangent to $S$ at $q$. Then, the line $\langle p, q \rangle$ is an asymptotic tangent of $S$ at $q$ if and only if $q \in \mathrm{D}^{p^2} S$, which by polar symmetry is equivalent to $p \in \mathrm{D}^{q^{n-2}} S$. On the other hand, being tangent to $S$, the line $\langle p, q \rangle$ is tangent to $\mathcal{B}^p S$ at $q$ if and only if it is tangent to $\mathrm{D}^p S$, which is equivalent to $p \in \mathrm{D}^{q^{n-2}}(\mathrm{D}^p S)$ (note that $\mathrm{D}^p S$ is a degree $n-1$ hypersurface). Now, $\mathrm{D}^{q^{n-2}}(\mathrm{D}^p S) = \mathrm{D}^p(\mathrm{D}^{q^{n-2}} S)$, and the two conditions $p \in \mathrm{D}^{q^{n-2}} S$ and $p \in \mathrm{D}^p(\mathrm{D}^{q^{n-2}} S)$ are equivalent by the Euler formula, see (A.6). □

*Proof of (2.9) by duality.* We shall use the duality statement that the orthogonal of an asymptotic tangent to $S$ at $q$ is an asymptotic tangent to $S^\vee$ at $(\mathbf{T}_q S)^\perp$, see (2.2.2).

Thus, in the conditions of the lemma, the line $\langle p, q \rangle$ is an asymptotic tangent of $S$ at $q$ if and only if the line $\langle p, q \rangle^\perp$ is an asymptotic tangent of $S^\vee$ at $(\mathbf{T}_q S)^\perp$. Since the plane $p^\perp$ contains the line $\langle p, q \rangle^\perp$, this is equivalent to the hyperplane section $p^\perp \cap S^\vee$ having a flex along the line $\langle p, q \rangle^\perp$ at the point $(\mathbf{T}_q S)^\perp$. In turn, since the projection $B_p$ of $\mathcal{B}^p S$ from $p$ is dual as a plane curve to $p^\perp \cap S^\vee$, this is equivalent to $B_p$ having a cusp at the projection of $q$ from $p$, and the latter condition is equivalent, as we have already noted, to the line $\langle p, q \rangle$ being tangent to $\mathcal{B}^p S$ at $q$. □

Let $q \in \mathcal{B}^p S$ be a smooth point. We have seen above that, on the one hand, $B_p$ has a flex at $\pi_p(q)$ if and only if $q$ is a parabolic point of $S$, and then the tangent $\mathbf{T}_q(\mathcal{B}^p S)$ is the unique asymptotic tangent of $S$ at $q$, and on the other hand, $B_p$ has a cusp at $\pi_p(q)$ if and only if the line $\langle p, q \rangle$ is the tangent $\mathbf{T}_q(\mathcal{B}^p S)$, and this happens if and only if $\langle p, q \rangle$ is an asymptotic tangent at $q$. With this respect, it is interesting to note the following.

**(2.10) Lemma.** *Let $q \in \mathcal{B}^p S$ be a smooth point. The tangent line $\mathbf{T}_q(\mathcal{B}^p S)$ is an asymptotic tangent of $S$ at $q$ if and only if either $q$ is a parabolic point, or $\langle p, q \rangle$ is an asymptotic tangent at $q$. In the latter case, $\langle p, q \rangle$ is the tangent line $\mathbf{T}_q(\mathcal{B}^p S)$.*

*Proof.* First note that since $q$ is a smooth point of $\mathcal{B}^p S$, the linear spaces $\mathbf{T}_q S$ and $\mathbf{T}_q(\mathrm{D}^p S)$ intersect along a line. This implies that $\mathbf{T}_q(\mathrm{D}^p S)$ has dimension 2, and since it equals $\mathrm{D}^{q^{n-2}} \mathrm{D}^p S = \mathrm{D}^p \mathrm{D}^{q^{n-2}} S$, the polar quadric of $S$ at $q$ is not the whole $\mathbf{P}^3$.

Assume that $q$ is parabolic. Then $\mathrm{D}^{q^{n-2}} S$ is a quadric cone, smooth at $q$ by Theorem (A.10), and $S$ has a unique flex tangent line at $q$, which is the line $L = \mathbf{T}_q(\mathrm{D}^{q^{n-2}} S) \cap \mathrm{D}^{q^{n-2}} S = \mathbf{T}_q S \cap \mathrm{D}^{q^{n-2}} S$ spanned by $q$ and the vertex of the cone $\mathrm{D}^{q^{n-2}} S$. The tangent line $\mathbf{T}_q(\mathcal{B}^p S)$ equals $\mathbf{T}_q S \cap \mathbf{T}_q(\mathrm{D}^p S)$, so we only need to show that the line $L$ is contained in $\mathbf{T}_q(\mathrm{D}^p S)$. Now, since $\mathbf{T}_q(\mathrm{D}^p S) = \mathrm{D}^p \mathrm{D}^{q^{n-2}} S$, it contains the singular point of the cone $\mathrm{D}^{q^{n-2}} S$, hence also the whole line $L$.

If, on the other hand, $\langle p, q \rangle$ is an asymptotic tangent at $q$, then it equals $\mathbf{T}_q(\mathcal{B}^p S)$ by Lemma (2.9) above, and we are done in this case.

Conversely, assume that the line $L = \mathbf{T}_q(\mathcal{B}^p S)$ is an asymptotic tangent at $q$, but $\langle p, q \rangle$ is not. Thus $p$ is off $L$, and $L$ is contained in the quadric $\mathrm{D}^{q^{n-2}} S$. Then, for all $r \in L$, the tangent space $\mathbf{T}_r(\mathrm{D}^{q^{n-2}} S)$ contains $L$, but also $p$ because $r \in L \subseteq \mathrm{D}^p(\mathrm{D}^{q^{n-2}} S)$. It follows that the plane $\langle p, L \rangle$ is tangent to $\mathrm{D}^{q^{n-2}} S$ at all points $r \in L$, and this implies that $\mathrm{D}^{q^{n-2}} S$ is a cone, i.e., $q$ is a parabolic point of $S$. □

## 3 – Zero-dimensional strata of the dual surface

In this section we study some zero-dimensional strata of the dual of a general smooth surface $S$ in $\mathbf{P}^3$; in particular, we compute the number of tritangent planes of $S$, see (3.17). The general strategy is to look at the cuspidal edges of the polar cone of $S^\vee$ with respect to a general point of $\check{\mathbf{P}}^3$; the main result in this direction is Theorem (3.7). The ordinary double edges of this cone are studied in Section 4.

### 3.1 – Local models

In Section 3.2 below we shall study a batch of formulae valid for surfaces in $\mathbf{P}^3$ with a given class of singularities that are modeled on those of the dual of a general enough smooth surface. Here we describe these singularities. A discussion about why these should be the appropriate singularities may be found in Paragraph (3.8) below.

**(3.1) Swallowtails.** We define the *swallowtail[2] singularity* (also referred to in this text as "point of type $\beta$") as the singularity at the origin of the discriminant locus in the semi-universal deformation space of an $A_3$ curve singularity, also known as an ordinary tacnode.

The semi-universal deformation of the tacnode is the family of curves over the affine space $\mathbf{A}^3_{a,b,c}$, of dimension 3 and with coordinates $(a, b, c)$, defined in $\mathbf{A}^2_{x,y} \times \mathbf{A}^3_{a,b,c}$ by the equation:

$$(3.1.1) \qquad y^2 = x^4 + ax^2 + bx + c.$$

For more details, please see [II]. The semi-universal deformation space of the tacnode is this $\mathbf{A}^3_{a,b,c}$, and the *discriminant locus* therein is the locus over which the curves in the family defined by (3.1.1) For all $(a, b, c)$, the curve in $\mathbf{A}^2_{x,y}$ defined by Equation (3.1.1) is singular if and only if the polynomial in $x$ on the right-hand-side has a double root. Therefore the equation of the discriminant locus is given by the classical discriminant of a univariate degree 4 polynomial, i.e.,

$$(3.1.2) \qquad \mathrm{Disc}(x^4 + ax^2 + bx + c) = 256c^3 - 128a^2c^2 + 144ab^2c - 27b^4 + 16a^4c - 4a^3b^2.$$

---

[2]in french, « queue d'arronde », one of seven elementary catastrophes in Thom's theory.

From this equation we see that the swallowtail singularity is a locally irreducible triple point, with tangent cone a triple plane.

Alternatively, the discriminant locus is the image of the map

$$(u,v) \longmapsto (x-u)^2(x-v)(x+2u+v)$$

(note that the absence of $x^3$ term on the right-hand-side of (3.1.1) is equivalent to the sum of the four roots of the polynomial in $x$ being 0). It follows that the equation of the discriminant locus may be obtained by eliminating $u$ and $v$ from the equations

$$\begin{cases} a = -3u^2 - 2uv - v^2 \\ b = 2u^3 + 4u^2v + 2uv^2 \\ c = -2u^3v - u^2v^2 \end{cases}$$

gotten by computing the coefficients of $(x-u)^2(x-v)(x+2u+v)$.

Similarly, the locus in the semi-universal deformation space of the tacnode over which the curves in the family defined by (3.1.1) have a cusp is the image of the map $(u,v) \mapsto (x-u)^3(x+3u)$, hence its equations may be found by eliminating $u$ and $v$ from the equations

$$a = -6u^2; \quad b = 8u^3; \quad c = -3u^4.$$

This gives the ideal of equations

(3.1.3) $$(a^2 + 12c, 8a^3 + 27b^2),$$

from which we see that the locus in consideration is an irreducible curve with a planar ordinary cusp at the origin, the generic point of which is a cuspidal double point of the discriminant locus. We shall call this curve the *cuspidal double curve* of the swallowtail singularity.

In turn, the locus over which the curves in the family defined by (3.1.1) have two nodes is the image of the map $(u,v) \mapsto (x-u)^2(x+u)^2$, hence its equations may be found by eliminating $u$ and $v$ from the equations

$$a = -2u^2; \quad b = 0; \quad c = u^4;$$

this gives the ideal of equations

(3.1.4) $$(b, a^2 - 4c),$$

from which we see that the locus in consideration is a smooth irreducible curve, the generic point of which is an ordinary double point of the discriminant locus. We shall call this curve the *ordinary double curve* of the swallowtail singularity.

Finally, one finds that the tangent cone of the intersection of the cuspidal and ordinary double curves is defined by the ideal $(b, c, a^2)$, so these two curves intersect with multiplicity 2.

**(3.2) Type $\gamma$ singularity.** Similarly to what we have done for the swallowtail in (3.1) above, we define the *singularity of type $\gamma$* as the singularity at the origin of the discriminant locus in the semi-universal deformation space of the multisingularity $(A_1, A_2)$ consisting of a node and an ordinary cusp.

The semi-universal deformations of a node and a cusp are respectively given by

(3.2.1) $$y^2 = x^2 + a \quad \text{and} \quad y^2 = x^3 + bx + c,$$

and we consider the disjoint union of their pull-back to the product $\mathbf{A}_a^1 \times \mathbf{A}_{b,c}^2$. The equation of the discriminant is merely the product of the two discriminants, i.e.,

(3.2.2) $$\operatorname{Disc}(x^2 + a) \cdot \operatorname{Disc}(x^3 + bx + c) = a \cdot (4b^3 + 27c^2).$$

At the origin, the discriminant locus has two local sheets, one smooth and one with a cuspidal double curve, and the two meet transversely.

The locus in $\mathbf{A}_a^1 \times \mathbf{A}_{b,c}^2$ above which the curve defined by (3.2.1) has a cusp parametrizes those deformations of the multisingularity $(A_1, A_2)$ along which the node smoothes and the cusp is maintained. Its has equations
$$b = c = 0.$$

It is a smooth curve, which we will call the *cuspidal double curve* of the type $\gamma$ singularity.

The locus in $\mathbf{A}_a^1 \times \mathbf{A}_{b,c}^2$ above which the curve defined by (3.2.1) has two nodes parametrizes those deformations of the multisingularity $(A_1, A_2)$ along which the node is maintained and the cusp deforms to a node, which has equations
$$a = 4b^3 + 27c^2 = 0.$$

It is a curve with a planar ordinary cusp at the origin, which we will call the *ordinary double curve* of the type $\gamma$ singularity. It is the intersection curve of the two sheets of the discriminant surface.

**(3.3) Type $t$ singularity.** Finally, we define the *singularity of type $t$* as the singularity at the origin of the discriminant locus in the semi-universal deformation space of the multisingularity $(A_1, A_1, A_1)$ consisting of three nodes.

Thus we consider three copies of the semi-universal deformations of a node,

(3.3.1) $$y^2 = x^2 + a \quad \text{and} \quad y^2 = x^2 + b \quad \text{and} \quad y^2 = x^2 + c,$$

and look at the disjoint union of their pull-back to the product $\mathbf{A}_a^1 \times \mathbf{A}_b^1 \times \mathbf{A}_c^1$. The equation of the discriminant is the product of the three discriminants, i.e.,

(3.3.2) $$\operatorname{Disc}(x^2 + a) \cdot \operatorname{Disc}(x^2 + b) \cdot \operatorname{Disc}(x^2 + c) = abc.$$

The discriminant locus has three smooth local sheets meeting transversely at the origin.

The locus in $\mathbf{A}_a^1 \times \mathbf{A}_b^1 \times \mathbf{A}_c^1$ above which the curve defined by (3.3.1) has two nodes is the union of the three coordinate lines $\mathbf{A}_a^1 \times \{0\} \times \{0\}$, $\{0\} \times \mathbf{A}_b^1 \times \{0\}$, and $\{0\} \times \{0\} \times \mathbf{A}_c^1$; it is defined by the ideal of equations $(ab, ac, bc)$, and has a triple point at the origin. Yet again, we shall call it the ordinary double curve of the type $t$ singularity.

The locus in $\mathbf{A}_a^1 \times \mathbf{A}_b^1 \times \mathbf{A}_c^1$ above which the curve defined by (3.3.1) has a cusp is empty.

## 3.2 – Polar cones and their cuspidal edges

In this section we describe the cuspidal edges of the polar cone of a surface in $\mathbf{P}^3$ which has the same kind of singularities as the dual of a sufficiently general smooth surface. We shall call this class of singularities d-ordinary singularities, since (plain) ordinary singularities classically refer to the kind of singularities found on a sufficiently general projection to $\mathbf{P}^3$ of a smooth surface in $\mathbf{P}^5$ (or in $\mathbf{P}^N$ for any $N \geqslant 5$), see (5.3).

The central result is Theorem (3.7), which is the keystone of the computation of the number of tritangent planes in the next section.

**(3.4) Surfaces with d-ordinary singularities.** We say that a surface $X \subseteq \mathbf{P}^3$ has *d-ordinary singularities* if it has a curve $B$ of ordinary double points, a curve $C$ of simple cuspidal double points, and no other singularities except from those arising from the singularities of the curves $B$ and $C$, by which we mean:
— $\beta$ intersection points between $B$ and $C$ at which $C$ is cuspidal (and $B$ is smooth), at which $X$ has a type $\beta$ (i.e., swallowtail) singularity as in (3.1) above;
— $\gamma$ intersection points between $B$ and $C$ at which $B$ is cuspidal (and $C$ is smooth), at which $X$ has a type $\gamma$ singularity as in (3.2) above;
— $t$ the number of triples points of the ordinary double curve $B$, at which $X$ has a type $t$ singularity as in (3.3) above.
Moreover, we let $b$ and $c$ denote the degrees of the curves $B$ and $C$ respectively.

We now consider, for the rest of the present Section 3.2, a surface $X \subseteq \mathbf{P}^3$ with d-ordinary singularities, with the notation as in (3.4) above.

Let $p$ be a general point of $\mathbf{P}^3$. It follows from a local computation, see Theorem (A.10), that the first polar $\mathrm{D}^p X$ contains simply the curves $B$ and $C$, with tangent cone equal to that of $X$ at the generic point of $C$. Therefore, the intersection $X \cap \mathrm{D}^p X$ contains $B$ and $C$ with multiplicities 2 and 3 respectively, and we need to refine the definition of apparent boundary given in (1.7).

**(3.5) Polar cones of a surface with d-ordinary singularities.** We let the *polar cone* $\mathcal{D}^p X$ of $X$ with respect to $p$ be the cone projecting from $p$ over the curve $X \cap \mathrm{D}^p X$. It contains the two cones projecting $B$ and $C$ from $p$ with multiplicities 2 and 3 respectively, and we let the *circumscribed cone* $\mathcal{C}^p X$ be the cone residual to these two cones in the polar cone; it is the closure of the set of all lines passing through $p$ and tangent to $X$ at a smooth point. Finally, we denote by $\mathcal{B}^p X$ the curve $(X \cap \mathcal{C}^p X)_{\mathrm{red}}$, which we call the *apparent boundary*.

We let $a$ be the degree of $\mathcal{C}^p X$. Then, the degree of the polar cone $\mathcal{D}^p X$ reads

$$(3.5.1) \qquad n(n-1) = a + 2b + 3c.$$

**(3.6) Cuspidal edges of the polar cone.** We now introduce more notation, and list various kinds of points that will appear in the intersection $X \cap \mathrm{D}^p X \cap \mathrm{D}^{p^2} X$, together with the singular points listed in (3.4). The lines joining $p$ to the points of $X \cap \mathrm{D}^p X \cap \mathrm{D}^{p^2} X$ are the cuspidal edges of the polar cone $\mathcal{D}^p X$, in analogy with what happens when $X$ is smooth, see Paragraph (2.8.5).

We set:
— $\kappa$ the number of cuspidal edges of the circumscribed cone $\mathcal{C}^p X$;
— $\rho$ the number of intersection points between the apparent boundary $\mathcal{B}^p X$ and the ordinary double curve $B$;
— $\sigma$ the number of intersection points between the apparent boundary $\mathcal{B}^p X$ and the cuspidal double curve $C$.
We shall introduce some more notation later in (4.18), but it is not relevant for the present considerations.

Salmon gives the following set of formulae [45, §610], which describes separately the intersections with the second polar $\mathrm{D}^{p^2} X$ of the three components of the polar cone $\mathcal{D}^p X$.

**(3.7) Theorem.** *For $X \subseteq \mathbf{P}^3$ a surface with d-ordinary singularities and $p \in \mathbf{P}^3$ a general point, and with the notation of (3.4) and (3.6), one has:*

(☙ i) $\qquad\qquad\qquad\qquad\qquad a(n-2) = \kappa + \rho + 2\sigma \; ;$

(☙ ii) $\qquad\qquad\qquad\qquad\qquad b(n-2) = \rho + 2\beta + 3\gamma + 3t \; ;$

(☙ iii) $\qquad\qquad\qquad\qquad\qquad c(n-2) = 2\sigma + 4\beta + \gamma.$

    He comments: « The reader can see without difficulty that the points indicated in these formulæ are included in the intersections of $[\mathrm{D}^{p^2}X$ with the three components of $\mathcal{D}^pX]$; but it is not so easy to see the reason for the numerical multipliers which are used in the formulæ. Although it is probably not impossible to account for these constants by a priori reasoning, I prefer to explain the method by which I was led to them inductively. » The latter method consists in considering the dual of a smooth cubic surface, for which all the quantities may be directly computed, and then in tuning the coefficients to make the formulae work.

    The remainder of the present Section 3.2 is dedicated to the proof of the above theorem. By Theorem (A.8), the intersection

$$X \cap \mathrm{D}^p X \cap \mathrm{D}^{p^2} X = \bigl(\mathcal{B}^p X + 2B + 3C\bigr) \cap \mathrm{D}^{p^2} X$$

is supported on the locus $\mathfrak{D}$ of points $q$ in $X$ such that the line $\langle p, q \rangle$ has intersection multiplicity at least 3 with $X$ at $q$. The proof consists in the following two points. (i) List the various kinds of points in the locus $\mathfrak{D}$; this indicates which terms (without the integer coefficients) appear on the right-hand-side of Formulae (☙); the condition that $X$ has d-ordinary singularities leaves us with finitely many possibilities. (ii) For each of those, determine the corresponding multiplicity in the intersection cycle $X \cap \mathrm{D}^p X \cap \mathrm{D}^{p^2} X$; this will tell us the integer coefficient with which each term on the right-hand-side of (☙) appears.

    The smooth points in the locus $\mathfrak{D}$ are the points of type $\kappa$, see (3.7.1) below. Some special double points of $X$ lie in $\mathfrak{D}$, and these are the points in the two intersections $\mathcal{B}^p X \cap B$ and $\mathcal{B}^p X \cap C$. Finally, all triple points automatically lie in $\mathfrak{D}$. Thus the points appearing in $\mathfrak{D}$ are those of the types $\kappa, \rho, \sigma, \beta, \gamma, t$. The remainder of the proof consists in a case by case local study of the polars of $X$ at all such points.

*(3.7.1) Type $\kappa$.* Cuspidal edges of $\mathcal{C}^p X$ appear in correspondence with points $q \in \mathcal{B}^p X$ such that the line $\langle p, q \rangle$ is tangent to $\mathcal{B}^p X$, as has been explained in (2.8.5), and such points $q$ are simple points in the intersection $\mathcal{B}^p X \cap \mathrm{D}^{p^2} X$ as in Lemma (2.9). For general $p$, these points lie on neither one of the double curves of $X$.

*(3.7.2) Type $\rho$.* Points of type $\rho$ are those points $q$ lying on the ordinary double curve $B$ such that the line $\langle p, q \rangle$ is an honest tangent to one of the two smooth local sheets of $X$ at $q$; let us call $X'_q$ this local sheet. We claim that (i) the first polar $\mathrm{D}^p X$ is smooth at $q$ with the same tangent plane as $X'_q$, and its intersection with $X$ locally at $q$ consists of $B$ and $\mathcal{C}_p X$, the latter entirely contained in $X'_q$, and (ii) the second polar $\mathrm{D}^{p^2} X$ is smooth at $q$ with general tangent plane, hence it intersects both $B$ and $\mathcal{C}_p X$ transversely at $q$.

    Let us justify this claim by a specific local computation, as it does not directly follow from Theorem (A.10). We may assume that $q = (0:0:0:1)$ in a homogeneous system of coordinates $(x:y:z:w)$, and $X$ is given locally at $q$ by the equation

$$f(x, y, z) = \bigl[x + g_2(x, y, z) + O(3)\bigr]\bigl[y + h_2(x, y, z) + O(3)\bigr] = 0$$

in the affine chart $w = 1$, where $g_2$ and $h_2$ are homogeneous of degree 2, and $O(3)$ denotes higher order terms. Then we may take $p = (1:0:0:0)$, and thus $\mathrm{D}^p = \partial_x$. Then,

$$\mathrm{D}^p f(x,y,z) = \left[x + g_2 + O(3)\right]\left[\partial_x h_2 + O(2)\right] + \left[1 + \partial_x g_2 + O(2)\right]\left[y + h_2 + O(3)\right]$$
$$= y + x\,\partial_x h_2 + y\,\partial_x g_2 + h_2 + O(3)$$

and part (i) of the claim follows. Finally,

$$\mathrm{D}^{p^2} f(x,y,z) = 2\,\partial_x h_2 + x\,\partial_x^2 h_2 + y\,\partial_x^2 g_2 + O(2),$$

and part (ii) follows.

*(3.7.3) Type $\sigma$*. Points of type $\sigma$ are those points $q$ lying on the cuspidal double curve $C$ such that the line $\langle p, q \rangle$ lies in the tangent cone of the unique, cuspidal, local sheet of $X$ at $q$. We claim that (i) the first polar $\mathrm{D}^p X$ has a double point at $q$, and the intersection of the tangent cones $\mathrm{TC}_q X$ and $\mathrm{TC}_q(\mathrm{D}^p X)$ is a quadruple line supported on $\mathbf{T}_q C$; moreover, (ii) the second polar $\mathrm{D}^{p^2} X$ has a simple point at $q$, and the intersection of its tangent plane with $[\mathrm{TC}_q X]_{\mathrm{red}}$ is the line $\mathbf{T}_q C$. By (i), since at $q$ the intersection of $\mathrm{D}^p X$ and $X$ is the sum of $3C$ and $\mathcal{B}^p X$, the apparent boundary $\mathcal{B}^p X$ is smooth at $q$, and tangent to the cuspidal edge $C$. Thus, (ii) tells us that the second polar $\mathrm{D}^{p^2} X$ is tangent to both curves $C$ and $\mathcal{B}^p X$ at $q$, which accounts for the multiplicities 2 of $\sigma$ in (✎ i) and (✎ iii).

As in (3.7.2) above, we justify our claim by a local computation, with similar notation. We let $q = (0:0:0:1)$, and consider the following equation for $X$ locally at $q$ in the affine chart $w = 1$:

$$f(x,y,z) = \left[y + g_2(x,y,z) + O(3)\right]^2 + \left[x + h_2(x,y,z) + O(3)\right]^3 = 0,$$

where $g_2$ and $h_2$ are homogeneous of degree 2. We may again take $p = (1:0:0:0)$, and then

$$\mathrm{D}^p f = 2\,\partial_x g_2 \cdot y + 3x^2 + O(3), \qquad \text{and} \qquad \mathrm{D}^{p^2} f = 2\,\partial_x^2 g_2 \cdot y + 6x + O(2),$$

which proves our claim.

*(3.7.4) Type $t$*. These are the points at which $X$ has three transverse smooth local sheets. We may directly apply Theorem (A.10), to the effect that at a type $t$ singular point $q$, (i) the first polar $\mathrm{D}^p X$ has an ordinary double point at $q$, and (ii) the second polar $\mathrm{D}^{p^2} X$ is smooth with general tangent plane, hence intersects $B$ transversely at $q$. This gives the multiplicity 3 for type $t$ in (✎ iv), as $q$ is a triple point of $B$. On the other hand, type $t$ points lie neither on $\mathcal{B}^p X$ nor on $C$.

*(3.7.5) Type $\beta$*. These are triple points of $X$, with a local equation equivalent to that in (3.1.2); in particular the tangent cone of $X$ at such a point is a triple plane. Thus, it follows from Theorem (A.10) that if $q$ is a type $\beta$ singular point, then (i) the first polar $\mathrm{D}^p X$ has a double point at $q$ with tangent cone a double plane supported on $(\mathrm{TC}_q X)_{\mathrm{red}}$, and (ii) the second polar $\mathrm{D}^{p^2} X$ is smooth with tangent plane equal to $(\mathrm{TC}_q X)_{\mathrm{red}}$.

The ordinary double curve $B$ is smooth at points of type $\beta$, as we have seen in Paragraph (3.1), and the second polar is tangent to it. This accounts for the multiplicity 2 in (✎ i). On the other hand, the cuspidal double curve $C$ has planar ordinary cusps at points of type $\beta$, but it turns out that it intersects $\mathrm{D}^{p^2} X$ with multiplicity 4 at such points: indeed, in the notation of (3.1), $C$ is defined by the ideal $(a^2 + 12c, 8a^3 + 27b^2)$, and the tangent cone of $\mathrm{D}^{p^2} X$ is the plane $c = 0$. This confirms the multiplicity of $\beta$ in (✎ iv), and shows that $(\mathrm{TC}_q X)_{\mathrm{red}}$ is the osculating plane of $C$ at swallowtail points $q$.

*(3.7.6) Type $\gamma$.* These are triple points $q$ of $X$, at which $X$ has two transverse local sheets, one smooth: $X_q^{\text{s}}$, and the other cuspidal: $X_q^{\text{c}}$. The first polar $\mathrm{D}^p X$ has a double point at $q$, with tangent cone the sum of $[\mathrm{TC}_q(X_q^{\text{c}})]_{\text{red}}$ and a general plane containing the line $\mathrm{TC}_q(B)$. The second polar $\mathrm{D}^{p^2} X$ is smooth at $q$, with tangent plane a general plane containing the line $\mathrm{TC}_q(B)$. Therefore it intersects $B$, which is cuspidal at $q$, with multiplicity 3, and it is transverse to $C$, which is smooth at $q$. This amounts for the multiplicities 3 and 1 in (🐛 ⸝) and (🐛 ⸜⸜⸜) respectively.

To clarify the situation, let us carry out the local computations proving the above assertions. We let $q = (0:0:0:1)$ and consider the two sheets $X_q^s$ defined respectively by the equations

$$z + O(2) = 0 \quad \text{and} \quad y^2 + x^3 + O(4) = 0.$$

The curve $B$ is the intersection of these two sheets, it has a cusp at $q$ with tangent cone defined by $z = y^2 = 0$. The curve $C$ is the singular locus of $X_q^{\text{c}}$, it is smooth at $q$ with tangent line defined by $y = x = 0$. We may take $p = (u_0 : u_1 : u_2 : 0)$. Then, by Theorem (A.10) the tangent cones of $\mathrm{D}^p X$ and $\mathrm{D}^{p^2} X$ are defined respectively by

$$\mathrm{D}^p(zy^2) = y(2u_1 z + u_2 y) \quad \text{and} \quad \mathrm{D}^{p^2}(zy^2) = 2u_1(u_1 z + 2u_2 y),$$

which proves our claims.

The proof of Theorem (3.7) is now complete. □

## 3.3 – Application to the case when $X = S^\vee$

We now let $S$ be a smooth, degree $n$, surface in $\mathbf{P}^3$, and call $X = S^\vee \subseteq \check{\mathbf{P}}^3$ its dual surface. We shall explore in this situation the geometric landscape described in Section 3.2 above.

**(3.8)** We take for granted that if $S$ is general enough, then $X$ has d-ordinary singularities (see (3.4)), and we assume that $S$ is indeed general enough.

The principle behind this implication is that, since hyperplane sections of $S$ form a 3-dimensional system, we expect that under suitable generality/transversality assumptions, only (multi-)singularities imposing at most three conditions may occur on hyperplane sections of $S$. The condition that a given multi-singularity imposes at most three conditions translates into the fact that its semi-universal deformation space has dimension at most three (for a general discussion of the expected codimension imposed by the presence of a multi-singularity, see [XIII, Section **]). We expect moreover that, still under the aforementioned generality assumptions, for all $H^\perp \in \check{\mathbf{P}}^3$, there is a submersion from an analytic or étale neighbourhood of $H^\perp$ in $\check{\mathbf{P}}^3$ to the semi-universal deformation space of the multi-singularity of $H \cap S$; this implies that, locally around $H^\perp$, the dual $S^\vee$ looks like the discriminant locus in the deformation space of the multi-singularity of $H \cap S$.

For a proof that if $S$ is sufficiently general, then things happen according to the above described expectations, see [6], [34], and [43]. See also [54] for a related result.

**(3.9) Notation.** As a rule, we use the decoration with a duality symbol to indicate that a symbol denotes a quantity associated with $X = S^\vee$. Thus, we set:
— $\check{n} = n(n-1)^2$, the degree of $X$.
We use the notation introduced in Section 3.2 for the singularities of $X = S^\vee$, with an additional "$\vee$" decoration.

**(3.10)** The double curves $\check{B}$ and $\check{C}$ are the curves in $S^\vee$ parametrizing the bitangent and parabolic tangent planes to $S$, respectively. The Gauss map of $S$ induces a $2:1$ map $\mathrm{Kn}_S \to \check{B}$ from the node-couple curve $\mathrm{Kn}_S \subseteq S$, and a birational map $\mathrm{He}_S \to \check{C}$ from the Hessian curve $\mathrm{He}_S \subseteq S$. We have seen in Corollary (2.7) that the degrees of the two curves $\check{B}$ and $\check{C}$ in $\check{\mathbf{P}}^3$ are:
— $\check{b} = \frac{1}{2} n(n-1)(n-2)(n^3 - n^2 + n - 12)$ ;
— $\check{c} = 4n(n-1)(n-2)$.

**(3.11)** In the context of Section 3.2, we need to choose a general point $\check{p} \in \check{\mathbf{P}}^3$ to consider the polar cone $\mathcal{D}^{\check{p}} X$ and the apparent boundary $\mathcal{B}^{\check{p}} X$. This point $\check{p}$ corresponds to a hyperplane $H_{\check{p}} = \check{p}^\perp$ in $\mathbf{P}^3$. The circumscribed cone $\mathcal{C}^{\check{p}}(S^\vee) \subseteq \check{\mathbf{P}}^3$ is the dual of the plane curve $H_{\check{p}} \cap S \subseteq \mathbf{P}^3$. This gives:
— $\check{a} = n(n-1)$.
Note that with our values of $\check{a}, \check{b}, \check{c}$, Formula (3.5.1) holds indeed.

It also follows from the above description of $\mathcal{C}^{\check{p}}(S^\vee)$ that the cuspidal edges of type $\check{\kappa}$ of the polar cone $\mathcal{D}^{\check{p}} S^\vee$ are the lines orthogonals to the osculating tangent lines of the plane curve $H_{\check{p}} \cap S$ (i.e., the orthogonals of the flex tangents of $H_{\check{p}} \cap S$). Thus,
— $\check{\kappa} = 3n(n-2)$.

**(3.12)** The apparent boundary $\mathcal{B}^{\check{p}} X$ parametrizes points $\varpi \in X = S^\vee$ such that $\check{p} \in \mathbf{T}_\varpi X$, i.e., dually, planes $\varpi^\perp \subseteq \mathbf{P}^3$ tangent to $S$ at a point $(\mathbf{T}_\varpi X)^\perp$ lying on $H_{\check{p}}$. It follows that the points of type $\check{\rho}$ (resp. $\check{\sigma}$) correspond to hyperplanes that are tangent to $S$ at some point on $H_{\check{p}}$ and at some other point (resp. at some point on $H_{\check{p}}$ which is parabolic). Thus, $\check{\rho}$ (resp. $\check{\sigma}$) is the degree of the intersection of $H_{\check{p}}$ with the node-couple curve $\mathrm{Kn}_S$ (resp. with the Hessian curve $\mathrm{He}_S$), and therefore:
— $\check{\rho} = n(n-2)(n^3 - n^2 + n - 12)$ ;
— $\check{\sigma} = 4n(n-2)$
(for $\check{\rho}$, see (2.8.6); $\check{\sigma}$ follows directly from the fact that $\mathrm{He}_S$ is the complete intersection of $S$ and $\mathrm{Hess}(S)$).

**(3.13)** Points of type $\check{\gamma}$ correspond to planes $H$ tangent to $S$ at two points, one of which is parabolic; correspondingly the curve $H \cap S$ has one cusp and one node. Points of type $\check{\gamma}$ are honest points on $\check{C}$, in account of the cusp in $H \cap S$; the presence of an additional, remote, node does not affect the local geometry of the curve $\hat{C}$. On the other hand, points of type $\check{\gamma}$ are cusps of $\check{B}$: if one watches a general point of $\check{B}$ move towards a point of type $\check{\gamma}$, one sees a family of hyperplane sections of $S$ with two nodes such that one the nodes degenerates to a cusp, and this is the reason why points of type $\check{\gamma}$ are cuspidal points of $\check{B}$, much in the same way as the curve in $S^\vee$ corresponding to cuspidal tangent sections of $S$ is a cuspidal curve of $S^\vee$. In less handwaving terms, the local geometry of $\check{B}$ and $\check{C}$ at points of type $\check{\gamma}$ is justified by the discussions in (3.2) and (3.8). The number $\check{\gamma}$ shall be derived indirectly in (3.16) below.

Points of type $\check{\beta}$ correspond to planes $H$ tangent to $S$ at one point, such that the section $H \cap S$ has a tacnode there. Indeed, points of type $\check{\beta}$ are cusps of $\check{C}$, so that as a general point of $\check{C}$ moves towards a point of type $\check{\beta}$, one must see a family of hyperplane sections of $S$ with a cusp which degenerates to a tacnode. Going to a point of type $\check{\beta}$ along $\check{B}$, one merely sees two nodes tending to a common limit point, which is harmless for $\check{B}$, as indeed points $\check{\beta}$ are smooth points of $\check{B}$. The number $\check{\beta}$ will be derived in (3.15) by considering the flecnodal locus in $S$.

Finally, points of type $\check{t}$ correspond to the famous tritangent planes of $S$. Their number we shall derive using formula (☞ ⋈).

**(3.14) The flecnodal locus.** The *flecnodal locus* is the locus of points $q \in S$ such that there exists a line having intersection multipliciy at least 4 with $S$ at $q$. It is a curve on $S$, and for a general point $q$ on this curve, the tangent section $\mathbf{T}_q S \cap S$ has an ordinary double point at $q$, such that one of the two local branches has a flex at $q$: the tangent line to this branch is the one line having intersection multiplicity 4 with $S$ at $q$.

It has been shown by Salmon [45, §588] that the flecnodal curve is cut out on $S$ by a polynomial of degree $11n - 24$, which we shall call the *flecnodal polynomial*. We describe his proof in [C, Section **]. A modern proof, and a generalization to arbitrary dimension, is given in [7]. One can also find a brief account in [27, Section 8], and another generalization in [1]. Remarkably, the flecnodal polynomial vanishes modulo the equation of $S$ if and only if $S$ is ruled (note that the analogous vanishing property for the Hessian characterizes developable ruled surfaces): Salmon attributes this result to Cayley; for a proof, see [27, Section 8], or the references given in [7].

In [45, §589–595], Salmon gives the following explicit formula for the flecnodal polynomial, which he attributes to Clebsch. Let $F$ be an equation for $S$, and denote by $\mathcal{H}$ its Hessian matrix. Set
$$\Theta = \sum_{i,j} |\mathcal{H}|_{i,j} \, \partial_i F \, \partial_j F,$$
where $|\mathcal{H}|_{i,j}$ denotes the cofactor formed by removing the $i$-th row and the $j$-th column of $\mathcal{H}$, i.e., $(-1)^{i+j}$ times the corresponding $3 \times 3$ minor. Set
$$\Phi = -\sum_{i,j} \sum_{\substack{i_1 < i_2 \\ j_1 < j_2}} |\mathcal{H}|_{i,j} \, |\mathcal{H}|_{i_1 i_2, j_1 j_2} \, \partial_{i i_1 i_2} F \, \partial_{j j_1 j_2} F$$
where, in addition to the previous notation, $|\mathcal{H}|_{i_1 i_2, j_1 j_2}$ is $(-1)^{i_1+i_2+j_1+j_2}$ times the obvious $2 \times 2$ minor of $\mathcal{H}$. Then, the flecnodal polynomial is

(3.14.1) $\qquad\qquad\qquad\qquad \Theta - 4\Phi \operatorname{Hess}(F).$

**(3.15) Intersection of the parabolic and flecnodal curves.** The points $q \in S$ such that the tangent section $\mathbf{T}_q S \cap S$ has a tacnode at $q$ are the intersection points of the parabolic curve with the flecnodal curve. Indeed, a point $q$ sits on the parabolic curve if and only if $\mathbf{T}_q S \cap S$ has a double point at $q$ with degenerate tangent cone (or a point of higher multiplicity), and it sits on the flecnodal curve if and only if there is a line in the plane $\mathbf{T}_q S$ having intersection multiplicity at least 4 with the curve $\mathbf{T}_q S \cap S$ at $q$. With our blanket transversality assumption from (3.8), this happens if and only if $\mathbf{T}_q S \cap S$ has a tacnode at $q$.

The upshot of the above discussion is that the number $\check{\beta}$ is the number of intersection points of the parabolic and flecnodal curves, which are cut out by polynomials of degrees $4(n-2)$ and $11n - 24$, respectively. There is however one last subtlety, namely that the latter two curves are everywhere tangent, see (3.15.1). Therefore, one finds:
— $\check{\beta} = 2n(n-2)(11n - 24)$.

*(3.15.1)* The fact that the parabolic and flecnodal curves are everywhere tangent has been proved by Salmon [45, §596]; a different proof has been given by McCrory and Shifrin [34] in modern times. In fact, Salmon even proves that the two surfaces cutting out the parabolic and flecnodal curves on $S$ are everywhere tangent.

His argument is the following. Taking into account the expression (3.14.1) for the flecnodal polynomial, it suffices to show that $\Theta$ is a square modulo the Hessian determinant $\operatorname{Hess}(F)$. Now, it follows from a determinantal identity that for all 4-tuples $(c_i) \in \mathbf{C}^4$,
$$\left(\sum |\mathcal{H}|_{i,j} \, \partial_i F \, \partial_j F\right)\left(\sum |\mathcal{H}|_{i,j} \, c_i \, c_j\right) = \left(\sum |\mathcal{H}|_{i,j} \, \partial_i F \, c_j\right)^2 \quad \mod \operatorname{Hess}(F).$$

From the latter identity, Salmon deduces without further indication that $\Theta$ is indeed a square modulo $\mathrm{Hess}(F)$. I reckon that some additional arguments ought to be provided, but let me not try to do it here. Funnily enough, the squared factor on the right-hand-side of the last identity has degree $7n - 15$, a quantity that will show up later in a seemingly unrelated manner, see (4.17).

**(3.16) Intersection of the two double curves.** We now analyze the intersection of the node-couple $\mathrm{Kn}_S$ and Hessian $\mathrm{He}_S$ curves on $S$. For each point $q$ in their intersection, the tangent plane $\mathbf{T}_q S$ is either of type $\check{\gamma}$ or of type $\check{\beta}$, and we shall see that the two curves are transverse if it is of type $\check{\gamma}$, and tangent if it is of type $\check{\beta}$.

To a tangent plane of type $\check{\gamma}$ there corresponds a section of $S$ with one node and one cusp, which gives two points on $\mathrm{Kn}_S$, but only one of those lies on $\mathrm{He}_S$. We add to our blanket transversality assumption that the two curves are transverse at this point.

To a tangent plane of type $\check{\beta}$ there corresponds a unique tangency point $q \in S$. On the one hand $q$ belongs to $\mathrm{Kn}_S$, and it is a ramification point of the double cover $\mathrm{Kn}_S \to \check{B}$, as $\mathbf{T}_q S \cap S$ may be seen as a bi-nodal curve with the two nodes coalesced. This implies that $T_q \mathrm{Kn}_S$ lies in the kernel of the differential of the Gauss map of $S$ at $q$. On the other hand $q$ is a ramification point of $\mathrm{He}_S \to \check{C}$ as well, since $\check{C}$ has a cusp at $(\mathbf{T}_q S)^\perp$, and this implies that $T_q \mathrm{He}_S$ lies in the kernel of the differential of the Gauss map of $S$ at $q$ (see also Lemma (4.9)). This implies that the two curves $\mathrm{Kn}_S$ and $\mathrm{He}_S$ are tangent at $q$, unless the differential of the Gauss map be $0$ at $q$, but this is not the case by our blanket transversality assumption.

In the upshot, we find that
$$\mathrm{Kn}_S \cdot \mathrm{He}_S = 2\check{\beta} + \check{\gamma}.$$

Moreover, using the fact that both $\mathrm{Kn}_S$ and $\mathrm{He}_S$ are complete intersections of $S$ and a hypersurface of known degree, we find

$$\mathrm{Kn}_S \cdot \mathrm{He}_S = (n-2)(n^3 - n^2 + n - 12) \cdot 4(n-2) \cdot n.$$

We may now plug in the formula for $\check{\beta}$ obtained in (3.13), and obtain:
— $\check{\gamma} = 4n(n-2)(n-3)(n^3 + 3n - 16)$.[3]

**(3.17) Formulae (✑) and the number of tritangent planes.** First of all let me mention that with the quantities we have been able to find so far, formulae (✑ ı) and (✑ ııı) are indeed verified (beware that $n$ should be replaced by $\check{n} = n(n-1)^2$).

Formula (✑ ıv) on the other hand yields the number of tritangent planes to a smooth surface of degree $n$ in $\mathbf{P}^3$:
— $\check{t} = \frac{1}{6} n(n-2)(n^7 - 4n^6 + 7n^5 - 45n^4 + 114n^3 - 111n^2 + 548n - 960)$.

## 4 – Hessian and node-couple developables

In this section we determine some more enumerative numbers related to 0-dimensional strata of a smooth surface $S \subseteq \mathbf{P}^3$ and its dual, thus complementing the analysis carried out in the previous Section 3. These numbers are not immediately related to the enumeration of hyperplane sections of $S$ with given singularities; for instance, we determine the number of flexes of the cuspidal double curve $\check{C} \subseteq S^\vee$ (this is $\beta_{\mathrm{He}}$ in (4.17)), and these correspond to planes cutting out an honest tacnodal curve on $S$. In general, we are concerned in this section with enumerative number associated to the Hessian and node-couple developables (for the definition, see Section 4.1 below).

---

[3]There is a typo in the formula for $\check{\gamma}$ in [45, §630] ($n^3 - 3n - 16$ instead of $n^3 + 3n - 16$), which curiously enough goes back and forth between the various editions of the book, as Ragni Piene has pointed out to me. The formula for $\check{\gamma}$ also appears in [45, §612], without the typo.

There are two main aspects to this study. On the one hand we shall describe (in Section 4.3) the ordinary double edges of the polar cone of a surface with d-ordinary singularities, thus complementing the analysis of cuspidal double edges in Section 3.2 above. On the other hand we will carry out (in Section 4.2) a detailed study of the geometry of a surface along its parabolic curve, which I believe is interesting for its own sake.

## 4.1 – Developables associated to curves on a surface

We refer to Appendix B for the basics on developable surfaces in $\mathbf{P}^3$.

**(4.1) Developables associated to curves on a surface.** Let $Y$ be a curve traced on a surface $S \subseteq \mathbf{P}^3$. Recalling that a developable surface in $\mathbf{P}^3$ merely amounts to a one-dimensional family of planes, one may consider the developable surface $\Sigma_Y \subseteq \mathbf{P}^3$ corresponding to the family of planes tangent to $S$ at points of $Y$. In classical terms: this is the family of planes touching $S$ at points of $Y$, and $\Sigma_Y$ is the surface enveloped by the family of tangent planes to $S$ at points of $Y$.

Let $\gamma_S : S \dashrightarrow S^\vee$ be the Gauss map of $S$. We are considering the family of planes parametrized by the curve $\gamma_S(Y) \subseteq \check{\mathbf{P}}^3$. By definition, the developable surface $\Sigma_Y$ is the dual surface $\gamma_S(Y)^\vee \subseteq \mathbf{P}^3$ of the curve $\gamma_S(Y) \subseteq \check{\mathbf{P}}^3$. There exists a curve $\Gamma \subseteq \mathbf{P}^3$ such that $\Sigma_Y$ is the tangential surface $\mathrm{Tan}(\Gamma)$, but there is no reason that $\Gamma$ should coincide with $Y$, and neither that it should lie on $S$.

**(4.2) Hessian and node-couple developables.** The *Hessian and node-couple developables* associated to a smooth surface $S \subseteq \mathbf{P}^3$ are the developable surfaces $\Sigma_Y$ as in (4.1) above with $Y$ the parabolic curve $\mathrm{He}(S)$ and node-couple curve $\mathrm{Kn}(S)$, repectively. They are the dual surfaces of the cuspidal and ordinary double curves of $S^\vee$, respectively, and we will denote them by $\Sigma_{\mathrm{He}}$ and $\Sigma_{\mathrm{Kn}}$.

**(4.3) Proposition.** *The generators of the node-couple developable surface are the lines joining pairs of conjugated points of the node-couple curve, i.e., the lines $\langle p, q \rangle$ such that there exists a plane which is tangent to $S$ at both $p$ and $q$.*

*Proof.* The node-couple developable is the dual of the curve $\gamma_S(\mathrm{Kn}(S)) \subseteq \check{\mathbf{P}}^3$, so its generators are the orthogonals of the tangents lines of $\gamma_S(K_S)$ (see Section B.1). Let $\varpi \in \gamma_S(K_S)$ be a general point. There are two points $p, q \in S$ such that the plane $H = \varpi^\perp \subseteq \mathbf{P}^3$ is tangent to $S$ at both $p$ and $q$. The tangent line $\mathbf{T}_\varpi(\gamma_S(K_S))$ is the intersection of the tangent planes to the two local sheets of $S^\vee$ at $\varpi$, respectively $p^\perp$ and $q^\perp$. Therefore,

$$\mathbf{T}_\varpi(\gamma_S(K_S))^\perp = (p^\perp \cap q^\perp)^\perp = \langle p, q \rangle.$$

□

We will determine the numerical characters of the node-couple developable in Section 4.4; we need the results of Section 4.3 on the ordinary double edges of the polar cone to do it.

In order to describe the generators of the developable $\Sigma_Y$ for a general curve $Y$, we introduce the following device.

**(4.4) The involution $\iota_p$ on $\mathbf{P}(T_pS)$.** Let $p$ be a smooth, non-parabolic, point of a surface $S \subseteq \mathbf{P}^3$. There is a self-map $\iota_p$ of the projective line $\mathbf{P}(T_pS)$ (indeed an involution by Proposition (4.6) below), which is defined as follows. Let $u \in \mathbf{P}(T_pS)$ be a tangent direction at $p$. The differential of the Gauss map sends $u$ to a direction $\gamma_S(u)$ in $\check{\mathbf{P}}^3$ tangent to $S^\vee$ at the point $(\mathbf{T}_pS)^\perp$: thus, $\gamma_S(u) \in \mathbf{P}\big(T_{(\mathbf{T}_pS)^\perp} S^\vee\big)$. The tangent direction $\gamma_S(u)$ spans a line $\langle \gamma_S(u) \rangle$ such that
$$(\mathbf{T}_pS)^\perp \subseteq \langle \gamma_S(u) \rangle \subseteq \mathbf{T}_{(\mathbf{T}_pS)^\perp} S^\vee.$$
By biduality, $(\mathbf{T}_{(\mathbf{T}_pS)^\perp} S^\vee)^\perp \cong \{p\}$. Therefore, the orthogonal of the line $\langle \gamma_S(u) \rangle$ is such that
$$p \in \langle \gamma_S(u) \rangle^\perp \subseteq \mathbf{T}_pS,$$
hence it defines a point $v \in \mathbf{P}(T_pS)$. We define $\iota_p$ as the map $u \mapsto v$.

If one prefers, $u$ is a point of $S$ infinitely near $p$. Then $\gamma_S(u)$ is a point of $S^\vee$ infinitely near the point $(\mathbf{T}_pS)^\perp$. Therefore, $\gamma_S(u)^\perp$ is a tangent plane of $S$ infinitely near the plane $\mathbf{T}_pS$. The intersection of these two planes is a line contained in $\mathbf{T}_pS$ and passing through $p$, and we let $\iota_p(u)$ be its tangent direction at $p$.

We will consider the map $\iota_p$ indifferently as a self-map of $\mathbf{P}(T_pS)$ or of the set of lines contained in $\mathbf{T}_pS$ and passing through $p$.

**(4.5) Lemma.** *Let $Y$ be a curve traced on a surface $S \subseteq \mathbf{P}^3$, and $p \in Y$ be a point which is non-parabolic on $S$. Let $\Sigma_Y$ be the developable surface defined as in (4.1). Then $p$ sits on $\Sigma_Y$, and the generator of $\Sigma_Y$ through $p$ is the line $\iota_p(\mathbf{T}_pY)$.*

*Proof.* By definition, $\Sigma_Y$ is the dual of the curve $\gamma_S(Y) \subseteq \check{\mathbf{P}}^3$. Therefore, by the general theory of developable surfaces (see Appendix B), the rulings of $\Sigma_Y$ are the orthogonals of the tangent lines of $\gamma_S(Y)$. The lemma then follows tautologically from the definition of $\iota_p$ in (4.4). □

**(4.6) Proposition.** *The map $\iota_p$ is the orthogonality involution of $\mathbf{P}(T_pS)$ defined by the second fundamental form at $p$. In particular, its fixed points are the tangent directions of the two flex tangents of $S$ at $p$.*

Since $p$ is assumed to be non-parabolic, the second fundamental form $II_p$ at $p$ is a non-degenerate quadratic form on the tangent plane $T_pS$ (up to scalar multiplication), see Appendix A.2. By definition, the orthogonality involution maps a vector line in the vector plane $T_pS$ to its orthogonal with respect to $II_p$, which is another vector line in $T_pS$.

*Proof.* We consider homogeneous coordinates $(x:y:z:w)$ in $\mathbf{P}^3$, such that $p = (0:0:0:1)$, and $\mathbf{T}_pS$ is the plane $z = 0$. We write the equation of $S$ in the form
$$f(x,y,z,w) = w^{n-1}z + w^{n-2}(ax^2 + 2rxy + by^2 + 2sxz + 2tyz + cz^2) + O(3) = 0.$$
It follows from the local computation in (A.13) that the second fundamental form at $p$ is given by the matrix

(4.6.1)
$$\begin{pmatrix} a & r \\ r & b \end{pmatrix}$$

in the coordinates $(x,y)$ on $T_pS$. On the other hand, the "large" Gauss map of $S$, defined on a Zariski open subset of $\mathbf{P}^3$, maps the point $(x:y:z:w)$ to the transpose of
$$\begin{pmatrix} \partial_x f \\ \partial_y f \\ \partial_z f \\ \partial_w f \end{pmatrix} = \begin{pmatrix} w^{n-2}(2ax + 2ry + 2sz) + \cdots \\ w^{n-2}(2rx + 2by + 2tz) + \cdots \\ w^{n-1} + w^{n-2}(2sx + 2ty + 2cz) + \cdots \\ (n-1)w^{n-2}z + \cdots \end{pmatrix}.$$

We write it in the affine chart $w = 1$; its differential at the origin (i.e., at $p$) is then

$$(4.6.2) \qquad \begin{pmatrix} 2a\,dx + 2r\,dy + 2s\,dz \\ 2r\,dx + 2b\,dy + 2t\,dz \\ 2s\,dx + 2t\,dy + 2c\,dz \\ (n-1)\,dz \end{pmatrix}.$$

Now, still in the affine chart $w = 1$, locally at $p$, $S$ is given by an implicit function $z = \varphi(x, y)$ with vanishing differential at the origin. In order to compute the differential of the "small" Gauss map of $S$, defined on a Zariski open subset of $S$, we need to compose the differential at the origin of $(x, y) \mapsto (x, y, \varphi(x, y))$ with that of the "large" Gauss map. This amounts to set $dz = 0$ in (4.6.2).

The upshot is that the tangent plane to $S$ at the point infinitely near $p$ in the direction $(u_0 : u_1)$ is given by the equation

$$(2au_0 + 2ru_1)x + (2ru_0 + 2bu_1)y + (2su_0 + 2tu_1)z = 0.$$

Therefore its intersection with $\mathbf{T}_p S$ is given by

$$(2au_0 + 2ru_1)x + (2ru_0 + 2bu_1)y = 0,$$

hence $\iota_p$ is the involution

$$(4.6.3) \qquad \begin{pmatrix} -r & -b \\ a & r \end{pmatrix} \in \mathrm{PGL}_2(\mathbf{C})$$

(it is indeed an involution because the trace of this matrix is 0, see [46, Chap. II, §D]). By a straightforward computation (in fact more or less the exact same computation that we have just performed), this is the orthogonality involution defined by the bilinear form given by the matrix (4.6.1), and this ends the proof. □

## 4.2 – The Hessian developable surface

In this section we study the family of planes in $\mathbf{P}^3$ tangent to $S$ at its parabolic points, and the associated developable surface in $\mathbf{P}^3$. We first need to adapt the observations of (4.4) and (4.6) to the case when $p$ is a parabolic point.

Let $p$ be a parabolic point of $S$. If $p$ is general, then the second fundamental form of $S$ at $p$ has rank 1, see Appendix A.2, hence the tangent section $\mathbf{T}_p S \cap S$ has a degenerate double point at $p$, and we call *double asymptotic (or flex) tangent line* of $S$ at $p$ the tangent cone $\mathrm{TC}_p(\mathbf{T}_p S \cap S)$.

**(4.7) Proposition.** *Let $p$ be a parabolic point of $S$ at which the second fundamental form does not vanish. Denote by $d\gamma_p$ the differential of the Gauss map of $S$ at $p$, and by $\Lambda$ the double asymptotic tangent line of $S$ at $p$. Then:*
*(i) the kernel of $d\gamma_p$ is the tangent space at $p$ of $\Lambda$;*
*(ii) the line in $\check{\mathbf{P}}^3$ passing through the point $\mathbf{T}_p S^\perp$ and going in the direction $\mathrm{im}(d\gamma_p)$ is the orthogonal of $\Lambda$.*

This tells, in particular, that if $p$ is a general parabolic point of $S$, the map $\iota_p$ constructed in (4.4) is essentially a contraction of $\mathbf{P}(T_p S)$ to its point corresponding to the double flex tangent line of $S$ at $p$.

*Proof.* This follows from the computations carried out in the proof of Proposition (4.6). Note that, in the notation introduced there, we may assume $a = r = 0$ since the second fundamental form at $p$ has rank 1. Then, we see from (4.6.2) that the kernel of $d\gamma$ is $y = 0$, which is the reduced tangent cone of the tangent section $\mathbf{T}_p S \cap S$ at $p$, i.e., the tangent space at $p$ of $\Lambda$. In addition, we see from (4.6.3) that for $u = (u_0 : u_1)$, the orthogonal of $\langle d\gamma(u) \rangle$ is $y = 0$ if $u_1 \neq 0$ (and the whole $\mathbf{T}_p S$ is $u_1 = 0$). This proves (ii). □

**(4.8) Remark.** Proposition (4.7) is a biduality statement, the analog of (2.2.2) for parabolic points.

Let $p$ be a general parabolic point. The line in $\check{\mathbf{P}}^3$ passing through the point $\gamma(p) = \mathbf{T}_p S^\perp$ and going in the direction $\mathrm{im}(d\gamma_p)$ is the tangent line $\mathbf{T}_{\mathbf{T}_p S^\perp} \check{C}$ to the cuspidal double curve of $S^\vee$ at $p$, denoted by $\check{C}$. Therefore, the proposition tells us that the orthogonal of the double flex tangent of $S$ at $p$ (denoted by $\Lambda$) is the line $\mathbf{T}_{\mathbf{T}_p S^\perp} \check{C}$. We give another proof of this in Appendix C.

Note that the double flex tangent $\Lambda$ is the line of maximum contact with $S$ at $p$, and intersects $S$ with multiplicity 3, and the tangent to the cuspidal edge, $\mathbf{T}_{\mathbf{T}_p S^\perp} \check{C}$, is the line of maximum contact with $S^\vee$ at the point $\mathbf{T}_p S^\perp$, and intersects $S^\vee$ with multiplicity 4 at this point (the intersection multiplicity is 2 for a general line, 3 for a general line in the tangent cone, and 4 for the tangent of the cuspidal edge).

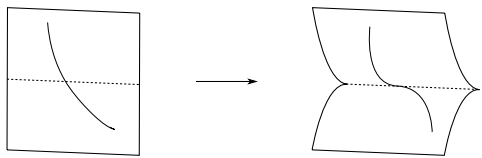

Figure 1: The Gauss map around the parabolic curve

**(4.9) Proposition.** *Let $p$ be a parabolic point of $S$ at which the second fundamental form does not vanish. The tangent line at $p$ to the curve $S \cap \mathrm{Hess}(S)$ coincides with the double flex tangent line of $S$ at $p$ if and only if the tangent section $S \cap \mathbf{T}_p S$ has a tacnode at $p$ (or a worse singularity).*

*Proof.* Let us write the equation of $S$ in the form
$$f(x,y,z,w) = w^{n-1}z + w^{n-2}\left(\tfrac{1}{2}y^2 + axz + byz + \tfrac{1}{2}cz^2\right) + w^{n-3}f_3(x,y,z) + O(4) = 0,$$
in a system of homogeneous coordinates $(x : y : z : w)$ such that $p = (0 : 0 : 0 : 1)$ ($n$ is the degree of $S$, and $f_3$ is homogeneous of degree 3). The tangent at $p$ to the Hessian curve is obtained by intersecting $\mathbf{T}_p S$, that is, the plane $z = 0$, with the tangent plane to the Hessian, which we will find with formula (A.11.1).

In the affine chart $w = 1$, we have
$$\partial_x f = az + \partial_x f_3 + O(3)$$
$$\partial_y f = y + bz + \partial_y f_3 + O(3),$$
$$\partial_z f = 1 + ax + by + cz + \partial_z f_3 + O(3)$$

so that, up to a multiplicative constant, $\mathrm{Hess}(f)$ is

$$\begin{vmatrix} \tfrac{n}{n-1}(z + \tfrac{1}{2}y^2 + \cdots) & az + \partial_x f_3 + \cdots & y + bz + \partial_y f_3 + \cdots & 1 + ax + by + cz + \partial_z f_3 + \cdots \\ az + \partial_x f_3 + \cdots & \partial_x^2 f_3 + \cdots & \partial_x \partial_y f_3 + \cdots & a + \partial_x \partial_z f_3 + \cdots \\ y + bz + \partial_y f_3 + \cdots & \partial_x \partial_y f_3 + \cdots & 1 + \partial_y^2 f_3 + \cdots & b + \partial_y \partial_z f_3 + \cdots \\ 1 + ax + by + cz + \partial_z f_3 + \cdots & a + \partial_x \partial_z f_3 + \cdots & b + \partial_y \partial_z f_3 + \cdots & c + \partial_z^2 f_3 + \cdots \end{vmatrix}.$$

Plugging in $z = 0$ and keeping for each entry only the terms of order at most 1, we obtain

$$\begin{vmatrix} 0 & 0 & y & 1+ax+by \\ 0 & \partial_x^2 f_3 & \partial_x\partial_y f_3 & a+\partial_x\partial_z f_3 \\ y & \partial_x\partial_y f_3 & 1+\partial_y^2 f_3 & b+\partial_y\partial_z f_3 \\ 1+ax+by & a+\partial_x\partial_z f_3 & b+\partial_y\partial_z f_3 & c+\partial_z^2 f_3 \end{vmatrix}.$$

The only linear term in this determinant is $\partial_x^2 f_3$, so the tangent line to $S \cap \text{Hess}(S)$ has equations

(4.9.1) $$z = \partial_x^2 f_3 = 0.$$

Therefore, it coincides with the double flex tangent line if and only if there is no $x^3$ term in $f_3$, which is equivalent to the tangent section of $S$ at $p$ having an equation of the form

$$\frac{1}{2} y^2 + y(a_{21} x^2 + a_{12} xy + a_{03} y^2) + O(4) = 0,$$

which is the local form of a tacnodal singularity. $\square$

**(4.10) Remark.** Let $p \in S$ be a point such that the tangent section $\mathbf{T}_p S \cap S$ has a tacnode at $p$ (in particular $p$ lies on the Hessian curve $\text{He}(S) = S \cap \text{Hess}(S)$). It follows from Proposition (4.9) that the image of the Hessian curve $\text{He}(S)$ by the Gauss map (i.e., the cuspidal double curve of $S^\vee$) has a cusp at $\gamma_S(p)$. (We had seen this fact already in (3.1) and (3.7.5), arguing that $\gamma_S(p)$ is a swallowtail of $S^\vee$).

Indeed, if the second fundamental form has rank one at $p$, then the kernel of the differential of the Gauss map at $p$ is the tangent space of the double flex tangent of $S$ at $p$ by Proposition (4.7), and this coincides with the tangent space of $S \cap \text{Hess}(S)$ by Proposition (4.9). On the other hand, if the second fundamental form at $p$ is 0, there is nothing to prove.

**(4.11) Proposition.** *Let $p \in S$ be a point such that the tangent section $\mathbf{T}_p S \cap S$ has a tacnode at $p$. The regression (i.e., cuspidal) edge of the Hessian developable $\Sigma_{\text{He}}$ passes through $p$, and its osculating plane there is $\mathbf{T}_p S$. In particular the two surfaces $\Sigma_{\text{He}}$ and $S$ are tangent at $p$.*

*Proof.* This follows from the fact, which we have observed in (3.7.5), that $p^\perp$ is the osculating plane of the cuspidal double curve of $S^\vee$, which we will denote here by $C \subseteq \check{\mathbf{P}}^3$. Note that it has a local parametrization of the form $t \mapsto (t^2, t^3, t^4)$ at $p$, and has intersection multiplicity 4 with $p^\perp$ at the point $(\mathbf{T}_p S)^\perp$, see (3.1) and (3.7.5).

Let $\Gamma \subseteq \mathbf{P}^3$ be the regression edge of $\Sigma_{\text{He}}$. By the general theory of developable surfaces, see Appendix B, $\Gamma$ parametrizes the osculating planes of $C$. Thus, since $p^\perp$ is an osculating plane of $C$, $p = (p^\perp)^\perp$ sits on $\Gamma$. Moreover, by the osc-duality for space curves, see again Appendix B, the osculating plane of $\Gamma$ at $p$ is the orthogonal of the point at which $p^\perp$ is osculating to $C$: the latter point is $\gamma_S(p) = (\mathbf{T}_p S)^\perp$, hence the osculating plane of $\Gamma$ at $p$ is $\mathbf{T}_p S$. Finally, the tangent plane to $\Sigma_{\text{He}}$ at $p$ is the osculating plane of $\Gamma$ at $p$, which we have just seen is $\mathbf{T}_p S$. $\square$

We may now give a direct geometric description of the Hessian developable.

**(4.12) Theorem.** *The generators of the Hessian developable surface are the double flex tangent lines at parabolic points of $S$. Its stationary planes are the tangent planes to $S$ at points where the tangent section has a tacnode.*

*Proof.* The identification of the generators follows from Proposition (4.7) and Lemma (4.5) (which works fine even at parabolic points). That of the stationary planes follows from Proposition (4.9), see also Remark (4.10). $\square$

The rest of this section is dedicated to the computation of the numerical characters of the Hessian developable.

**(4.13) Lemma.** *Let $S, T$ be two surfaces in $\mathbf{P}^3$ of degrees $s, t$ respectively. The degree of the ruled surface generated by the flex tangents of $S$ at the points of $S \cap T$ equals $st(3s - 4)$.*

Note that there is no reason that this ruled surface be developable in general, although it is indeed when $T = \mathrm{Hess}(S)$.

*Proof.* Let $f, g$ be homogeneous equations of $S$ and $T$ respectively. The ruled surface we are interested in is the projection on the first factor of the incidence graph

$$\{(p, p') : p' \in S \cap T \text{ and } p \text{ sits on a flex tangent of } S \text{ at } p'\},$$

which is the complete intersection in $\mathbf{P}^3 \times \mathbf{P}^3$ defined by the following four bihomogeneous polynomials in the variables $(p, p')$: $f(p')$, $g(p')$, $\mathrm{D}^p f(p')$, $\mathrm{D}^{p^2} f(p')$, of bidegrees $(0, s)$, $(0, t)$, $(1, s-1)$, and $(2, s-2)$ respectively.

Our ruled surface is therefore the hypersurface in $\mathbf{P}^3$ defined by the homogeneous equation obtained by eliminating the variable $p'$, i.e., by the polynomial

$$\mathrm{Res}_{s,t,s-1,s-2}(f, g, \mathrm{D}^p f, \mathrm{D}^{p^2} f),$$

taken as a resultant of homogeneous polynomials in the variable $p'$, with coefficients homogeneous polynomials in the variable $p$. By a standard homogeneity property of the resultant, see [C, Theorem **], it is a homogeneous polynomial of degree

$$\bigl[0 \times t(s-1)(s-2)\bigr] + \bigl[0 \times s(s-1)(s-2)\bigr] + \bigl[1 \times st(s-2)\bigr] + \bigl[2 \times st(s-1)\bigr] = st(3s-4)$$

in the variable $p$. □

**(4.14) Corollary.** *The degree of the Hessian developable surface equals $2n(n-2)(3n-4)$.*

*Proof.* We apply the above lemma with $T = \mathrm{Hess}(S)$. This gives the degree $4n(n-2)(3n-4)$, which however has to be divided by two, as the obtained equation is double in this case, being the two flex tangents coincident at parabolic points. □

**(4.15)** The class of the Hessian developable is the number of planes tangent to $S$ at a parabolic point, and passing through a fixed general point in $\mathbf{P}^3$. It is therefore the intersection number of $S$, $\mathrm{Hess}(S)$, and $\mathrm{D}^{p'} S$ for a general $p' \in \mathbf{P}^3$, which equals $4n(n-1)(n-2)$.

**(4.16)** The number of stationary planes of the Hessian developable is, by Theorem (4.12), the number of points in the intersection of the parabolic curve $S \cap \mathrm{Hess}(S)$ with the flecnodal locus. We have seen in (3.15) that this number is $2n(n-2)(11n-24)$.

**(4.17)** We may now find the remaining characteristic numbers using the formulae of Appendix B. We use the notation of (B.6) and (B.17), with an additional decoration 'He' to avoid confusion; see also the warning below. We have found so far:
— $r_{\mathrm{He}} = 2n(n-2)(3n-4)$, the degree (or rank) of the developable, see Corollary (4.14);
— $n_{\mathrm{He}} = 4n(n-1)(n-2)$, its class, see (4.15);
— $\alpha_{\mathrm{He}} = 2n(n-2)(11n-24)$, its number of stationary planes, see (4.16).
Knowing three quantities is enough to derive all the others by direct applications of the formulae in (B.19) and (B.20). This way one finds, as in [45, §608]:
— $m_{\mathrm{He}} = 4n(n-2)(7n-15)$, the degree of the regression edge of $\Sigma_{\mathrm{He}}$;

— $\beta_{\mathrm{He}} = 10n(n-2)(7n-16)$, the number of flexes of the cuspidal double curve $\check{C}$ of $S^\vee$;
— $g_{\mathrm{He}} = 2n(n-2)(4n^4 - 16n^3 + 20n^2 - 27n + 39)$, the number of apparent double points of $\check{C} \subseteq S^\vee$;
— $h_{\mathrm{He}} = 2n(n-2)(196n^4 - 1232n^3 + 2580n^2 - 1861n + 137)$, the number of apparent double points of the regression edge of $\Sigma_{\mathrm{He}}$.

*(4.17.1) Warning.* For ease of reference I have done my best to stick to Salmon's notation, as other have done before me, e.g., Piene [38], Ronga [42]. There are however some compatibility problems between the notation of (B.6) and (B.17) for developable surfaces, and that of (3.4), (3.6), and (4.20) for the singularities of the dual. For instance, $\alpha_{\mathrm{He}} = \check{\beta}$, both the number of stationary planes in the Hessian system and of swallowtails of $S^\vee$, and $g_{\mathrm{He}} = \check{h}$, the number of apparent double points of the cuspidal double curves of $S^\vee$.

*(4.17.2) References.* Some of these numbers have been computed by Kulikov in an interesting series of two papers, with intersection theoretic computations. In particular, the number $\alpha_{\mathrm{He}}$ is computed in [28] ("number of points $P_3$"), and the number $\beta_{\mathrm{He}}$ is computed ("number of points type $P_3$") in [29]; there the regression edge of the Hessian developable is called the focal curve, and "points of type $P_3$" are the cusps of this curve, hence correspond to the flexes of $C = \Sigma_{\mathrm{He}}^\vee$.

Moreover, [28] provides some invariants of the flecnodal curve that we have not considered here. "Points $H_3$" are the ordinary double points of the flecnodal curve; they are the points $q \in S$ such that both local branches at $q$ of the tangent section $\mathbf{T}_q S \cap S$ have a flex at $q$. "Points $H_4$" are the points where the flecnodal curve is tangent to an asymptotic tangent line of $S$; they are the points $q \in S$ such that the tangent section $\mathbf{T}_q S \cap S$ has a branch at $q$ with a higher order inflection point. Kulikov obtains:
— $\#H_3 = 5n(7n^2 - 28n + 30)$;
— $\#H_4 = 5n(n-4)(7n-12)$.
These numbers are given by Salmon as well, at [45, §604].

## 4.3 – Ordinary double edges of the polar cone

In this section we consider the same setup as in Section 3.2, and use the notation introduced there. In particular, we let $X$ be a degree $n$ surface in $\mathbf{P}^3$ with d-ordinary singularities, and call $B$ and $C$ its ordinary and cuspidal double curves. We consider a general point $p \in \mathbf{P}^3$, and denote by $\mathcal{D}^p X$ the polar cone of $X$ with respect to $p$, and by $\mathcal{B}^p X$ the apparent boundary.

**(4.18) Ordinary double edges of the polar cone.** In this section we study the ordinary double edges of the polar cone $\mathcal{D}^p X$; they are the lines passing through $p$ and intersecting $X$ with multiplicity 2 in two distinct points. The locus of points $q \in S$ such that the line $\langle p, q \rangle$ verifies this condition is cut out on $X \cap \mathrm{D}^p X$ by a degree $(n-2)(n-3)$ hypersurface $R^p X$: this has been proved by Salmon, and his construction is explained in [C, Section **]. This hypersurface has been invoked in (2.8.3) already.

We shall conduct our study of ordinary double edges by analyzing the various kinds of intersection points of this hypersurface $R^p X$ with the curve $X \cap \mathrm{D}^p X$. This strategy is similar to that followed in Section 3.2, with the hypersurface $R^p X$ playing the role that the second polar $\mathrm{D}^{p^2} X$ played there.

**(4.19) Definition.** *The* number of apparent double points *of a curve in $\mathbf{P}^3$ is the number of double points of its projection from a general point.*

*The* number of apparent only double points *of two curves $Y, Y'$ in $\mathbf{P}^3$ is the number of intersection points of their projections from a general point that do not come from an actual intersection point of $Y$ and $Y'$. We denote it by $[Y, Y']^\circ$.*

**(4.20) Notation.** We add the following notation to that introduced in Section 3.2:
— $\delta$ is the number of ordinary double edges of the circumscribed cone $\mathcal{C}_p X$ (these are the honest bitangent lines to $X$);
— $k$ is the number of apparent double points of the ordinary double curve $B$ of $X$;
— $h$ is the number of apparent double points of the cuspidal double curve $C$ of $X$;
— $i$ is the number of intersection points of $B$ and $C$ smooth for both curves (when $X = S^\vee$ for a smooth $S$, this number vanishes).

The following formulae, given by Salmon without any comment, describe the intersection of the hypersurface $R^p X$ mentioned in (4.18) above with the various components of the curve

$$X \cap \mathrm{D}^p X = \mathcal{B}^p X + 2B + 3C,$$

of total degree $n(n-1)(n-2)(n-3)$. Recall that $a, b, c$ denote the degrees of the curves $\mathcal{B}^p X, B, C$ (called $A, B, C$ by Salmon), respectively.

**(4.21) Lemma.** *The following relations hold:*

(☙i)   $\qquad a(n-2)(n-3) = 2\delta + 2[\mathcal{B}^p X \cdot B]^\circ + 3[\mathcal{B}^p X \cdot C]^\circ$ ;
(☙ii)  $\qquad b(n-2)(n-3) = 4k + [\mathcal{B}^p X \cdot B]^\circ + 3[B \cdot C]^\circ$ ;
(☙iii) $\qquad c(n-2)(n-3) = 6h + [\mathcal{B}^p X \cdot C]^\circ + 2[B \cdot C]^\circ$.

*Proof.* Let $Y, Y'$ be two curves (possibly the same) among $\mathcal{B}^p X, B, C$, and denote by $m_Y, m_{Y'}$ their multiplicities in the cycle $X \cap \mathrm{D}^p X$. An intersection point between two local branches of the projections $\pi_p(Y)$ and $\pi_p(Y')$ gives a double edge of the polar cone with multiplicity $m_Y m_{Y'}$, hence two contact points $y \in Y$ and $y' \in Y'$ with $X$, each of which should contribute with the same multiplicity $m_Y m_{Y'}$ to the intersection $X \cap \mathrm{D}^p X \cap R^p X$. Taking into account the multiplicities of $Y$ and $Y'$ in $X \cap \mathrm{D}^p X$, one finds that $R^p X$ should intersect $Y$ with multiplicity $m_{Y'}$ at $y$ and $Y'$ with multiplicity $m_Y$ at $y'$. This explains all the coefficients in the formulae. □

Some claims in the above proof arguably deserve a more rigorous treatment, but let me not try to fix this here. Note that the actual intersection points of $\mathcal{B}^p X, B, C$ do not contribute to the intersection $X \cap \mathrm{D}^p X \cap R^p X$, as they give rise to cuspidal double edges of the polar cone instead of ordinary double edges.

For two curves $Y$ and $Y'$ in $\mathbf{P}^3$, the total number of intersection points of the projections $\pi_p(Y)$ and $\pi_p(Y')$ is $\deg(Y) \deg(Y')$. The following set of formulae, yet again given coldbloodedly by Salmon, indicates the quantities to be subtracted from $\deg(Y) \deg(Y')$ to obtain $[Y, Y']^\circ$, in the case when $Y$ and $Y'$ are two components of $X \cap \mathrm{D}^p X$.

**(4.22) Theorem.** *The following relations hold:*

(✠i)   $\qquad [\mathcal{B}^p X \cdot B]^\circ = ab - 2\rho$ ;
(✠ii)  $\qquad [\mathcal{B}^p X \cdot C]^\circ = ac - 3\sigma$ ;
(✠iii) $\qquad [B \cdot C]^\circ = bc - 3\beta - 2\gamma - i$.

*Proof.* In our setup, every actual intersection point of two curves among $\mathcal{B}^p X, B, C$ is of one of the types $\rho, \sigma, \beta, \gamma, i$, so only the coefficients are to be explained. The latter are the intersection multiplicities at the corresponding points of the projections from $p$ of the two curves under consideration. For points of type $i$, it is 1 by the generality of $p$. The rest of the proof consists in examining one by one the other possible types.

*(4.22.1) Points of type $\rho$.* These are the points $q \in \mathcal{B}^p X \cap B$. At such a point, the surface $X$ has two smooth local sheets $X'_q$ and $X''_q$, $B = X'_q \cap X''_q$ and the circumscribed curve $\mathcal{B}^p X$ is contained in one of them, say $X'_q$, and transverse to $B$. The curve $\mathcal{B}^p X$ is the ramification curve of the projection of $X'_q$, hence the projections of $\mathcal{B}^p X$ and $B$ are tangent at the projection of $q$, and thus intersect with multiplicity 2.

*(4.22.2) Points of type $\sigma$.* These are the points $q \in \mathcal{B}^p X \cap C$. At such a point, the surface $X$ has one cuspidal local sheet $X_q$, and $C = \mathrm{Sing}(X_q)$. We have seen in (3.7.3) that $\mathcal{B}^p X$ and $C$ are smooth and tangent at $q$, and it follows from the same argument as in (4.22.1) that their projections from $p$ are tangent. But the claim is that, in fact, the two projections intersect with multiplicity 3.

It should be possible to prove this with a local computation, using the setup of (3.7.3), but I have not been able to do it in general. The difficulty is that, while the curve defined by the two equations $f$ and $\mathrm{D}^p f$ splits as $\mathcal{B}^p X + 3C$ locally at $q$ (in the notation of (3.7.3)), neither of the two summands is a Cartier divisor on $X$, and therefore it is difficult to find their equations. It is however doable in practice for an arbitrary, explicit, choice of the polynomials $g_2$ and $h_2$, using one's favourite computer algebra system; I have done it for several randomly chosen polynomials, using Macaulay2; I recommend to do this over a finite field. Finally, once one knows equations for $C$ and $\mathcal{B}^p X$, one finds equations for their projections by eliminating the variable $x$ (because $p = (1:0:0:0)$), and it only remains to check with one's own eyes that the intersection multiplicity at the projection of $q$ is indeed 3.

In any event, the number $[\mathcal{B}^p X \cdot C]^\circ$ occurs only in (☙ ㅣ) and (☙ ㅣㅣㅣ), in which we already know all the other numbers involved. One may therefore use them to find $[\mathcal{B}^p X \cdot C]^\circ$, and thus obtain an indirect proof of Formula (☙).

*(4.22.3) Points of type $\beta$.* These are the swallowtail singularities of $X$. At these points $B$ is smooth and $C$ is cuspidal, and we can see on the local model in (3.1) that they have the same reduced tangent cone. It follows that after projection they intersect with multiplicity 3 at these points.

*(4.22.4) Points of type $\gamma$.* These are the points at which $X$ has two transverse local sheets, one smooth: $X^s_q$, and the other cuspidal: $X^c_q$. There $C$ is smooth and $B$ is cuspidal, and these two curves have transverse reduced tangent cones. It follows that after projection they intersect with multiplicity 2 at these points.

Theorem (4.22) is now proved. $\square$

Finally, putting together formulas (☙) and (☙) gives the following.

**(4.23) Corollary.**

$$\begin{aligned}
(\text{☙ ㅣ}) \qquad & a(n-2)(n-3) = 2\delta + 2ab + 3ac - 4\rho - 9\sigma \;; \\
(\text{☙ ㅣㅣ}) \qquad & b(n-2)(n-3) = 4k + ab + 3bc - 9\beta - 6\gamma - 3i - 2\rho \;; \\
(\text{☙ ㅣㅣㅣ}) \qquad & c(n-2)(n-3) = 6h + ac + 2bc - 6\beta - 4\gamma - 2i - 3\sigma.
\end{aligned}$$

## 4.4 – Application to the case when $X = S^\vee$

We now proceed as in Section 3.3, i.e., we let $X$ be the dual of a smooth, sufficiently general, degree $n$ surface $S \subseteq \mathbf{P}^3$, and apply the formulas of Section 4.3 above to $X$. We maintain the setup and notation introduced in Section 3.3.

**(4.24)** As in (3.11), the circumscribed cone $\mathcal{C}^{\check{p}}(S^\vee)$ is the dual in $\check{\mathbf{P}}^3$ of the plane curve $\check{p}^\perp \cap S$. The number $\check{\delta}$ of ordinary edges of this cone is therefore the number of bitangents of the plane curve $\check{p}^\perp \cap S$, which by (1.14) is:
— $\check{\delta} = \frac{1}{2}n(n-2)(n-3)(n+3)$.

**(4.25)** The numbers $\check{h}$ and $\check{k}$ are the invariants $g$ (introduced in (B.17)) of the Hessian and node-couple developables respectively. The former has been computed in (4.17) and the latter we don't know yet:
— $\check{h} = g_{\text{He}} = 2n(n-2)(4n^4 - 16n^3 + 20n^2 - 27n + 39)$ ;
— $\check{k} = g_{\text{Kn}}$.

**(4.26) Examination of formulae ($\underset{\approx}{\text{ǔǔǔ}}$ ◊◊).** The numbers $\check{n}, \check{a}, \check{b}, \check{c}$ have been given in Section 3.3, and the numbers $\check{\rho}, \check{\sigma}, \check{\beta}, \check{\gamma}$ have already been computed. As we mentioned earlier, the number $\check{\iota}$ vanishes. One may thus check by a direct computation that the first and third identities of (4.23) indeed hold. The second one, on the other hand, gives the formula:
— $\check{k} = \frac{1}{8}n(n-2)(n^{10}-6n^9+16n^8-54n^7+164n^6-288n^5+547n^4-1058n^3+1068n^2-1214n+1464)$.

**(4.27) Numerical characters of the node-couple developable.** At this point, we know the following three numerical characters of the node-couple developable: $n_{\text{Kn}} = \check{b}$, the degree of the curve $\check{B} \subseteq \check{\mathbf{P}}^3$, $g_{\text{Kn}} = \check{k}$, the number of apparent double points of $\check{B}$, and $\alpha_{\text{Kn}} = \check{\gamma}$, the number of cusps of $\check{B}$. We may therefore derive all the others using the formulas of Section B.3; however, these ought to be suitably modified in order to take into account the presence of $\check{t}$ ordinary triple points on the curve $C_{\text{Kn}} = \check{B}$. This may be taken care of in various ways, but the idea is always the same, namely that an ordinary triple point counts exactly as three nodes.

For instance, let us adapt the formula $r = n(n-1) - 2g - 3\alpha$ from (B.19) in order to compute the rank $r_{\text{Kn}}$ of the node-couple developable. The adapted formula is

$$r_{\text{Kn}} = n_{\text{Kn}}(n_{\text{Kn}} - 1) - 2\big(g_{\text{Kn}} + 3\#(\text{triple points})\big) - 3\alpha_{\text{Kn}}$$
$$= \check{b}(\check{b} - 1) - 2\check{k} - 6\check{t} - 3\check{\gamma},$$

and it gives
— $r_{\text{Kn}} = n(n-2)(n-3)(n^2 + 2n - 4)$.
We leave this task for the remaining numerical characters an exercise to the not-exhausted-yet reader.

Finally, I would like to mention the article [11], in which the authors perform an interesting study of the family of bitangent lines of a very general quartic surface, which involves in particular the node-couple developable.

# 5 – Formulae for surfaces projected in $\mathbf{P}^3$

In this section, we prove formulae for the enumeration of curves with prescribed singularities in a possibly incomplete linear system of dimension $r \leqslant 3$ on a smooth surface; this will take the form of enumerating curves on a surface in $\mathbf{P}^3$ with ordinary singularities, see (5.3). These formulae are from a different age with respect to those of Salmon, and involve invariants of a different nature, e.g., the arithmetic genus; by the way, beware the difference in notation with respect to the previous sections. The material presented here is taken from [15, Capitolo V].

**(5.1) Setup.** We consider a linear system $|V|$ of dimension $r \leqslant 3$ on a smooth surface $S$. We will usually make some transversality assumptions on $|V|$, which are typically granted if $|V|$ is a general subsystem of a very ample linear system. In particular, we always assume that the general member of $|V|$ is smooth and irreducible. We set:
— $C$ the linear equivalence class of divisors in $|V|$ (and sometimes a member of $|V|$ by abuse of notation) ;
— $n$ the degree of $|V|$, i.e., $n = C^2$ ;
— $\pi$ the common arithmetic genus of all members of $|V|$.
We have the following first identity, given by the adjunction formula:
$$(5.1.1) \qquad K_S \cdot C = 2\pi - 2 - n.$$

**(5.2) Warm up: the $1$-dimensional case.** Let $r = 1$, and assume that $|V|$ has $n$ simple base points, and all its singular members are 1-nodal curves. In this context, we define the class $\check{n}$ of $|V|$ as the number of its singular members. In this case, the topological formula discussed in details in [XI, Section **] gives:
— $\check{n} = e(S) + n - e(C) \cdot e(\mathbf{P}^1) = c_2(S) + n + 4\pi - 4$.

Next we treat the case $r = 3$. For $r = 2$ we shall only consider a cheap version in which we assume that $|V|$ is a subsystem of a linear system of dimension at least 3.

## 5.1 – The $3$-dimensional case: surfaces with ordinary singularities

In this section we let $r = 3$, and make the following assumption, granted by the general projection theorem, see [20], if $|V|$ is a general subssytem of a very ample linear system.

**(5.3)** We assume that the map defined by $|V|$ is a birational morphism from $S$ to a degree $n$ surface $S^\flat \subseteq \mathbf{P}^3$ with *ordinary singularities* (in the usual sense), i.e., $S^\flat$ is smooth except for an ordinary double curve $D$ with $\tau$ pinch points and $t$ triple points: locally at a pinch point $S^\flat$ has a Whitney umbrella singularity, i.e., it looks like the hypersurface $x^2 - zy^2 = 0$ in the affine space with coordinates $x, y, z$; locally at a triple point it looks like the hypersurface $xyz = 0$, so these points are non-planar triple points for $D$; locally at all other points of $D$, $S^\flat$ looks like the hypersurface $xy = 0$; in particular, $D$ is smooth except for the triple points. For simplicity, we assume moreover that $D$ is irreducible.

The morphism $\nu : S \to S^\flat$ is the normalization map. We let $\Delta \subseteq S$ be the pre-image of $D \subseteq S^\flat$; there is a double cover $\Delta \to D$, which ramifies exactly over the $\tau$ pinch points. By adjunction theory (see, e.g., [22, Exer. II.8.5] and [II, Section **]), one has $\omega_S \cong \mathcal{O}_S(-\Delta) \otimes \nu^* \omega_{S^\flat}$, and $\nu_* \omega_S = \mathcal{I}_D|_{S^\flat}$, where $\mathcal{I}_D$ is the sheaf of ideals defining $D$ in $\mathbf{P}^3$. Being $S^\flat$ a divisor in $\mathbf{P}^3$, we have $\omega_{S^\flat} \cong \omega_{\mathbf{P}^3}(S^\flat)\big|_{S^\flat} = \mathcal{O}_{S^\flat}(n-4)$ by the standard adjunction formula.

**(5.4)** In addition to the already introduced notation (see also (5.1)), we let:
— $\tau$ be the number of pinch points of $S^\flat$ ;
— $t$ be the number of triple points of $S^\flat$ ;
— $d$ be the degree of the double curve $D \subseteq \mathbf{P}^3$ ;
— $\rho$ be the geometric genus of $D$ ;
— $\Pi$ be the geometric genus of $\Delta \subseteq S$.

**(5.5)** The curves $C$ are represented in $\mathbf{P}^3$ as hyperplane sections of $S^\flat$, hence as degree $n$ plane curves with $d$ nodes. Therefore one has
$$(5.5.1) \qquad \pi = \binom{n-1}{2} - d.$$

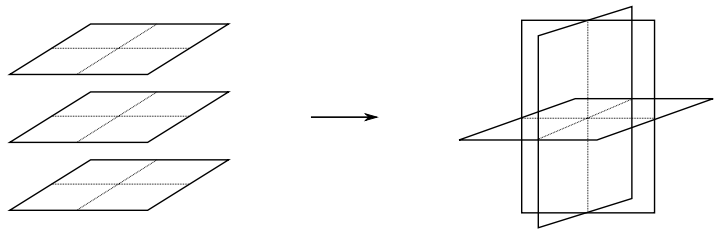

Figure 2: Normalization of $S^\flat$ at a triple point

**(5.6) Effective degree of first polars.** We compute the class $\check{n}$ of $S^\flat$, that is the number of singular curves in a general subpencil of $|V|$, following the method explained in (1.8).

For a general point $p \in \mathbf{P}^3$, the first polar $\mathrm{D}^p S^\flat$ is a degree $n-1$ hypersurface passing simply through $D$, and we define the apparent boundary $\mathcal{B}^p S^\flat$ as the Zariski closure of $(\mathrm{D}^p S^\flat \cap S^\flat) - D$, as in (3.5). The curve $\mathcal{B}^p S^\flat$ passes through all pinch points: this follows from a local computation, which we reckon the reader who has survived this far can make on his own.

Therefore, for a general pair of points $p, q \in \mathbf{P}^3$ the intersection of the two curves $\mathcal{B}^p S^\flat$ and $\mathcal{B}^q S^\flat$ consists of the $\tau$ pinch points of $S^\flat$ plus the $\check{n}$ points $a$ such that $p, q \in \mathbf{T}_a S$. To compute this intersection number $\tau + \check{n}$, we pull-back the curves to the normalization $S$. They are sections of
$$\mathcal{I}_D \otimes \mathcal{O}_{S^\flat}(n-1) \cong \mathcal{I}_D \otimes \omega_{S^\flat} \otimes \mathcal{O}_{S^\flat}(3),$$
hence pull-back to curves in the class $K_S + 3C$. We thus find the identity

(5.6.1) $$\check{n} + \tau = (K_S + 3C)^2 = K_S^2 + 3n + 12\pi - 12,$$

using (5.1.1) to substitute the intersection number $K_S \cdot C$.

**(5.7) Genus of the curve of neutral pairs of $|V|$.** One first observes that the curve $\Delta$ has three nodes above each triple point of $S^\flat$, as is best understood looking at Figure 2. Since it is otherwise smooth, its arithmetic genus is $p_a(\Delta) = p_g(\Delta) + 3t = \Pi + 3t$. On the other hand, using the already seen identity $K_S = (n-4)C - \Delta$, the adjunction formula gives
$$2(\Pi + 3t) - 2 = 2p_a(\Delta) - 2 = (K_S + \Delta) \cdot \Delta = (n-4)C \cdot \Delta = 2(n-4)d.$$
Therefore,

(5.7.1) $$p_a(\Delta) = \Pi + 3t = (n-4)d + 1.$$

On the other hand, the map $\Delta \to D$ lifts to a $2:1$ covering $\bar{\Delta} \to \bar{D}$ between the normalizations of $\Delta$ and $D$. It is ramified at $\tau$ points, so the Riemann–Hurwitz formula gives

(5.7.2) $$2\Pi - 2 = 2(2\rho - 2) + \tau.$$

**(5.8) Intersections of apparent boundaries with the double curve.** Let $p \in \mathbf{P}^3$ be a general point. We consider the apparent boundary $\mathcal{B}^p S^\flat$ and its intersection with the double curve $D$. To see the situation clearly, it is best to work on the normalization $S$. One has
$$\nu^* \mathcal{B}^p S^\flat \cdot \Delta = (K_S + 3C) \cdot \big((n-4)C - K_S\big) = -K_S^2 + 2\pi(n-7) + 2n^2 - 7n + 14.$$

Geometrically this intersection amounts to the $\tau$ pinch points plus the points $q \in D$ such that the tangent plane to one of the two local sheets of $S^\flat$ at $q$ passes through $p$. Let $a$ be the number of the latter points. It may be computed by considering the intersection of the second polar $D^{p^2}S^{\flat}$ with $D$. It follows from Theorem (A.8) that this intersection consists of the $a$ points we are interested in plus the $t$ triples points, and since those are triple for $D$ one finds

$$a + 3t = (n-2)d.$$

Finally, we end up with

(5.8.1) $$\tau + [(n-2)d - 3t] = -K_S^2 + 2\pi(n-7) + 2n^2 - 7n + 14.^4$$

**(5.9) Postulation formula.** We shall also use the following formula for the arithmetic genus of $S$, which we shall discuss in Section 5.2:

(5.9.1) $$p_a(S) = \binom{n-1}{3} - (n-4)d + \rho + 2t - 1.$$

**(5.10) Conclusion.** We may now express the quantities $\check{n}$, $\tau$, $t$, $d$, $\rho$, $\Pi$ in terms of the four invariants $n$, $\pi$, $p_a(S)$, $K_S^2$. The formula for $d$ is readily given by (5.5.1); we use it so substitute for $d$ in the five equations (5.6.1), (5.7.1), (5.7.2), (5.8.1), (5.9.1), thus obtaining

$$\begin{pmatrix} 1 & & & & 1 \\ & 2 & 6 & & \\ & 2 & 0 & -4 & -1 \\ & -3 & 0 & & 1 \\ & -2 & -1 & & 0 \end{pmatrix} \begin{pmatrix} \check{n} \\ \Pi \\ t \\ \rho \\ \tau \end{pmatrix} = \begin{pmatrix} 3n + 12\pi + K_S^2 - 12 \\ (n-4)(n-1)(n-2) - 2(n-4)\pi + 2 \\ -2 \\ -\frac{1}{2}n^3 + \frac{9}{2}n^2 - 11n + 16 - K_S^2 + \pi(3n - 16) \\ \frac{1}{6}(n-1)(n-2)(-2n+9) + \pi(n-4) - p_a(S) - 1 \end{pmatrix}.$$

This finally gives:
— $\check{n} = n + 4\pi + 12p_a(S) - K_S^2 + 8$ ;
— $d = \frac{1}{2}(n-1)(n-2) - \pi$ ;
— $\rho = \frac{1}{2}(n^2 - 7n) + \pi(n-12) + 9p_a(S) - 2K_S^2 + 22$ ;
— $\Pi = n^2 - 6n + 2\pi(n-10) + 12p_a(S) - 3K_S^2 + 33$ ;
— $t = \frac{1}{6}(n^3 - 9n^2 + 26n) - \pi(n-8) - 4p_a(S) + K_S^2 - 12$ ;
— $\tau = 2n + 8\pi - 12p_a(S) + 2K_S^2 - 20$.

**(5.11) Choice of four independent invariants.** I have chosen to follow Enriques' notation for ease of reference. I have however substituted the old-fashioned *linear genus* $p^{(1)}(S)$ with $K_S^2$. Indeed, $p^{(1)}(S)$ is defined as the *virtual* genus of canonical sections of $S$, that is the genus of the curves in $|K_S|$ obtained from the adjunction formula, even though $|K_S|$ may in fact be empty. In other words $p^{(1)}$ is defined by the formula

$$2p^{(1)} - 2 = (K_S + K_S) \cdot K_S,$$

hence $p^{(1)} = K_S^2 + 1$.

Nowadays the custom is to use $C^2$, $K_S \cdot C$, $K_S^2$, $c_2(S)$ rather than $n$, $\pi$, $p_a(S)$, $K_S^2$. To obtain the above formulas in this form, simply use

(5.11.1) $$C^2 = n; \qquad K_S \cdot C + C^2 = 2\pi - 2; \qquad 1 + p_a(S) = \frac{1}{12}\big(K_S^2 + c_2(S)\big).$$

The former two relations we have already encountered, and the latter is the well-known Noether formula (see, e.g., [22, Example A.4.1.2]; a modern proof in classical spirit is given in [39]; we shall prove it in (5.18)).

---

[4] There is a typo in [15, Formula (6), p. 175], but the final formulas on p. 176 agree with ours in (5.10).

## 5.2 – The postulation formula

In this section we discuss the postulation formula (5.9.1), which gives an expression for the arithmetic genus $p_a(S)$ in terms of the projection $S^\flat$ in $\mathbf{P}^3$. Here, we give a first simple proof of the formula, and then explain the classical definition of the arithmetic genus and why it coincides with the modern definition, which gives another simple proof of the formula. The postulation formula is discussed in [39] as well. The interested reader may also wish to consult [55, Chap. III], in addition to [15, Cap. IV].

**(5.12) Definition.** *Let $X$ be a variety. The* arithmetic genus *of $X$ is*

$$p_a(X) := (-1)^{\dim X}\bigl(\chi(\mathcal{O}_X) - 1\bigr).$$

*If $X$ is smooth, its* geometric genus *is*

$$p_g(X) := h^0(X, \omega_X).$$

If $X$ is a degree $n$ hypersurface in $\mathbf{P}^N$, one finds

(5.12.1) $$p_a(X) = h^0(\mathcal{O}_{\mathbf{P}^N}(n - N - 1)) = \binom{n-1}{N}.$$

**(5.13) Arithmetic and geometric genera of a curve.** For a smooth curve there is no difference between the arithmetic and geometric genera; it is important to note that the same does not hold in higher dimension. One defines the geometric genus of an arbitrary reduced curve as the genus of its normalization. The arithmetic and geometric genera of a curve $X$ may be related by means of the exact sequence

(5.13.1) $$0 \to \mathcal{O}_X \to \nu_*\mathcal{O}_{\bar{X}} \to \bigoplus_{p \in \mathrm{Sing}(X)}(\mathbf{C}_p)^{\oplus \delta(p)} \to 0,$$

where $\nu : \bar{X} \to X$ is the normalization map, and the third term is a skyscraper sheaf supported on the singular locus of $X$, with $\delta(p)$ an integer depending on the type of singularity, see [II, Section **]. Taking Euler characteristics in (5.13.1), one obtains

$$p_g(X) = p_a(X) - \sum_{p \in \mathrm{Sing}(X)} \delta(p).$$

*(5.13.2) Example.* Locally at a space triple point, the exact sequence (5.13.1) looks like

$$0 \to \mathbf{C}[x,y,z]/(xy,xz,yz) \to \begin{array}{c} \mathbf{C}[x] \oplus \mathbf{C}[y] \oplus \mathbf{C}[z] \\ (P,Q,R) \end{array} \begin{array}{c} \to \\ \mapsto \end{array} \begin{array}{c} \mathbf{C}^2 \\ \bigl(P(0) - Q(0), Q(0) - R(0)\bigr) \end{array} \to 0$$

so that a space triple point makes the genus drop by 2.

**(5.14) Proof of the postulation formula, I.** Let us apply the ideas of (5.13) in order to compute the arithmetic genus of a smooth surface $S$ in terms of a general projection $S^\flat$ in $\mathbf{P}^3$. Thus, $S^\flat$ has ordinary singularities as described in (5.3), and we shall use the same notation as in Section 5.1. There is an exact sequence

$$0 \to \mathcal{O}_{S^\flat} \to \nu_*\mathcal{O}_S \to \begin{array}{c} \nu_*\mathcal{O}_\Delta \\ (f_1, f_2) \end{array} \begin{array}{c} \to \\ \mapsto \end{array} \begin{array}{c} \mathcal{O}_D \\ f_1 - f_2 \end{array} \to 0,$$

and it gives, taking Euler characteristics:

(5.14.1) $$p_a(S) = p_a(S^\flat) - p_a(\Delta) + p_a(D).$$

Now, the arithmetic genus of $S^\flat$ is that of a degree $n$ hypersurface in $\mathbf{P}^3$, see (5.12.1), that of $\Delta$ is $(n-4)d + 1$ by (5.7.1), and that of $D$ is $\rho + 2t$ by (5.13.2), so the relation (5.14.1) gives the postulation formula (5.9.1). □

To explain the classical definition of the arithmetic genus of a surface, we will use the following lemma.

**(5.15) Lemma.** *For all $i = 0, 1, 2$, one has $H^i(S, \omega_S) \cong H^i(S^\flat, \mathcal{O}_{S^\flat}(n-4) \otimes \mathcal{I}_D)$.*

*Proof.* We use the relations between the canonical line bundles of $S$, $S^\flat$, and $\mathbf{P}^3$, see (5.3). By the projection formula, we have

$$(5.15.1) \qquad \nu_* \omega_S = \nu_*\big(\nu^* \mathcal{O}_{S^\flat}(n-4) \otimes \mathcal{O}_S(-\Delta)\big) = \mathcal{O}_{S^\flat}(n-4) \otimes \mathcal{I}_D.$$

Moreover one has $R^i \nu_*(\omega_S) = 0$ for all $i = 1, 2$, because $\nu$ is finite, so $H^i(S, \omega_S) \cong H^i(S^\flat, \nu_* \omega_S)$ for all $i = 0, 1, 2$ by Leray's spectral sequence (or merely [22, Exer. III.8.1]), and the result follows. □

**(5.16) Classical definition of the arithmetic genus of a surface.** Let $S$ be a smooth surface, and $S^\flat$ be a general projection of $S$ in $\mathbf{P}^3$, as in (5.14). By Lemma (5.15), one has

$$p_g(S) = h^0(\mathbf{P}^3, \mathcal{O}_{\mathbf{P}^3}(n-4) \otimes \mathcal{I}_D),$$

i.e., the geometric genus of $S$ is the number of independent degree $n-4$ hypersurfaces containing $D$; in classical words, the latter are called degree $n-4$ hypersurfaces adjoint to $S^\flat$. To compute this number, one may consider the restriction exact sequence

$$(5.16.1) \qquad 0 \to \mathcal{O}_{\mathbf{P}^3}(n-4) \otimes \mathcal{I}_D \to \mathcal{O}_{\mathbf{P}^3}(n-4) \to \mathcal{O}_D(n-4) \to 0.$$

If the linear series cut out on $D$ by $|\mathcal{O}_{\mathbf{P}^3}(n-4)|$ is complete and non-special, one finds that $h^0(\mathbf{P}^3, \mathcal{O}_{\mathbf{P}^3}(n-4) \otimes \mathcal{I}_D)$ equals

$$\begin{aligned} h^0\big(\mathcal{O}_{\mathbf{P}^3}(n-4)\big) - h^0\big(\mathcal{O}_D(n-4)\big) &= \tbinom{n-3}{3} - \big(1 - p_a(D) + \deg \mathcal{O}_D(n-4)\big) \\ &= \tbinom{n-3}{3} - 1 + (\rho + 2t) + (n-4)d \end{aligned}$$
(5.16.2)

(the first equality is given by the Riemann–Roch theorem, and the second one by the computation of $p_a(D)$ in (5.14)). One recognizes in (5.16.2) the right-hand-side of the postulation formula (5.9.1).

Classically, the arithmetic genus of $S$ is defined as the *virtual number* of independent degree $n-4$ hypersurfaces adjoint to $S^\flat$, i.e., the value that $h^0(\mathbf{P}^3, \mathcal{O}_{\mathbf{P}^3}(n-4) \otimes \mathcal{I}_D)$ would assume if the linear series cut out on $D$ by $|\mathcal{O}_{\mathbf{P}^3}(n-4)|$ were complete and non-special.[5] Thus, in the classical definition, the postulation formula (5.9.1) is taken as the *definition* of the arithmetic genus.

**(5.17) Equivalence of the definitions: proof of the postulation formula, II.** We shall now prove that the definition of the arithmetic genus given in (5.12) is equivalent to the classical definition given above. The matter is to prove the postulation formula, which we have already

---

[5]or "under the postulation that the linear series cut out on $D$ by $|\mathcal{O}_{\mathbf{P}^3}(n-4)|$ is complete and non-special", hence the name "postulation formula".

done; we shall give another proof, based on the considerations of (5.16). The exact sequence (5.16.1) gives

$$\chi\bigl(\mathcal{O}_{\mathbf{P}^3}(n-4) \otimes \mathcal{I}_D\bigr) = \chi\bigl(\mathcal{O}_{\mathbf{P}^3}(n-4)\bigr) - \chi\bigl(\mathcal{O}_D(n-4)\bigr)$$
$$= \binom{n-3}{3} - 1 + \rho + 2t + (n-4)d,$$

regardless of any assumption on the linear series cut out on $D$ by $|\mathcal{O}_{\mathbf{P}^3}(n-4)|$. Now, I claim that

$$\chi\bigl(\mathcal{O}_{\mathbf{P}^3}(n-4) \otimes \mathcal{I}_D\bigr) = \chi(\mathcal{O}_S) - 1,$$

which will prove the equivalence of the classical and cohomological definitions of the arithmetic genus.

To prove the claim, consider the exact sequence

$$0 \to \omega_{\mathbf{P}^3} \cong \mathcal{O}_{\mathbf{P}^3}(n-4) \otimes \mathcal{I}_{S^\flat} \to \mathcal{O}_{\mathbf{P}^3}(n-4) \otimes \mathcal{I}_D \to \mathcal{O}_{S^\flat}(n-4) \otimes \mathcal{I}_D \to 0,$$

with $\mathcal{I}_{S^\flat}$ the sheaf of ideals defining $S^\flat$ in $\mathbf{P}^3$. It implies

$$\chi\bigl(\mathcal{O}_{\mathbf{P}^3}(n-4) \otimes \mathcal{I}_D\bigr) = -1 + \chi\bigl(\mathcal{O}_S(n-4) \otimes \mathcal{I}_D\bigr).$$

Now, one has $\chi(S^\flat, \mathcal{O}_{S^\flat}(n-4) \otimes \mathcal{I}_D) = \chi(S, \omega_S)$ by Lemma (5.15), and moreover $\chi(S, \omega_S) = \chi(S, \mathcal{O}_S)$ by Serre duality: this proves the claim. □

## 5.3 – The 2-dimensional case: multiple planes

In this section we consider a linear system $|V|$ on a smooth surface $S$, of dimension $r = 1$ or $2$, and we assume that $|V|$ is a general linear subsystem of a 3-dimensional linear system $|V^\sharp|$ satisfying the conditions of (5.3). In particular, this is the case if $|V|$ is a general linear subsystem of a very ample linear system.

We shall apply the results obtained in Section 5.1 to $|V^\sharp|$. We will use the notation introduced there, by taking $S^\flat \subseteq \mathbf{P}^3$ to be the surface defined by $|V^\sharp|$.

**(5.18) One-dimensional systems: proof of the Noether formula.** Let $r = 1$. We have computed the class of $|V|$ in (5.2). On the other hand, the discriminant of $|V|$ is merely the section of that of $|V^\sharp|$ by the line $|V|$, hence the class of $|V|$ equals that of $|V^\sharp|$, which is given in (5.10). Equating these two expressions for $\check{n}$, one finds

$$c_2(S) + n + 4\pi - 4 = n + 4\pi + 12p_a(S) - K_S^2 + 8,$$

which is equivalent to the Noether formula, stated in (5.11.1). □

**(5.19)** We now let $r = 2$ for the rest of this section. Then $|V|$ defines an $n:1$ covering $S \to \mathbf{P}^2$ with branch curve $B \subseteq \mathbf{P}^2$, which may be incarnated as the projection of $S^\flat \subseteq \mathbf{P}^3$ from a general point $p$. We consider:
— $N$ the degree of the branch curve $B \subseteq \mathbf{P}^2$ ;
— $P$ the geometric genus of $B$ ;
— $\kappa$ the number of points $q \in S$ such that there exist $C, C' \in |V|$ intersecting with multiplicity $\geqslant 3$ at $q$ ;
— $\iota$ the number of cuspidal members of $|V|$ ;
— $\delta$ the number of unordered pairs $q, q' \in S$ such that there exist $C, C' \in |V|$ tangent at both $q$ and $q'$ ;
— $\beta$ the number of 2-nodal members of $|V|$.

The dual of $B$ is the discriminant of $|V|$, so the class of $B$ is $\check{n}$, which assumes the same value as in (5.2) and (5.10), see (5.18). Moreover, $\delta$, $\kappa$, $\beta$, $\iota$ are respectively the numbers of nodes, cusps, bitangents, flexes of $B$. We take once again for granted that $B$ and its dual have no other singularities by our generality assumption.

**(5.20) Characters of the apparent boundary of $S^\flat$.** The branch curve $B$ is the image of the apparent boundary $\mathcal{B}^p S^\flat$ by the map $S^\flat \to \mathbf{P}^2$, and in fact $\mathcal{B}^p S^\flat$ is the normalization of $B$. This implies that $N$ and $P$ are respectively the degree and genus of $\mathcal{B}^p S^\flat$, and since we have computed that its class in $S$ is $K_S + 3C$ one finds:
— $N = (K_S + 3C) \cdot C = 2n + 2\pi - 2$ ;
— $P = 9\pi + K_S^2 - 8$.

**(5.21) Conclusion.** Since we know in addition the class $\check{n}$ of $B$, we may derive all remaining numerical characters by applying the Plücker formulae of Section 1.3 to $B$. One finds:
— $\delta = 2\big[(n+\pi)^2 - 5n - 17\pi - 2K_S^2 + 6p_a(S) + 22\big]$ ;
— $\kappa = 3(n + 6\pi + K_S^2 - 4p_a(S) - 10)$ ;
— $\beta = \frac{1}{2}\big[(n + 4\pi - K_S^2 + 12p_a(S) - 1)^2 + 15n - 6\pi - 17K_S^2 + 132p_a(S) + 57\big]$ ;
— $\iota = 24(\pi + p_a(S))$.

**(5.22) Number of cuspidal curves in a net.** It is interesting to derive the quantity $\iota$ directly by projective methods. We end this section by sketching how this goes.

The parabolic curve of $S^\flat$ and the apparent boundary $\mathcal{B}^p S^\flat$ intersect at the pinch points, and at those points $q$ on the parabolic curve such that the tangent plane at $q$ contains $p$; the quantity $\iota$ is the number of points of the latter kind.

To find the intersection number of these two curves we compute the class of the pull-back of the parabolic curve to $S$. The degree of the Hessian is $4(n-2)$, and a local computation shows that its intersection with $S^\flat$ contains the double curve $D$ with multiplicity 4. Thus, arguing as in (5.6), one finds that the pull-back of the parabolic curve has class $4K_S + 8C$. The class of the pull-back of $\mathcal{B}^p S^\flat$ has been computed in (5.6).

On the other hand, we claim that the two curves intersect with multiplicity 2 at each pinch point. Indeed, if $q$ is a pinch point with preimage $q' \in S$, then the members of $|V|$ passing through $q'$ form a pencil of curves which are all singular at $q'$. The condition to be cuspidal in such a pencil is that the determinant of the $2 \times 2$ matrix of second derivatives vanishes, which is a degree 2 condition. We thus find 2 cuspidal curves at $q'$ in $|V|$, hence the claim.

Finally, we have

$$(5.22.1) \qquad (K_S + 3C) \cdot (4K_S + 8C) = \iota + 2\tau,$$

which indeed gives $\iota = 24(\pi + p_a)$ with the expression for $\tau$ in (5.10).

# A – Polynomial calculus

We give here a brief recap on polarity, together with a few facts about the Hessian, so that the reader unfamiliar with this may conveniently consult the relevant material. We refer the reader to [13, Chapter 1] and [5, § 5.4 and 5.6] for more details; see also [16].

## A.1 – Polarity

**(A.1) Notation and plea for the notation.** Let $\mathbf{k}$ be a field of characteristic zero, and let $F \in \mathbf{k}[X_0, \ldots, X_n]$ be a homogeneous degree $d$ polynomial in $n+1$ indeterminates.

For all $\mathbf{x} = (x_0, \ldots, x_n)$ and $\mathbf{y} = (y_0, \ldots, y_n)$ in $\mathbf{k}^{n+1}$, we write $\mathrm{D}_{\mathbf{x}}F$ for the differential of $F$ at the point $\mathbf{x}$ (which is a linear functional on $\mathbf{k}^{n+1}$), and $\mathrm{D}_{\mathbf{x}}F(\mathbf{y})$ for the value of the latter at the vector $\mathbf{y}$, so that

$$(\mathrm{A}.1.1) \qquad \mathrm{D}_{\mathbf{x}}F(\mathbf{y}) = y_0\,\partial_0 F(\mathbf{x}) + \cdots y_n\,\partial_n F(\mathbf{x}),$$

where $\partial_i F$ is the polynomial $\partial F / \partial X_i$ for all $i = 0, \ldots, n$.

The quantity (A.1.1) is a function in the two variables $\mathbf{x}$ and $\mathbf{y}$, and somehow a main point in polarity theory is to consider it as a function of the variable $\mathbf{x}$ depending on the parameter $\mathbf{y}$. With this in mind we introduce the following alternative notation:

$$(\mathrm{A}.1.2) \qquad \bigl(\mathrm{D}^{\mathbf{y}}F\bigr)(\mathbf{x}) = \bigl(\mathrm{D}_{\mathbf{x}}F\bigr)(\mathbf{y}) = y_0\,\partial_0 F(\mathbf{x}) + \cdots y_n\,\partial_n F(\mathbf{x}).$$

Although it may be disturbing at first, as I have witnessed on various occasions, we shall see that it is both meaningful and useful. One could be tempted at this point to introduce $\mathrm{D}^{\mathbf{y}}_{\mathbf{x}}F$ as well, but let me refrain from doing this.

**(A.2) The polar pairing.** In the notation of (A.1) above, for all $\mathbf{x} = (x_0, \ldots, x_n) \in \mathbf{k}^{n+1}$,

$$(\mathrm{A}.2.1) \qquad \mathrm{D}^{\mathbf{x}}F = x_0\,\partial_0 F + \cdots x_n\,\partial_n F \in \mathbf{k}[X_0, \ldots, X_n]$$

is a homogeneous polynomial of degree $d-1$ (if nonzero), which depends linearly on $\mathbf{x}$ (note that the $X_i$ are indeterminates whereas the $x_i$ are scalars). Considering this polynomial as $(\mathrm{D}^{\mathbf{x}})(F)$, where $\mathrm{D}^{\mathbf{x}}$ denotes the differential operator

$$\mathrm{D}^{\mathbf{x}} = x_0\partial_0 + \cdots x_n\partial_n,$$

for all $\mathbf{x}_1, \ldots, \mathbf{x}_k \in \mathbf{k}^{n+1}$, we set

$$(\mathrm{A}.2.2) \qquad \mathrm{D}^{\mathbf{x}_1 \cdots \mathbf{x}_k} = \mathrm{D}^{\mathbf{x}_1} \circ \cdots \circ \mathrm{D}^{\mathbf{x}_k}.$$

Extending (A.2.2) by $\mathbf{k}$-linearity, one obtains a perfect bilinear pairing

$$\mathrm{D} : \mathrm{Sym}^k(\mathbf{k}^{n+1}) \otimes \mathrm{Sym}^d(\check{\mathbf{k}}^{n+1}) \longrightarrow \mathrm{Sym}^{d-k}(\check{\mathbf{k}}^{n+1})$$

$$\sum\nolimits_{i,j} \Xi_i \otimes F_j \longmapsto \sum\nolimits_{i,j} \mathrm{D}^{\Xi_i} F_j$$

extending the natural pairing between $\mathbf{k}^{n+1}$ and its dual $\check{\mathbf{k}}^{n+1}$.

**(A.3) The Taylor–Newton formula.** In the formalism of polarity it has the following pleasant incarnation. Let $\mathbf{x} = (x_0, \ldots, x_n)$ and $\mathbf{y} = (y_0, \ldots, y_n)$ be in $\mathbf{k}^{n+1}$. As above, $F$ is a degree $d$ polynomial. One has:

$$F(\mathbf{x} + \mathbf{y}) = \sum_{k=0}^{d} \sum_{k_0 + \cdots + k_n = k} \frac{1}{k_0! \cdots k_n!} y_0^{k_0} \cdots y_n^{k_n} \left( \frac{\partial^{k_0}}{\partial X_0^{k_0}} \cdots \frac{\partial^{k_n}}{\partial X_n^{k_n}} F \right)(\mathbf{x})$$

$$= \sum_{k=0}^{d} \frac{1}{k!} \left( \left( y_0 \frac{\partial}{\partial X_0} + \cdots + y_n \frac{\partial}{\partial X_n} \right)^k F \right)(\mathbf{x})$$

$$= \sum_{k=0}^{d} \frac{1}{k!} \mathrm{D}^{\mathbf{y}^k} F(\mathbf{x}).$$

**(A.4) Polar hypersurfaces.** Consider the hypersurface $Z = V(F) \subseteq \mathbf{P}^n$. For $\hat{a} = (a_0, \ldots, a_n)$ and $\mathbf{x} = (x_0, \ldots, x_n)$ in $\mathbf{k}^{n+1}$, and $k \in \mathbf{N}$,

$$\mathrm{D}^{\hat{a}^k} F(\mathbf{x}) = \big((a_0 \partial_0 + \cdots + a_n \partial_n)^k F\big)(\mathbf{x})$$

is a bihomogeneous polynomial of bidegree $(k, d-k)$ in the variables $(a_0, \ldots, a_n)$ and $(x_0, \ldots, x_n)$. The hypersurface $V(\mathrm{D}^{\hat{a}^k} F)$ depends only on $Z$ and the point $a = (a_0 : \ldots : a_n) \in \mathbf{P}^n$, see [5, Thm. 5.4.2]: we call it the *k-th polar of $Z$ with respect to $a$*, and denote it by $\mathrm{D}^{a^k} Z$. It is a degree $d - k$ hypersurface. We will often abuse notation and omit the distinction between $\hat{a}$ and $a$. In addition, we will call $V(\mathrm{D}^{\hat{a}^k} F)$ a hypersurface even when $\mathrm{D}^{\hat{a}^k} F = 0$, and thus $V(\mathrm{D}^{\hat{a}^k} F) = \mathbf{P}^n$.

We shall also use the following notion: the *polar k-ic of $Z$ at $a$*, denoted by $\mathrm{D}_{a^k} Z$, is the degree $k$ hypersurface defined by the polynomial

$$\mathrm{D}^{\mathbf{x}^k} F(\hat{a}) = \big((x_0 \partial_0 + \cdots + x_n \partial_n)^k F\big)(\hat{a})$$

in the variables $(x_0, \ldots, x_n)$.[6] In other words, we interchange the roles of $\hat{a}$ and $\mathbf{x}$ when we pass from $\mathrm{D}^{a^k} Z$ to $\mathrm{D}_{a^k} Z$; thus, by definition, for all $a, b \in \mathbf{P}^n$,

(A.4.1) $$b \in \mathrm{D}_{a^k} Z \quad \Longleftrightarrow \quad a \in \mathrm{D}^{b^k} Z.$$

**(A.5) Proposition** (polar symmetry)**.** *Consider a degree $d$ hypersurface $Z \subseteq \mathbf{P}^n$. Let $a, b \in \mathbf{P}^n$, and $k \in [\![1, d-1]\!]$. One has the equivalence:*

$$a \in \mathrm{D}^{b^k} Z \quad \Longleftrightarrow \quad b \in \mathrm{D}^{a^{d-k}} Z.$$

In particular this says that the polar $k$-ic at $a$ (i.e., $\mathrm{D}_{a^k} Z$) coincides with the $d - k$-th polar with respect to $a$ (i.e., $\mathrm{D}^{a^{d-k}} Z$).

*Proof.* Let $F$ be a homogeneous polynomial defining the hypersurface $Z$. Consider two sets of $n + 1$ indeterminates $\mathbf{X} = (X_0, \ldots, X_n)$ and $\mathbf{Y} = (Y_0, \ldots, Y_n)$. By the Taylor–Newton Formula we have

$$F(\mathbf{X} + \mathbf{Y}) = \sum_{k=0}^{d} \frac{1}{k!} \big(\mathrm{D}^{\mathbf{Y}^k} F\big)(\mathbf{X}) \quad = \quad F(\mathbf{Y} + \mathbf{X}) = \sum_{k=0}^{d} \frac{1}{k!} \big(\mathrm{D}^{\mathbf{X}^k} F\big)(\mathbf{Y}).$$

Identifying the degree $k$ pieces in the indeterminates $\mathbf{X}$, we find

$$\frac{1}{(d-k)!} \mathrm{D}^{\mathbf{Y}^{d-k}} F(\mathbf{X}) \;=\; \frac{1}{k!} \mathrm{D}^{\mathbf{X}^k} F(\mathbf{Y}),$$

which gives the result at once. □

**(A.6) Proposition.** *Let $Z$ be a degree $d$ hypersurface in $\mathbf{P}^n$. For all $a \in \mathbf{P}^n$ and $k = 1, \ldots, d-1$, one has:*

$$a \in Z \quad \Longleftrightarrow \quad a \in \mathrm{D}^{a^k} Z.$$

*Proof.* Let $F$ be an equation of $Z$. One has $\mathrm{D}^{a^k} F(a) = d(d-1) \cdots (d-k+1) F(a)$ by the Euler formula. □

---

[6] to avoid possible confusion, let me take the risk of being pedantic: consider $F$ as a polynomial in the indeterminates $Y_0, \ldots, Y_n$, and let $X_0, \ldots, X_n$ be a second set of indeterminates; then, specializing the polynomial $\left(X_0 \frac{\partial}{\partial Y_0} + \cdots + X_n \frac{\partial}{\partial Y_n}\right)^k F(Y_0, \ldots, Y_n)$ by specializing $Y_i$ to $a_i$ for all $i$, one obtains a polynomial in the indeterminates $X_0, \ldots, X_n$: this is the defining equation of $\mathrm{D}_{a^k} Z$.

**(A.7)** If $b \in Z = V(F)$ is a smooth point, then $(D_b F)(\mathbf{X}) = (D^{\mathbf{X}} F)(b)$ is a linear homogeneous polynomial defining the tangent hyperplane to $Z$ at $b$, hence polar hyperplane $D_b Z$ is the tangent hyperplane to $Z$ at $b$, and moreover, for all $a \in \mathbf{P}^n$,

$$Z_{\mathrm{sm}} \cap D^a Z = \{b \in Z_{\mathrm{sm}} : \mathbf{T}_b Z \ni a\},$$

where $Z_{\mathrm{sm}}$ denotes the smooth locus of $Z$. This generalizes to the following fundamental property.

For all hypersurfaces $Z \subseteq \mathbf{P}^n$ and lines $\ell \subseteq \mathbf{P}^n$, and for all $p \in \ell$, we let $i(Z, \ell)_p$ be the multiplicity with which $p$ appears in $Z \cap \ell$.

**(A.8) Theorem.** *Let $Z$ be a degree $d$ hypersurface, $a \in Z$, and $b \in \mathbf{P}^n$. For all integer $s \geqslant 0$, one has*

$$i(Z, \langle a, b \rangle)_a \geqslant s+1 \quad \iff \quad \forall k = 1, \ldots, s, \quad b \in D_{a^k} Z$$
$$\iff \quad \forall k = 1, \ldots, s, \quad a \in D^{b^k} Z.$$

*Proof.* Let $F$ be a homogeneous equation for $Z$, and let $\hat{a}, \hat{b} \in \mathbf{k}^{n+1}$. Consider an indeterminate $T$. We have

$$F(\hat{a} + T\hat{b}) = \sum_{k=0}^{d} \frac{1}{k!} \left( D^{(T\hat{b})^k} F \right)(\hat{a}) = \sum_{k=0}^{d} \frac{1}{k!} \left( D^{\hat{b}^k} F \right)(\hat{a}) \cdot T^k$$

by the Taylor–Newton formula. This gives the result since, by definition, $i(Z, \langle a, b \rangle)_a \geqslant s+1$ if and only if $T^{s+1}$ divides $F(\hat{a} + T\hat{b})$ (note that $D^{\hat{b}^0} F(\hat{a}) = F(\hat{a}) = 0$ by hypothesis). $\square$

It turns out that for $a \in Z$, all polars of $Z$ with respect to $a$ (equivalently, all polar hypersurfaces $D_{a^k} Z$ at $a$) are tangent at $a$. Actually, $Z$ and its polar $k$-ic at $a$, $D_{a^k} Z$ have the same polar $s$-ics at $a$ for all $s \leqslant k$, as the following identities show. Let $d = \deg(Z)$. Using polar symmetry repeatedly, one finds

$$D_{a^s}(D_{a^k} Z) = D^{a^{k-s}}(D_{a^k} Z) = D^{a^{k-s}}(D^{a^{d-k}} Z) = D^{a^{d-s}} Z = D_{a^s} Z.$$

This has the following remarkable consequence.

**(A.9) Corollary.** *Let $Z$ be a hypersurface in $\mathbf{P}^n$, and $a$ be a point of $Z$. The $n$ hypersurfaces $D_a Z, \ldots, D_{a^n} Z$ intersect with multiplicity at least $n!$ at $a$ (with the convention that the multiplicity is infinite if the intersection is positive dimensional at $a$).*

*Proof.* We consider only the case in which $D_a Z, \ldots, D_{a^n} Z$ define a local complete intersection at $a$, otherwise the result is trivial. Then the intersection $D_a Z, \ldots, D_{a^{n-1}} Z$ consists of $(n-1)!$ lines each intersecting $Z$ with multiplicity at least $n$ at $a$, as follows from Theorem (A.8). Since the polar $k$-ics at $a$ of $D_{a^n} Z$ are the same as those of $Z$, as indicated above, each of these lines also intersect $D_{a^n} Z$ with multiplicity at least $n$ at $a$, and the result follows. $\square$

The polar hyperplane $D_a Z$ is indeed a hyperplane only if $a$ is a smooth point of $Z$, and then it is properly the tangent hyperplane to $Z$ at $a$, otherwise $D_a Z$ is the whole $\mathbf{P}^n$. The following result generalizes this: by (A.10.1), if $Z$ has multiplicity exactly $m$ at $a$, then $D_{a^m} Z$ coincides with the tangent cone $\mathrm{TC}_a(Z)$.

**(A.10) Theorem.** *Let $Z$ be a degree $d$ hypersurface, and $a$ be a point of $Z$ of multiplicity $m$. Consider an integer $r \leqslant d$.*

*(A.10.1)* The polar $r$-ic at $a$, $\mathrm{D}_{a^r}Z$, is the whole $\mathbf{P}^n$ if and only if $r < m$. If $r \geqslant m$, then $\mathrm{D}_{a^r}Z$ has multiplicity $m$ at $a$, and it has the same tangent cone as $Z$ at $a$: $\mathrm{TC}_a(\mathrm{D}_{a^r}Z) = \mathrm{TC}_a(Z)$.
*(A.10.2)* Let $b \in \mathbf{P}^n - \{a\}$. If $r < m$, the $r$-th polar with respect to $b$, $\mathrm{D}^{b^r}Z$, has multiplicity at least $m - r$ at $a$; the multiplicity is exactly $m - r$ if $b$ is general, and in this case the tangent cone of $\mathrm{D}^{b^r}Z$ at $a$ equals the $r$-th polar with respect to $b$ of the tangent cone $\mathrm{TC}_a(Z)$: $\mathrm{TC}_a\bigl(\mathrm{D}^{b^r}Z\bigr) = \mathrm{D}^{b^r}\bigl(\mathrm{TC}_a(Z)\bigr)$.

Note that (A.10.2) tells us that $\mathrm{D}^{b^r}Z$ contains $a$ if $r \leqslant m - 1$, and is singular at $a$ if $r \leqslant m - 2$.

*Proof.* Let us assume, without loss of generality, that $a = (1{:}0{:}\ldots{:}0)$. Then
$$F = X_0^{d-m} F_m(X_1, \ldots, X_n) + \cdots + F_d(X_1, \ldots, X_n)$$
where each $F_k$ is homogeneous of degree $k$, and $F_m$ is non-zero. One has $\mathrm{D}_{a^r}Z = \mathrm{D}^{a^{d-r}}Z$ by polar symmetry, and $\mathrm{D}^{a^{d-r}} = (\partial/\partial X_0)^{d-r}$, hence $\mathrm{D}_{a^r}Z$ is defined by the equation
$$\bigl(\tfrac{\partial}{\partial X_0}\bigr)^{d-r} F = \tfrac{(d-m)!}{(r-m)!} X_0^{r-m} F_m + \cdots + \tfrac{(d-r)!}{0!} F_r;$$
it is zero if and only if $r < m$, and if $r \geqslant m$ we see that $\mathrm{D}_{a^r}Z$ has multiplicity $m$ at $a$ with tangent cone defined by $F_m$, which proves (A.10.1).

Now let $b = (b_0, b_1, \ldots, b_n)$ be such that the corresponding point in $\mathbf{P}^n$ is distinct from $a$. We set $b^\natural = (0, b_1, \ldots, b_n)$ (it is non-zero because $b$ and $a$ are distinct in $\mathbf{P}^n$), and $b^r = (b^\natural)^r + \Xi$. Then $\mathrm{D}^{b^r} = \mathrm{D}^{(b^\natural)^r} + \mathrm{D}^\Xi$, where $\mathrm{D}^{(b^\natural)^r}$ is a polynomial in $\partial/\partial X_1, \ldots, \partial/\partial X_n$ only, while $\mathrm{D}^\Xi$ is a polynomial in $\partial/\partial X_0, \partial/\partial X_1, \ldots, \partial/\partial X_n$ in which all terms involve $\partial/\partial X_0$ with a positive exponent. Therefore, the piece of $\mathrm{D}^{b^r}F$ with the highest degree in $X_0$ is
$$\mathrm{D}^{(b^\natural)^r}(X_0^{d-m} F_m) = X_0^{d-m} \mathrm{D}^{(b^\natural)^r}(F_m) = X_0^{d-m} \mathrm{D}^{b^r}(F_m).$$

It follows that $\mathrm{D}^{b^r}Z$ has multiplicity at least $m - r$ at $a$; this multiplicity is exactly $m - r$ if and only if $\mathrm{D}^{b^r}(F_m)$ is non-zero, and in this case $\mathrm{D}^{b^r}(F_m)$ defines the tangent cone of $\mathrm{D}^{b^r}Z$ at $a$, so that the latter equals $\mathrm{D}^{b^r}(\mathrm{TC}_a(Z))$. Now $\mathrm{D}^{b^r}(F_m) = \mathrm{D}^{(b^\natural)^r}(F_m)$ is non-zero for general $b$, and (A.10.2) is proven. $\square$

## A.2 – Hessian and second fundamental form

Let $Z$ be a hypersurface in $\mathbf{P}^n$, defined by the homogeneous polynomial $F \in \mathbf{k}[X_0, \ldots, X_n]$.

**(A.11) Hessian of a hypersurface.** The *Hessian hypersurface* of $Z$ is the hypersurface $\mathrm{Hess}(Z)$ defined by the homogeneous polynomial
$$\mathrm{Hess}(F) = \det\bigl(\partial_i \partial_j F\bigr)_{0 \leqslant i, j \leqslant n}.$$

It has degree $(n+1)(n-2)$ (if non-zero).

Using the Euler formula and the standard properties of determinants, one may obtain the alternative expression

(A.11.1) $$\mathrm{Hess}(F) = \tfrac{(d-1)^2}{X_0^2} \begin{vmatrix} \tfrac{d}{d-1} \cdot F & \partial_1 F & \cdots & \partial_n F \\ \partial_1 F & \partial_1^2 F & \cdots & \partial_1 \partial_n F \\ \vdots & \vdots & \ddots & \vdots \\ \partial_n F & \partial_1 \partial_n F & \cdots & \partial_n^2 F \end{vmatrix},$$

which is useful to study the intersection of $Z$ and $\mathrm{Hess}(Z)$.

**(A.12) Second fundamental form.** Let $p$ be a smooth point of $Z$. The *second fundamental form* of $Z$ at $p$, denoted by $II_p$, is a quadric in the tangent space $T_pZ$ which is defined as follows. Analytically locally around $p$, and in an affine chart in which we identify $p$ with $(0,\ldots,0)$ and $T_pZ$ with $\mathbf{A}^{n-1} \times \{0\}$, $Z$ is the image of a map

$$\varphi : (u_1,\ldots,u_{n-1}) \mapsto \bigl(u_1,\ldots,u_{n-1},\varphi(u_1,\ldots,u_{n-1})\bigr)$$

defined on a neighbourhood of the origin in the tangent space $T_pZ$, such that the order 0 and order 1 terms in the Taylor expansion of $\varphi$ vanish. The second fundamental form is the order 2 term of the Taylor expansion of $\varphi$ modulo multiplication by a non-zero scalar; it does not depend on the choice of the parametrization.

**(A.13) Local forms.** Assume $F$ is in the form

$$X_0^{d-1}X_n + X_0^{d-2}\Bigl(\sum_{1 \leqslant i,j \leqslant n} a_{ij} X_i X_j\Bigr) + O(3),$$

so that $p = (1:0:\ldots:0) \in Z$, and $\mathbf{T}_pZ = \{X_n = 0\}$. We assume that the matrix $(a_{ij})_{1 \leqslant i,j \leqslant n}$ is symmetric. By computing the first terms of the Taylor expansion of

$$F\bigl(1,u_1,\ldots,u_{n-1},\varphi(u_1,\ldots,u_{n-1})\bigr),$$

one sees that the second fundamental form of $Z$ at $p$ is the quadratic form on $T_pZ \simeq \mathbf{C}^{n-1}$ defined by the symmetric matrix $(a_{ij})_{1 \leqslant i,j \leqslant n-1}$, modulo scalar multiplication.

Now let us show the relation between the determinant of this matrix and the Hessian determinant. We will assume $n = 3$ in order to have smaller matrices, but it is straightforward to do the same computation for arbitrary $n$. By direct computation, one finds that the Hessian determinant of $F$ at $p = (1:0:0:0)$ is

$$\begin{vmatrix} 0 & 0 & 0 & d-1 \\ 0 & 2a_{11} & 2a_{12} & 2a_{13} \\ 0 & 2a_{12} & 2a_{22} & 2a_{23} \\ d-1 & 2a_{13} & 2a_{23} & 2a_{33} \end{vmatrix} = \begin{vmatrix} 0 & d-1 \\ d-1 & 2a_{33} \end{vmatrix} \cdot \begin{vmatrix} 2a_{11} & 2a_{12} \\ 2a_{12} & 2a_{22} \end{vmatrix}.$$

**(A.14) Corollary.** *The second fundamental form at $p$ is non-degenerate if and only if the Hessian determinant is non-zero.*

**(A.15) Second fundamental form and polar quadric.** Still in the notation of (A.13) above, the polar quadric of $Z$ at $p$, $\mathrm{D}_{p^2}Z$, is defined by the equation

$$\mathrm{D}^{p^{d-2}}F = \bigl(\tfrac{\partial}{\partial X_0}\bigr)^{d-2} F = \tfrac{(d-1)!}{1!} X_0 X_n + \tfrac{(d-2)!}{0!}\Bigl(\sum_{1 \leqslant i,j \leqslant n} a_{ij} X_i X_j\Bigr).$$

We thus see that the intersection of the polar quadric $\mathrm{D}_{p^2}Z$ with the projective tangent hyperplane $\mathbf{T}_pZ = V(X_n)$ is the projective completion of the second fundamental form $II_p$.

Note that $II_p$ is an affine quadric in $T_pZ$ which is defined by a homogenous degree 2 polynomial, hence it is singular at the origin, i.e., at $p$. On the other hand, we have seen in Theorem (A.10) that $\mathbf{T}_pZ = \mathbf{T}_p(\mathrm{D}_{p^2}Z)$, which explains why $\mathrm{D}_{p^2}Z \cap \mathbf{T}_pZ$ is singular at $p$.

# B – Developable surfaces

In this appendix we review the basic theory of developable surfaces in the projective space of dimension 3, following the beautiful [14, Chap. V].

## B.1 – Synthetic description

**(B.1) Definition.** A *developable surface* $S \subseteq \mathbf{P}^3$ is a ruled surface such that the tangent plane to $S$ is constant along the lines of the ruling. In other words, for all lines $\Lambda$ of the ruling, and for all $p, p' \in \Lambda$, one has $\mathbf{T}_pS = \mathbf{T}_{p'}S$. An important aspect of developable surfaces is that they are the surfaces that can be unfolded into the plane with no distortion; we will however not consider this idea any farther in this text.

In classical terminology, if $S$ is a ruled surface, the lines of the rulings are called *generators*, and a generator $\Lambda$ is said to be *torsal* if the tangent plane to $S$ is constant along $\Lambda$. Thus developable surfaces are those ruled surfaces for which all generators are torsal.

**(B.2)** Developable surfaces in $\mathbf{P}^3$ are exactly those surfaces in $\mathbf{P}^3$ that have their dual which is not a hypersurface, see [49, Thm. 1.18]. It therefore follows from biduality that any non-degenerate developable surface is projective dual to a curve $C$ in $\check{\mathbf{P}}^3$. In particular it has an ordinary double curve, corresponding to bitangent planes to $C$, and a cuspidal double curve, corresponding to osculating planes to $C$ (unless $C$ is degenerate, in which case $S$ is a cone).

This description may be refined as follows.

**(B.3)** The following fact refines the description of (B.2); for a proof, see [17, Section 0.5], or [5, 5.8.14]. Also [23, §30] is very enlightening. Let $S$ be a developable surface in $\mathbf{P}^3$. On every generator there is a singular point at which $S$ does not have a well-defined tangent plane (beware that there are other singular points as well, see Proposition (B.14)). If $S$ is not a cone, the locus of the singular points on the various generators is a curve $\Gamma$, called the *regression edge*, which is none other than the cuspidal double curve of $C^\vee$ mentioned in (B.2). At a general point $p \in \Gamma$, the tangent line $\mathbf{T}_p\Gamma$ is the generator of $S$ through $p$, and the osculating plane $\mathrm{Osc}_p\Gamma$ is the plane tangent to $S$ along the generator $\mathbf{T}_p\Gamma$.

This tells us that $S$ is the *tangential surface* of $\Gamma$, namely
$$S = \mathrm{Tan}(\Gamma) := \overline{\bigcup\nolimits_{p \in \Gamma_{\mathrm{sm}}} \mathbf{T}_p\Gamma}.$$

Moreover, the curve $C = S^\vee \subseteq \check{\mathbf{P}}^3$ is what I call the *oscdual*[7] curve of $\Gamma \subseteq \mathbf{P}^3$, namely
$$C = \overline{\bigcup\nolimits_{p \in \Gamma^\circ}(\mathrm{Osc}_p\Gamma)^\perp},$$

where $\Gamma^\circ$ denotes the open subset of $\Gamma$ of points $p$ such that $\mathrm{Osc}_p\Gamma$ is a plane, i.e., $\Gamma$ and $\mathbf{T}_p\Gamma$ intersect with multiplicity exactly 2 at $p$.

The dual of $\Gamma$ is a developable surface in $\check{\mathbf{P}}^3$, the generators of which are the lines $\mathbf{T}_p\Gamma^\perp$ for all smooth points $p \in \Gamma$. It follows from biduality and the above description that this is the tangential surface of the curve $C$. In particular, $\Gamma$ is the oscdual curve of $C$: this is a reflexivity statement for oscduality.

**(B.4)** While we shall not prove here that a developable surface in $\mathbf{P}^3$ is either a cone or the tangential surface $\mathrm{Tan}(\Gamma)$ of a curve, it is easy and instructive to see why tangential surfaces are developable. Assume $\Gamma$ is parametrized locally by a map $\phi : u \mapsto \bigl(\alpha(u), \beta(u), \gamma(u)\bigr)$ in some affine chart of $\mathbf{P}^3$. Then $\mathrm{Tan}(\Gamma)$ is parametrized locally by
$$\Phi : (u, s) \longmapsto \begin{pmatrix} \alpha(u) \\ \beta(u) \\ \gamma(u) \end{pmatrix} + s \begin{pmatrix} \alpha'(u) \\ \beta'(u) \\ \gamma'(u) \end{pmatrix},$$

---

[7]Piene has introduced, more generally, the concept of $m$-dual varieties in [40]; in her notation, the oscdual of $\Gamma$ is the 2-dual variety $\Gamma_2^\vee$. See also Pohl's [41].

and the differential of this parametrization at the point $(u_0, s_0)$ is determined by

$$\frac{\partial}{\partial u} \longmapsto \begin{pmatrix} \alpha'(u_0) \\ \beta'(u_0) \\ \gamma'(u_0) \end{pmatrix} + s_0 \begin{pmatrix} \alpha''(u_0) \\ \beta''(u_0) \\ \gamma''(u_0) \end{pmatrix} \quad \text{and} \quad \frac{\partial}{\partial s} \longmapsto \begin{pmatrix} \alpha'(u_0) \\ \beta'(u_0) \\ \gamma'(u_0) \end{pmatrix}.$$

We thus see that, if $\Gamma$ doesn't have an inflectionary behaviour at the point $\phi(u_0)$, the differential of $\Phi$ fails to be of maximal rank if and only if $s_0 = 0$, and otherwise its image is constant as $s_0$ moves, equal to the plane spanned by $\phi'(u_0)$ and $\phi''(u_0)$. The upshot is that the tangent plane to $\mathrm{Tan}(\Gamma)$ is constant along its generator $\mathbf{T}_{\phi(u_0)}\Gamma$, equal to the osculating plane $\mathrm{Osc}_{\phi(u_0)}\Gamma$, except at the very point $\phi(u_0)$ where it is not well-defined, and thus $\Gamma$ is a cuspidal curve of $\mathrm{Tan}(\Gamma)$.

**(B.5)** Finally, although we shall not use it, note that developable surfaces are characterized by the property that they are contained in the zero locus of their associated Hessian determinant. This is an instance of a result valid for hypersurfaces of any dimension, which has been proved by B. Segre [47]; see also [44, Lemma 7.2.7] for a proof attributed to Kollár by the author.

## B.2 – Analysis of contacts

We carry on our analysis of developable surfaces, building upon the description given in the previous section, and following the general philosophy of biduality. Recall that $S = \mathrm{Tan}(\Gamma) \subseteq \mathbf{P}^3$ is the tangential surface of the curve $\Gamma \subseteq \mathbf{P}^3$, and $C = S^\vee \subseteq \check{\mathbf{P}}^3$ is the dual curve of the developable surface $S$, as well as the oscdual of $\Gamma$. We leave aside for the moment the peculiar behaviour that the curve $\Gamma \subseteq \mathbf{P}^3$ may have at isolated points in some special situations.

**(B.6) Characteristic numbers, I.** We let:
— $m = \deg \Gamma$ ;
— $n = \deg C$ ;
— $r = \deg\bigl(\mathrm{Tan}(\Gamma)\bigr)$.

Classically, developable surfaces were considered as systems packaging together three 1-dimensional families in $\mathbf{P}^3$ of planes, points, and lines respectively; these are the planes parametrized by points of $C$, the points of $\Gamma$, and the tangent lines of $\Gamma$. In this perspective, $m$ is the number of points of the system contained in a general plane of $\mathbf{P}^3$, $n$ is the number of planes of the system passing through a general point of $\mathbf{P}^3$, and $r$ is the number of lines of the system meeting a general line of $\mathbf{P}^3$. The integers $m, n, r$ are called respectively the order, class, and rank of the system.

**(B.7)** Knowing the family of planes given by $C \subseteq \check{\mathbf{P}}^3$, one may find the generators of the developable surface $\mathrm{Tan}(\Gamma) \subseteq \mathbf{P}^3$ by taking the intersections in $\mathbf{P}^3$ of pairs of infinitely near members of $C$.

Indeed, the generators of $\mathrm{Tan}(\Gamma)$ are the orthogonals of the tangent lines of $C$. Now for $p, q \in C$, the orthogonal of the line $\langle p, q \rangle \subseteq \check{\mathbf{P}}^3$ is the intersection of the two planes $p^\perp$ and $q^\perp$ in $\mathbf{P}^3$. Letting $q$ tend to $p$ yields the asserted statement.

**(B.8)** Let $o$ and $H$ be a general point and a general plane in $\mathbf{P}^3$. The duality between projections and hyperplanes sections (see Proposition (1.5)) tells us that the projection of $\Gamma$ from $o$ is dual as a plane curve to the hyperplane section of $\Gamma^\vee = \mathrm{Tan}(C)$ by $o^\perp$ and, dually, that the section of $S = \mathrm{Tan}(\Gamma)$ by $H$ is dual to the projection of $C$ from the point $H^\perp$. This may be established using the following elementary statement, see Proposition (B.11).

**(B.9) Lemma.** *Let $o \in \mathbf{P}^3$ be a point off $\Gamma$, and consider the projection $\pi_o : \mathbf{P}^3 \dashrightarrow \mathbf{P}^3/o$. Let $p \in \Gamma$ be a smooth point such that the line $\langle o, p \rangle$ intersects $\Gamma$ only at $p$.*
*(i) If the line $\langle o, p \rangle$ is tangent to $\Gamma$ at $p$, then the projected curve $\pi_o(\Gamma)$ has a cusp at $\pi_o(p)$.*
*(ii) If the line $\langle o, p \rangle$ is not tangent to $\Gamma$ at $p$, there is a unique plane $H$ containing $o$ and tangent to $\Gamma$ at $p$, and its projection $\pi_o(H)$ is the tangent line to $\pi_o(\Gamma)$ at $\pi_o(p)$. If $H$ osculates $\Gamma$ to the order at least $3$, then $\pi_o(H)$ is a flex tangent.*

This lemma is pretty clear from the differential point of view, given that the kernel of the differential of $\pi_o$ at the point $p$ is the tangent direction of the line $\langle o, p \rangle$.

Another interesting interpretation of this result is in terms of intersection multiplicities. Call $\Gamma' = \pi_o(\Gamma)$, $p' = \pi_o(p)$, and $H' = \pi_o(H)$, assuming that $H$ is a hyperplane in $\mathbf{P}^3$ containing $p$, so that $H'$ is a line in $\mathbf{P}^2$. Then, for local intersection multiplicities at $p$ and $p'$, the projection formula tells us that
$$(H \cdot \Gamma)_p = (H' \cdot \Gamma')_{p'}.$$

Thus, if $H$ is an ordinary (resp. osculating) tangent plane of $\Gamma$ at $p$, then $H'$ is an ordinary (resp. flex) tangent of $\Gamma'$ at $p'$. Moreover, if the line $\langle o, p \rangle$ is tangent to $\Gamma$ at $p$, then for all planes $H$ passing through $o$ and $p$ one has $(H \cdot \Gamma)_p > 1$, hence $(H' \cdot \Gamma')_{p'} > 1$ for all lines passing through $p'$: this implies that $\Gamma'$ is singular at $p'$, and it has a cusp there because it is unibranch, as follows from our assumption that $\langle o, p \rangle$ intersects $\Gamma$ only at $p$.

Finally, the lemma may also be proved by local computations. Let me now give a sample of these. Assume $p = (1 : 0 : 0 : 0)$ and $\Gamma$ is parametrized locally at $p$ by $t \mapsto (1 : t : t^2 : t^{3+k})$ for some non-negative integer $k$. Let $o = (a : b : c : d)$, and assume that $b \neq 0$ (the other possibilities may be considered analogously). Then the projection from $o$ writes
$$(x : y : z : w) \mapsto (x - \tfrac{a}{b} y : z - \tfrac{c}{b} y : w - \tfrac{d}{b} y),$$
so that $\Gamma' = \pi_o(\Gamma)$ is locally parametrized at $p' = \pi_o(p)$ by
$$t \mapsto \left(1 : \frac{t^2 - \tfrac{c}{b} t}{1 - \tfrac{a}{b} t} : \frac{t^{3+k} - \tfrac{d}{b} t}{1 - \tfrac{a}{b} t}\right).$$

If $d = 0$, then $o$ sits on the osculating plane of $\Gamma$ at $p$. If moreover $c \neq 0$, then $\Gamma'$ has a local parametrization at $p'$ of the form $t \mapsto (1 : -\tfrac{c}{b} t + \sum_{m \geqslant 2} \alpha_m t^m : t^{3+k} + \sum_{m \geqslant 4+k} \beta_m t^m)$, hence it has a flex at $p'$. If $d = c = 0$, then $o$ sits on the tangent line to $\Gamma$ at $p$, and the local parametrization shows that $\Gamma'$ has cusp at $p'$: an ordinary cusp (i.e., an $A_2$ double point) if $k = 0$, a ramphoid cusp (i.e., an $A_4$ double point) if $k = 1$, and so on.

**(B.10)** If we project a smooth curve $\Gamma \subseteq \mathbf{P}^3$ from a general point $o$, its image is a plane nodal curve $\Gamma'$, and the nodes correspond to the finitely many secant lines of $\Gamma$ (i.e., lines of the form $\langle p, q \rangle$ with $p, q \in \Gamma$) passing through $o$.[8] The projection $\Gamma'$ has moreover finitely many flexes (in addition to the original flexes of $\Gamma$), corresponding to the finitely many osculating hyperplanes of $\Gamma$ passing through $o$ (i.e., to the finitely intersection points of the oscdual curve $C = \text{Tan}(\Gamma)^\vee$ with the hyperplane $o^\perp$).

Assume now that we project from a general point $o$ of the tangential surface $\text{Tan}(\Gamma)$. Then the projection $\Gamma'$ has one cusp, corresponding to the unique generator of $\text{Tan}(\Gamma)$ passing through $o$, and one less node than in the general situation, because the geometric genus of $\Gamma'$ must always be the same as that of $\Gamma$ (another interpretation of this is that one of the secant lines of $\Gamma$ in the general situation has become a tangent line in the present situation). This implies in particular that the class of the plane curve $\Gamma'$ in the present situation is one less than in the general situation, by the Plücker formula $(P_1)$, see (1.12).

---

[8]see [22, IV, Theorem 3.10] for a complete proof, and [20] for a vast generalization.

Moreover, the projection $\Gamma'$ from a point $o \in \operatorname{Tan}(\Gamma)$ has two less flexes than the projection from a general point. This follows from the Plücker formula $(P_2)$, arguing once again that a node has become a cusp. From a more geometric point of view, we have seen in (B.7) that the generators of $\operatorname{Tan}(\Gamma)$ are the intersections of two infinitely near osculating planes of $\Gamma$; these are the two osculating planes of $\Gamma$ through $o$ that no longer give flexes of $\Gamma'$. Note that these two infinitely near osculating planes of $\Gamma$ correspond to the tangency point of the oscdual curve $C = \operatorname{Tan}(\Gamma)^\vee$ with the plane $o^\perp$.

**(B.11) Proposition.** *The section of* $\operatorname{Tan}(\Gamma)$ *by a plane* $H$ *is the dual curve of the plane curve* $\pi_{H^\perp}(C)$.

*Proof.* For all $p \in \operatorname{Tan}(\Gamma) \cap H$, there is a point $q \in \Gamma$ such that $p = H \cap \mathbf{T}_q\Gamma$. By duality, the orthogonal line $(\mathbf{T}_q\Gamma)^\perp$ is the tangent to $C$ at the point $(\operatorname{Osc}_q\Gamma)^\perp$. It follows that the plane $p^\perp$ is tangent to $C$ at $(\operatorname{Osc}_q\Gamma)^\perp$, as it contains $(\mathbf{T}_q\Gamma)^\perp$; moreover, as it passes through the point $H^\perp$, its projection from $H^\perp$ is a tangent line of $\pi_{H^\perp}(C)$. Conversely, all tangent lines to $\pi_{H^\perp}(C)$ are obtained in this way. □

As an illustration of the proposition, note that the plane curve $H \cap \operatorname{Tan}(\Gamma)$ has cusps at the points where $H$ meets the cuspidal edge of $\operatorname{Tan}(\Gamma)$, i.e., $\Gamma$ itself, while the dual plane curve $\pi_{H^\perp}(C)$ correspondingly has flex points given by the osculating planes to $C$ passing through $H^\perp$, i.e., the planes $q^\perp$ for $q \in H \cap \Gamma$.

**(B.12) Proposition.** *The two tangential surfaces* $\operatorname{Tan}(\Gamma)$ *and* $\operatorname{Tan}(C)$ *have the same degree.*

*Proof.* As we have seen in (B.6), the degree of $\operatorname{Tan}(\Gamma)$ is the number of lines $\mathbf{T}_p\Gamma$ meeting a given, general, line $\Lambda \subseteq \mathbf{P}^3$. Since the orthogonal $(\mathbf{T}_p\Gamma)^\perp$ is the tangent line to $C$ at the point $(\operatorname{Osc}_p\Gamma)^\perp$, this is also the number of lines $\mathbf{T}_{p'}C$ meeting the line $\Lambda^\perp$, i.e., the degree of $\operatorname{Tan}(C)$. □

*(B.12.1) Remark.* Proposition (B.12) is a twisted version of the fact that the class of a (non-developable) ruled surface equals its degree.

To compare the two situations, let's recall the proof. Let $S \subseteq \mathbf{P}^3$ be a ruled surface, and assume $\check{S}$ is a surface. The class of $S$ is the number of tangent planes of $S$ containing a given, general, line $\Lambda \subseteq \mathbf{P}^3$. Since $\check{S}$ is a surface, the planes tangent to $S$ are exactly the planes that contain a ruling of $S$. It follows that the class of $S$ is the number of rulings of $S$ that meet $\Lambda$. But the latter number also equals the degree of $S$.

If $S$ is developable, it is no longer true that its tangent planes are those containing a ruling (indeed there is only one tangent plane per ruling), and the number of planes tangent to $S$ containing a line $\Lambda$ is zero, as $S^\vee$ actually has codimension 2 in $\check{\mathbf{P}}^3$.

It is useful to analyze the sections of the developable surface by planes containing a generator.

**(B.13) Lemma.** *Let the hyperplane* $H \subseteq \mathbf{P}^3$ *contain the line* $\mathbf{T}_p\Gamma$ *for some* $p \in \Gamma$, *and write* $H \cap \operatorname{Tan}(\Gamma)$ *as* $\mathbf{T}_p\Gamma + R$. *In general, the residual section* $R$ *has one more flex and two less cusps than the general plane section of* $\operatorname{Tan}(\Gamma)$. *The extra flex point is* $p$, *and* $\mathbf{T}_pR = \mathbf{T}_p\Gamma$.

*Proof.* The point $H^\perp \in \check{\mathbf{P}}^3$ sits on the line $(\mathbf{T}_p\Gamma)^\perp$ which is tangent to $C$ at the point $(\operatorname{Osc}_p\Gamma)^\perp$, so the curve $\pi_{H^\perp}(C)$ has one more cusp and two less flexes than the projection of $C$ from a general point, as explained in Paragraph (B.10). Dually, this gives the first assertion about $R$. In addition, note that the class of $\pi_{H^\perp}(C)$ is one less than that of a general projection, agreeing with the fact that $\deg(R) = \deg(\operatorname{Tan}\Gamma) - 1$.

The tangent cone to $\pi_{H^\perp}(C)$ at its extra cusp is the projection from $H^\perp$ of the plane $\mathrm{Osc}_{(\mathrm{Osc}_p\Gamma)^\perp}C$; its orthogonal is the extra flex point of $R$, and it is the point $p$. Moreover, the tangent line to $R$ at $p$ is the orthogonal of the cuspidal point of $\pi_{H^\perp}(C)$; the latter is the projection from $H^\perp$ of the line $(\mathbf{T}_p\Gamma)^\perp$, hence $\mathbf{T}_p R$ equals $\mathbf{T}_p\Gamma$. □

We may now prove the following. Before we start, note that being $\mathrm{Tan}(\Gamma)$ ruled, its double points lie on two generators: if they are distinct then the double point is on the ordinary double curve, otherwise it is on the cuspidal edge. Given a generator $\mathbf{T}_p\Gamma$, $p$ is the only point of $\mathbf{T}_p\Gamma$ lying on the cuspidal edge, unless $\mathbf{T}_p\Gamma$ is tangent to $\Gamma$ at some other point.

**(B.14) Proposition.** *The degree of the ordinary double curve of* $\mathrm{Tan}(\Gamma)$ *equals* $r - 4$.

*Proof.* Consider the section of $\mathrm{Tan}(\Gamma)$ by a plane $H$ containing a generator $\mathbf{T}_p\Gamma$ as in Lemma (B.13): it consists of the line line $\mathbf{T}_p\Gamma$ and a residual curve $R$ of degree $r - 1$. They meet with multiplicity 3 at $p$, and at $r - 4$ further points $p_1, \ldots, p_{r-4}$, which are the intersection points of $\mathbf{T}_p\Gamma$ and the ordinary double curve of $\mathrm{Tan}(\Gamma)$. □

Note that, in the notation of the above proof, since $p_1, \ldots, p_{r-4}$ are the intersection points of $\mathbf{T}_p\Gamma$ and the ordinary double curve, they are independent of the choice of $H$ in the pencil of planes containing $\mathbf{T}_p\Gamma$.

For all $i = 1, \ldots, r - 4$, there is a point $q_i \in \Gamma$ such that $p_i = \mathbf{T}_p\Gamma \cap \mathbf{T}_{q_i}\Gamma$. It follows that the plane $p_i^\perp \subseteq \check{\mathbf{P}}^3$ is tangent to $C$ at the two points $(\mathrm{Osc}_p\Gamma)^\perp$ and $(\mathrm{Osc}_{q_i}\Gamma)^\perp$, which explains why it should be an ordinary double point of the surface $\mathrm{Tan}(\Gamma) = C^\vee$ according to the general principles of duality.

**(B.15) Lemma.** *Let* $p \in \Gamma$. *The section of* $\mathrm{Tan}(\Gamma)$ *by its tangent plane* $\mathrm{Osc}_p(\Gamma)$ *decomposes as* $2\mathbf{T}_p\Gamma + R$, *where the curve* $R$ *is tangent to* $\mathbf{T}_p\Gamma$ *at* $p$, *and passes through all the ordinary double points* $p_1, \ldots, p_{r-4}$ *of* $Tan(\Gamma)$ *lying on* $\mathbf{T}_p\Gamma$.

**(B.16) Remark.** The pencil of curves cut out on $\mathrm{Tan}(\Gamma)$ by the planes containing $\Lambda = \mathbf{T}_p\Gamma$ is a family of plane curves with prescribed contact $3p + p_1 + \cdots + p_{r-4}$ with $\Lambda$, which degenerates to a curve containing $\Lambda$ itself when the plane tends to the tangent plane $\mathrm{Osc}_p(\Gamma)$. This is exactly the situation considered by Caporaso and Harris in [8] and studied in the first half of this volume.

*Proof of Lemma (B.15).* Let $H = \mathrm{Osc}_p\Gamma$. The section $H \cap \mathrm{Tan}(\Gamma)$ must decompose as $2\mathbf{T}_x\Gamma + R$, with $R$ the plane dual of the projected curve $\pi_{H^\perp}(C)$. This time the point $H^\perp \in \check{\mathbf{P}}^3$ sits on the curve $C$, and the projection $\pi_{H^\perp}(C)$ is smooth at the point image of $H^\perp$, namely $\pi_{H^\perp}(\mathbf{T}_{H^\perp}C)$, with tangent line $\pi_{H^\perp}(\mathrm{Osc}_{H^\perp}C)$ there. It follows that $p = (\mathrm{Osc}_{H^\perp}C)^\perp$ is on the plane dual $R$, with tangent $\mathbf{T}_p R = (\mathbf{T}_{H^\perp}C)^\perp = \mathbf{T}_p\Gamma$. This proves the first assertion.

For the second, letting once again $p_i = \mathbf{T}_p\Gamma \cap \mathbf{T}_{q_i}\Gamma$, we have seen that the plane $p_i^\perp$ is tangent to $C$ at the two points $H^\perp$ and $(\mathrm{Osc}_{q_i}\Gamma)^\perp$. Thus $\pi_{H^\perp}(p_i^\perp)$ is the tangent line to $\pi_{H^\perp}(C)$ at $\pi_{H^\perp}((\mathrm{Osc}_{q_i}\Gamma)^\perp)$, hence $p_i \in R$ since the latter is the plane dual of $\pi_{H^\perp}(C)$. □

Finally, let us check by duality that the degree of $R$ is indeed $r - 2$ in the above context. Let $g$ be the genus of $C$. The projection of $C$ from a general point of $\check{\mathbf{P}}^3$ is a nodal degree $n$ plane curve of genus $g$, which thus has $\delta = \frac{1}{2}(n-1)(n-2) - g$ nodes, hence class $r = n(n-1) - 2\delta$. If however it gets projected from a general point on itself, it still gives a genus $g$ nodal curve, but with degree $n - 1$, so that the number of nodes becomes $\delta' = \frac{1}{2}(n-2)(n-3) - g = \delta - (n-2)$. The class of this curve is then $(n-1)(n-2) - 2\delta'$, which is $r - 2$ as required.

## B.3 – Plücker formulae

We now explain how the Plücker formulae for plane curves may be used to derive relations between the numerical characters of a developable surface, following Cayley [9]. As for plane curves, it turns out to be appropriate to let the curves $\Gamma$ and $C$ have various kinds of singularities which are stable under oscduality. We thus consider a slightly more general situation than what we did in Subsection B.2, by allowing the stationary behaviours indicated below.

**(B.17) Characteristic numbers, II.** In addition to the quantities introduced in (B.6), we let
— $\alpha$ the number of ordinary cusps of $C$, i.e., points at which $C$ has a local parametrization of the form $t \mapsto (t^2, t^3, t^4)$ ;
— $\beta$ the number of ordinary cusps of $\Gamma$ ;
— $x$ the degree of the ordinary double curve of $\mathrm{Tan}(\Gamma)$;
— $y$ the degree of the ordinary double curve of $\mathrm{Tan}(C)$;
— $g$ the number of apparent double points of $C$;
— $h$ the number of apparent double points of $\Gamma$.

The quantities $\alpha$ and $\beta$ are classically referred to as the number of stationary planes and points of the system. The quantity $x$ is the number of points in a general plane $H \subseteq \mathbf{P}^3$ that lie on two distinct tangent lines of $\Gamma$, see Proposition (B.14) and the comment before it; dually, $y$ is the number of planes passing through a general point, and containing two distinct tangent lines of $\Gamma$. The quantity $h$ is the number of lines meeting $\Gamma$ in two points (i.e., secants of $\Gamma$) passing through a general point of $\mathbf{P}^3$; dually, $g$ is the number of lines lying in a general plane $H \subseteq \mathbf{P}^3$ that are contained in two osculating planes of $\Gamma$.

**(B.18) Lemma.** *Consider the osculatory Gauss map $\Gamma \to C$, mapping a general point $p \in \Gamma$ to the point $(\mathrm{Osc}_p \Gamma)^\perp \in \check{\mathbf{P}}^3$. It sends a point where $\Gamma$ has a local parametrization of the form $t \mapsto (t, t^2, t^4)$ (resp. $t \mapsto (t^2, t^3, t^4)$) to a point where $C$ has a local parametrization of the form $t \mapsto (t^2, t^3, t^4)$ (resp. $t \mapsto (t, t^2, t^4)$).*

The upshot is that the osculatory Gauss map turns points of types $\alpha$ and $\beta$ of $\Gamma$ to points of types $\beta$ and $\alpha$ of $C$, respectively.

*Proof.* This is a local computation: to a local branch of $\Gamma$ parametrized by $t \in \mathbf{C} \mapsto f(t) \in \mathbf{C}^4$, the osculatory Gauss maps associates the local branch of $C$ parametrized by $t \in \mathbf{C} \mapsto (f \wedge \partial_t f \wedge \partial_t^2 f)(t) \in \bigwedge^3 \mathbf{C}^4 \cong \check{\mathbf{C}}^4$. In the former case, the result follows from the identity

$$\begin{pmatrix} 1 \\ t \\ t^2 \\ t^4 \end{pmatrix} \wedge \begin{pmatrix} 0 \\ 1 \\ 2t \\ 4t^3 \end{pmatrix} \wedge \begin{pmatrix} 0 \\ 0 \\ 2 \\ 12t^2 \end{pmatrix} = \begin{pmatrix} -6t^4 & 16t^3 & -12t^2 & 2 \end{pmatrix}.$$

We leave the computation in the latter case to the reader. □

**(B.19)** Let us apply the Plücker formulae to a general hyperplane section $H \cap \mathrm{Tan}(\Gamma)$; as a plane curve, it is dual to $\pi_{H^\perp}(C)$. The former is a curve of degree $r$, with $x$ nodes and $m$ cusps, and the latter has degree $n$ and $g$ nodes; for general $H$, its only cusps come from cusps of $C$, which are in the number $\alpha$. We thus get the following identities:

$$n = r(r-1) - 2x - 3m; \qquad r = n(n-1) - 2g - 3\alpha;$$
$$\alpha = 3r(r-2) - 6x - 8m; \qquad m = 3n(n-2) - 6g - 8\alpha,$$

hence

$$\alpha - m = 3(n - r); \qquad 2(g - x) = (n - r)(n + r - 9).$$

**(B.20)** Dually, one may apply the Plücker formulae to the pair of dual plane curves $\pi_o(\Gamma)$ and $o^\perp \cap \operatorname{Tan}(C)$ for a general point $o \in \mathbf{P}^3$. This gives:

$$r = m(m-1) - 2h - 3\beta; \qquad m = r(r-1) - 2y - 3n;$$
$$n = 3m(m-2) - 6h - 8\beta; \qquad \beta = 3r(r-2) - 6y - 8n,$$

hence

$$n - \beta = 3(r - m) \qquad 2(y - h) = (r - m)(r + m - 9).$$

**(B.21)** Eventually, one may combine the identities from (B.19) and (B.20) to obtain:

$$\alpha - \beta = 2(n - m), \qquad x - y = n - m,$$
$$2(g - h) = (n - m)(n + m - 7).$$

**(B.22)** It seems natural to allow the system to have stationary lines as well, i.e., the curve $\Gamma$ to have points at which it is locally parametrized by $t \mapsto (t, t^3, t^4)$. I have chosen not to do it from the start, in order to give a more straightforward application of the Plücker formulae for plane curves. Let me now add two statements with which it is possible to adapt the above arguments in order take the possibility of stationary lines in consideration, which I will leave to the reader.

First, a local computation as in the proof of (B.18) shows that the osculatory Gauss map sends a point at which $\Gamma$ has a stationary tangent line to a point at which $C$ has a stationary tangent line. The second statement is the following.

**(B.23) Lemma.** *The stationary lines of the system give cuspidal lines to both tangential surfaces* $\operatorname{Tan}(\Gamma)$ *and* $\operatorname{Tan}(C)$.

*Proof.* Let $p \in \Gamma$ be a point at which $\Gamma$ is locally parametrized by $t \mapsto (t, t^3, t^4)$, and $\Lambda = \mathbf{T}_p\Gamma$ be the corresponding (stationary) generator. We shall prove that for a general plane $H \subseteq \mathbf{P}^3$, the section $H \cap \operatorname{Tan}(\Gamma)$ has a cusp at the point $H \cap \Lambda$, or equivalently the curve $\pi_{H^\perp}(C)$ has a flex at $\pi_{H^\perp}\big((\operatorname{Osc}_p\Gamma)^\perp\big)$ along the line $\pi_{H^\perp}(\Lambda^\perp)$.

As we have seen, the curve $C$ has a local parametrization of the form $t \mapsto (t, t^3, t^4)$ at the point $(\operatorname{Osc}_p\Gamma)^\perp$. Thus after projection from the general point $H^\perp$, the curve $\pi_{H^\perp}(C)$ has a local parametrization of the form $t \mapsto (t, t^3)$. The result follows since $\Lambda^\perp$ is the tangent line to $C$ at $(\operatorname{Osc}_p\Gamma)^\perp$. □

## B.4 – Generalized Plücker formulae

In this section I outline the generalization of the previous formulae to curves in a projective space of arbitrary dimension, which is due to Piene [36]. This is not used in the course of the text, but it would be a pity not to mention it. Note that there also exists a duality between the $i$-th ranks of an arbitrary variety in projective space and those of its dual, which is due to Urabe [50], see also [26], and [2], but this I will be reasonable and not talk about it here.

**(B.24)** So far we have only considered extrinsic invariants of the curve $\Gamma \subseteq \mathbf{P}^3$. Adding its geometric genus $p_g$ in the picture, one finds:

(B.24.1) $$r = 2m + 2p_g - 2 - \beta$$

(recall the notation from (B.6) and (B.17)).

*Proof.* Let $\Lambda \subseteq \mathbf{P}^3$ be a general line. We want to compute the number of tangent lines to $\Gamma$ meeting $\Lambda$, or equivalently the number of planes tangent to $\Gamma$ and containing $\Lambda$. To this end we consider the projection from $\Lambda$; it induces a degree $m$ covering $\bar\Gamma \to \mathbf{P}^1$, where $\bar\Gamma$ is the normalization of $\Gamma$. By the Riemann–Hurwitz formula, it has $2m + 2p_g - 2$ ramification points, of which $\beta$ come from the ramification points of the normalization map of $\Gamma$. The other ones correspond to planes tangent to $\Gamma$ and containing $\Lambda$. □

*(B.24.2) Remark.* The same reasoning may be applied to the curve $C$; one gets in this way the identity
$$r = 2n + 2p_g - 2 - \alpha.$$
Indeed, the osculatory Gauss map gives a birationality between $\Gamma$ and $C$, so these two curves have the same geometric genus, and $r$ is the degree of both $\mathrm{Tan}(C)$ and $\mathrm{Tan}(\Gamma)$ by Proposition (B.12).

Putting the two identities together, one finds the formula $2(n-m) = \alpha - \beta$ which we have already gotten in (B.21).

**(B.25) New setup.** From now on, we consider a curve $\Gamma \subseteq \mathbf{P}^N$ of degree $m$ and geometric genus $p_g$. Let $\bar\Gamma \to \Gamma$ be its normalization, and $f : \bar\Gamma \to \mathbf{P}^N$ its composition with $\Gamma \hookrightarrow \mathbf{P}^N$. For all $p \in \bar\Gamma$ one may choose a local parameter $t$ at $p$ on $\bar\Gamma$ and a system of homogeneous coordinates $(x_0 : x_1 : \ldots : x_N)$ on $\mathbf{P}^N$ such that $f : \bar\Gamma \to \mathbf{P}^N$ is given locally by
$$t \mapsto \left(1 : t^{1+s_1} : \ldots : t^{N+s_1+\cdots+s_N}\right)$$
up to higher order terms (i.e., the $i$-th coordinate of $f$ is a formal power series in $t$ whose leading term is $t^{i+s_1+\cdots+s_i}$), and $s_1, \ldots, s_N$ are non-negative integers depending only on $p$. The point $p$ is *ordinary* if $s_1 = \cdots = s_N = 0$, and *stationary* or *hyperosculating* otherwise.

For all $i = 0, \ldots, N-1$, the $i$-th osculating subspace of $\Gamma$ at $p \in \bar\Gamma$, denoted $\mathrm{Osc}_p^i(\Gamma)$, is the $(N-i)$-plane in $\mathbf{P}^N$ defined by the equations $x_N = \cdots = x_{N-i+1} = 0$ in a system of homogeneous coordinates as above. We shall consider the $i$-th *osculating variety* of $\Gamma$, denoted by $\mathrm{Osc}^i(\Gamma)$, that is the union of the $i$-th osculating subspaces of $\Gamma$ at all the points of $\bar\Gamma$. For $i = 0$ this is $\Gamma$ itself, for $i = 1$ this is the tangential surface of $\Gamma$ that we have considered in the previous subsections; for $i = N-1$ this fills the whole $\mathbf{P}^N$, hence it is best seen as a curve $C \subseteq \check{\mathbf{P}}^N$ which I call the *osc-dual* curve of $\Gamma$.

**(B.26) The $i$-th ranks.** For all $i = 0, \ldots, N-1$, let $r_i$ be the degree of the $i$-th osculating variety of $\Gamma$; it is called the *$i$-th rank* of $\Gamma$. For $i = N-1$, this is by definition the degree of the osc-dual curve $C \subseteq \check{\mathbf{P}}^N$. We define integers $k_i$ by
$$k_i := \sum\nolimits_{p \in \Gamma} s_i(p),$$
where $s_i(p)$ is defined as in (B.25).

A standard computation on the Chern classes of the bundles of principal parts of the embedding of $\bar{\Gamma}$ in $\mathbf{P}^N$, see [13, Prop. 10.4.13], gives the following formulae:

(B.26.1) $\quad r_i = (i+1)\big(m + i(p_g - 1)\big) - \sum_{j=1}^{i} \big((i-j+1)k_j\big) \quad$ for all $i = 0, \ldots, N-1$

(B.26.2) $\quad \sum_{j=1}^{N} \big((N-j+1)k_j\big) = (N+1)\big(m + N(p_g - 1)\big).$

Formula (B.26.1) for $i = 0$ is trivial, as $r_0 = \deg(\Gamma) = m$ by definition. Formula (B.26.2) has to be understood as the incarnation for $i = N$ of formula (B.26.1), setting $r_N = 0$. For $i = 1$, one gets a formula which is exactly (B.24.1). Eventually, note that formula (B.26.2) may be interpreted as counting the stationary points of $\Gamma$, with the multiplicity

$$s_N + 2s_{N-1} + \cdots + (N-1)s_2 + Ns_1.$$

**(B.27) Generalized Plücker formulae.** As a straightforward corollary of the formulas in (B.26), one has

(B.27.1) $\qquad r_{i-1} - 2r_i + r_{i+1} = 2p_g - 2 - k_{i+1}, \qquad$ for all $i = 0, \ldots, N-1,$

with the convention that $r_N = r_{-1} = 0$.

*(B.27.2) Plane curves.* For $N = 2$, we are considering a plane curve $\Gamma$ of degree $m$. The number $k_1$ (resp. $k_2$) counts, with multiplicities, the cusps (resp. the flexes) of $\Gamma$, so we will call it $\kappa$ (resp. $\iota$) as usual. Formula (B.27.1) gives the two identities

$$-2m + \check{m} = 2p_g - 2 - \kappa \, ; \qquad m - 2\check{m} = 2p_g - 2 - \iota.$$

The second formula is merely the first applied to $C$, which in this case is the plain dual of $\Gamma$. Assuming $\Gamma$ has only $\delta$ nodes and $\kappa$ (ordinary) cusps as singularities, these may be recovered from the usual formulae

$$\check{m} = m(m-1) - 2\delta - 3\kappa \, ; \qquad 2p_g - 2 = m(m-3) - 2\delta - 2\kappa$$

(the first one is the computation of the class using polars for a singular curve as in (1.12), and the second one is the adjunction formula for a (possibly) singular plane curve. One gets the first formula by writing

$$\check{m} - 2m = m(m-1) - 2\delta - 3\kappa - 2m = m(m-3) - 2\delta - 3\kappa = 2p_g - 2 - \kappa.$$

*(B.27.3) Space curves.* For $N = 3$ we are considering a space curve $\Gamma \subseteq \mathbf{P}^3$ as in the previous subsections. In the notation therein, the ranks $r_0, r_1, r_2$ are respectively the degree $m$, the rank $r$, and the class $n$. On the other hand, $k_1, k_2, k_3$ are respectively the number $\beta$ of cusps of $\Gamma$, the number of stationary lines (let us call it $\gamma$), and the number $\alpha$ of cusps of $C$. Formula (B.27.1) gives the three identities

$$-2m + r = 2p_g - 2 - \beta \, ; \qquad m - 2r + n = 2p_g - 2 - \gamma \, ; \qquad r - 2n = 2p_g - 2 - \alpha.$$

**(B.28) Piene duality** [36, Thm. 5.1]**.** *Let $\Gamma \subseteq \mathbf{P}^N$ be a curve, and $C \subseteq \check{\mathbf{P}}^N$ its osc-dual curve. The osc-dual curve of $C$ is $\Gamma$, and the following relations hold for all $i = 0, \ldots, N-1$:*

(B.28.1) $\qquad\qquad\qquad\qquad r_i(C) = r_{N-i-1}(\Gamma)$

(B.28.2) $\qquad\qquad\qquad\qquad k_i(C) = k_{N-i-1}(\Gamma).$

Formula (B.28.1) generalizes the fact that for $N = 3$ the two tangential surfaces $\mathrm{Tan}(\Gamma)$ and $\mathrm{Tan}(C)$ have the same degree, see Proposition (B.12), and Formula (B.28.2) the fact that the osc-Gauss map turns stationary planes, lines and points of $\Gamma$ into stationary points, lines, and planes of $C$ respectively, see Lemma (B.18).

**(B.29) Interpretation as strata of the dual hypersurface.** Underlying the merely numerical relations of (B.28) is the following geometric duality, which generalizes (B.3).

For all $i > 0$, we say that a hyperplane $H \subseteq \mathbf{P}^N$ is *$i$-osculating* to $\Gamma$ at $p \in \bar{\Gamma}$ if it contains $\mathrm{Osc}_p^i(\Gamma)$. If $p \in \bar{\Gamma}$ is an ordinary point of $\Gamma$, this is equivalent to the divisor $f^*H - (i+1)p$ on $\bar{\Gamma}$ being effective. We let $\Gamma^{\vee i}$ be the closed subset parametrizing hyperplanes $i$-osculating to $\Gamma$ at some point on $\bar{\Gamma}$; thus $\Gamma^{\vee 1}$ is the plain dual hypersurface $\Gamma^\vee$, and $\Gamma^{\vee N-1}$ is the osc-dual curve $C$.

For all $i = 0, \ldots, N-1$, one has $\Gamma^{\vee i} = \mathrm{Osc}^{N-1-i}(C)$, and moreover $\mathrm{Osc}^i(\Gamma)^\vee = \Gamma^{\vee i+1}$. In fact, for all $p \in \bar{\Gamma}$, the linear space $\mathrm{Osc}_p^i(\Gamma)^\perp$ is the $(N-1-i)$-th osculating space of $C$ at $p$ (note that the map $p \in \bar{\Gamma} \mapsto [\mathrm{Osc}_p^{N-1}(\Gamma)] \in \check{\mathbf{P}}^N$ is the normalization of $C$).

For all $i = 0, \ldots, N-1$, let $f_i : \bar{\Gamma} \to \mathbf{G}(\mathbf{P}^N, i)$ and $f^{\vee i} : \bar{\Gamma} \to \mathbf{G}(\check{\mathbf{P}}^N, i))$ be the maps $p \mapsto \mathrm{Osc}_p^i(\Gamma)$ and $p \mapsto \mathrm{Osc}_p^i(C)$ respectively. The map $f^{\vee N-1-i}$ equals the composition of $f_i$ with the orthogonality map $\Lambda \in \mathbf{G}(\mathbf{P}^N, i) \mapsto \Lambda^\perp \in \mathbf{G}(\check{\mathbf{P}}^N, N-1-i)$. It follows that the two curves $f_i(\bar{\Gamma})$ and $f^{\vee N-1-i}(\bar{\Gamma})$ have the same degree under the Plücker embeddings of $\mathbf{G}(\mathbf{P}^N, i)$ and $\mathbf{G}(\check{\mathbf{P}}^N, N-1-i)$ respectively; these degrees are $r_i(\Gamma)$ and $r_{N-1-i}(C)$ respectively (see, e.g., [21, Example 19.11]), so we find Formula (B.28.1).

## B.5 – De Jonquières' formula

The next result is even more general than the previous ones. See [4, Chap. VIII] for a proof and pointers to the literature.

**(B.30) Theorem** (de Jonquières formula)**.** *Let $m$ be a positive integer, and $m_1, \ldots, m_a$ be non–negative integers such that $m = \sum_{1 \leqslant s \leqslant a} s \cdot m_s$. The virtual number of divisors having $m_s$ points of multiplicity $s$ for all $s = 1, \ldots, a$ in a given linear series of degree $m$ and dimension $i = \sum_{1 \leqslant s \leqslant a}(s-1) \cdot m_s$ on a smooth, genus $g$ curve is the coefficient of the $t_1^{m_1} \cdots t_a^{m_a}$ term in the formal power series*

$$\left(1 + 1^2 t_1 + 2^2 t_2 + \cdots + a^2 t_a\right)^g \cdot \left(1 + 1\,t_1 + 2t_2 + \cdots + at_a\right)^{m-i-g}.$$

**(B.31) Link with the $i$-th ranks.** We now show how the de Jonquières formula may be used to derive the $i$-th ranks of a curve $\Gamma \subseteq \mathbf{P}^N$ of degree $m$ and geometric genus $g$ as in (B.25), thus recovering formula (B.26.1).

By definition, $r_i(\Gamma) = r_{N-1-i}(C) = \deg(\mathrm{Osc}^i(\Gamma))$ is the number of hyperplanes $i$-osculating to $\Gamma$ in a general $i$-dimensional linear subsystem of $|\mathcal{O}_\Gamma(1)|$. By de Jonquières' formula, the number of hyperplanes $H$ such that $f^*H - (i+1)p$ is effective for some $p \in \bar{\Gamma}$ in a general

$i$-dimensional linear subsystem of $|\mathcal{O}_\Gamma(1)|$ is the coefficient of the monomial $t_1^{m-i-1}t_{i+1}$ in

$$\begin{aligned}\left(1 + 1^2 t_1 + (i+1)^2 t_{i+1}\right)^g &\cdot \left(1 + 1\, t_1 + (i+1) t_{i+1}\right)^{m-i-g}\\ &= \left(1 + 1^2 t_1 + [(i+1) + i(i+1)] t_{i+1}\right)^g \cdot \left(1 + 1\, t_1 + (i+1) t_{i+1}\right)^{m-i-g}\\ &= \sum_{j=0}^{g} \binom{g}{j} \left(1 + t_1 + (i+1) t_{i+1}\right)^{(g-j)+m-i-g} \cdot \left(i(i+1) t_{i+1}\right)^j,\end{aligned}$$

which is also the coefficient of $t_1^{m-i-1}t_{i+1}$ in

$$\left(1 + t_1 + (i+1) t_{i+1}\right)^{m-i} + i(i+1) t_{i+1} \cdot g \left(1 + t_1 + (i+1) t_{i+1}\right)^{m-i-1},$$

that is

(B.31.1) $\qquad (m-i)(i+1) + i(i+1) \cdot g = (i+1)\big(m + (g-1)i\big).$

To obtain $r_i$ from the number (B.31.1), one ought to subtract the excess contributions of the various stationary points of $\Gamma$. For $i = 1$ for instance, for all cusp $p$ of $\Gamma$ (i.e., for all $p \in \bar{\Gamma}$ such that $s_1(p) > 0$), the whole $p^\perp \subseteq \check{\mathbf{P}}^N$ gives rise to hyperplanes contributing to the number (B.31.1), so that the latter actually is the degree of the reducible hypersurface sum of $\Gamma^{\vee_1}$ and the various $p^\perp$ for all cusps $p$ of $\Gamma$. Taking this into consideration, one recovers Formula (B.24) for $r_1$.

For $i = 2$, also the points $p \in \bar{\Gamma}$ such that $s_1(p) = 0$ and $s_2(p) > 0$ contribute to (B.31.1), to the effect that $(\mathrm{Osc}_p^1(\Gamma))^\perp$ occurs together with $\Gamma^{\vee_2}$. In general, all $p \in \bar{\Gamma}$ such that $s_j(p) > 0$ for some $j \leqslant i$ contribute to (B.31.1), and to conclude we need to show that each contributes with the multiplicity

$$s_i + 2\, s_{i-1} + \cdots + i\, s_1$$

in de Jonquières' formula. I have not found this particular computation in the literature, and will leave this as an open exercise to the reader.

## C – Biduality for asymptotic tangent lines by enumerative means

Let $\Lambda$ be an asymptotic tangent line to a surface $S \subseteq \mathbf{P}^3$ at a smooth point $p$. In this appendix, we show that the line $\Lambda^\perp$ is an asymptotic tangent line to the dual surface $S^\vee$ at the point $(\mathbf{T}_p S)^\perp$, by computing the intersection multiplicity of $\Lambda^\perp$ and $S^\vee$ at this point. We also prove several variants, see Theorem (C.2).

To compute this intersection multiplicity, we consider the pencil of hyperplane sections of $S$ parametrized by $\Lambda^\perp$, and count its singular members with the topological formula involving Euler numbers described in [XI, Section **]. The intersection multiplicity we are looking for is the multiplicity of the tangent section $\mathbf{T}_p S \cap S$ as a singular member of this pencil (indeed the singular members of $\Lambda^\perp$ correspond to the intersection points of $\Lambda^\perp$ with $S^\vee$). The upshot is the following lemma.

**(C.1) Lemma.** *Let $S \subseteq \mathbf{P}^3$ be a smooth surface, $\Lambda$ be a line, and $H$ be a hyperplane containing $\Lambda$. Let $\tilde{S} \to S$ be the minimal sequence of blow-ups at reduced points such that the proper transform $\tilde{\Lambda}$ of the pencil $\Lambda^\perp$ is base-point-free. Let $\tilde{C}_H$ and $\tilde{C}_{\mathrm{gen}}$ be the members of $\tilde{\Lambda}$ corresponding to $C_H = H \cap S$ and a general member $C_{\mathrm{gen}}$ of $\Lambda^\perp$, respectively. The intersection multiplicity of $\Lambda^\perp$ and $S^\vee$ at $H^\perp$ equals the difference in Euler numbers $e(\tilde{C}_H) - e(\tilde{C}_{\mathrm{gen}})$.*

The assumption that $S$ is smooth is only to simplify matters: one may adapt the lemma to a surface with given singularities, using the same method. The lemma follows from a direct application of [XI, Lemma **] to the morphism $\tilde{S} \to \mathbf{P}^1$ given by $\tilde{\Lambda}$. The same idea appears in [35] and [3]. With this tool at hand, we may now prove the following biduality statement.

**(C.2) Theorem.** *Let $S \subseteq \mathbf{P}^3$ be a smooth surface, and $p \in S$. The orthogonal line to*
*(a) a tangent at $p$,*
*(b) one of the two asymptotic tangents at a non-parabolic $p$,*
*(c) the asymptotic tangent at a general parabolic point $p$,*
*(d) the asymptotic tangent at a parabolic point $p$ of type $\beta$,*
*respectively, is*
*(a) a tangent at $(\mathbf{T}_p S)^\perp$,*
*(b) one of the two asymptotic tangents at $(\mathbf{T}_p S)^\perp$,*
*(c) the tangent at $(\mathbf{T}_p S)^\perp$ of the cuspidal double curve of $S^\vee$,*
*(d) the common reduced tangent cone at $(\mathbf{T}_p S)^\perp$ of both the ordinary and cuspidal double curves of $S^\vee$.*

*Proof.* We proceed with a case by case analysis, and apply Lemma (C.1) in all cases. We use similar notation (and abuse of notation) in all cases, which we introduce in the first case, and then use freely in the others.

(a) Let $\Lambda$ be a general tangent line to $S$, and call $p$ the contact point. Then all members of $\Lambda^\perp$ are tangent at $p$, and we have to blow-up twice at $p$ to remove the base points: we let $\varepsilon_1 : S_1 \to S$ be the blow-up at $p$, with exceptional divisor $E_1$, and $\varepsilon_2 : S_2 \to S_1$ be the blow-up at the point $p_1 \in E_1$ corresponding to the tangent direction at $p \in S$ defined by $\Lambda$, with exceptional divisor $E_2$.

Locally over $p \in S$, the $\tilde{S} \to S$ of Lemma (C.1) is $\varepsilon = \varepsilon_2 \varepsilon_1 : S_2 \to S$. Let $C$ be the hyperplane class on $S$. Forgetting about the base points of $\Lambda^\perp$ away from $p$, the base-point-free pencil $\tilde{\Lambda}$ of Lemma (C.1) is the complete linear system $|\varepsilon^* C - \varepsilon_2^* E_1 - E_2|$. Our task now is to compute the transform $\tilde{C}_0$ of the curve $C_0 = S \cap \mathbf{T}_p S$ in this pencil.

Since $C_0 \subseteq S$ has a double point at $p$, we have

$$\varepsilon_1^* C_0 = C_0 + 2E_1,$$

where by abuse of notation $C_0$ also denotes the proper transform of $C_0$ in $S_1$. Next, one has

$$\varepsilon^* C_0 = \varepsilon_2^*(C_0 + 2E_1) = C_0 + 2E_1 + 2E_2$$

with the same abuse of notation, this time for both $C_0$ and $E_1$. Finally we find

$$\tilde{C}_0 = \varepsilon^* C_0 - \varepsilon_2^* E_1 - E_2 = \varepsilon^* C_0 - (E_1 + E_2) - E_2 = C_0 + E_1.$$

This is a reduced curve with two ordinary double points, hence it counts with multiplicity 2 as a singular member of $\Lambda^\perp$ (see [XI, subsec. **] for the details of the computation). It follows that $\Lambda^\perp$ is tangent to $S^\vee$ at $(\mathbf{T}_p S)^\perp$.

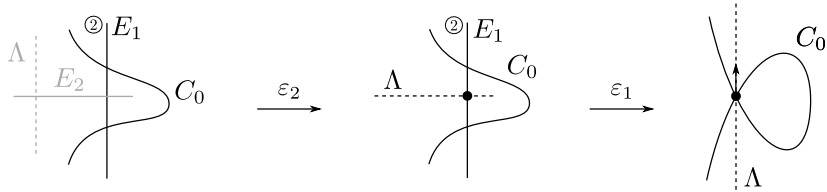

Figure 3: Pencil with base locus an ordinary tangent line $\Lambda$

(b) Let $\Lambda$ be an asymptotic tangent at a general point $p$. Then the members of $\Lambda^\perp$ are all osculating to the order 3 at $p$ (in fact, we know they all have a flex along $\Lambda$), and we need to blow-up three times to remove the base points, with exceptional divisors $E_1, E_2, E_3$. The base-point-free pencil $\tilde\Lambda$ is now

$$|\varepsilon^* C - \varepsilon_3^* \varepsilon_2^* E_1 - \varepsilon_3^* E_2 - E_3| = |\varepsilon^* C - E_1 - 2E_2 - 3E_3|.$$

The curve $C_0 = S \cap \mathbf{T}_p S$ has an ordinary double point at $p$, with one local branch simply tangent to $\Lambda$. We thus have successively

$$\varepsilon_1^* C_0 = C_0 + 2E_1$$
$$\varepsilon_2^* \varepsilon_1^* C_0 = C_0 + 2E_1 + 3E_2$$
$$\varepsilon_3^* \varepsilon_2^* \varepsilon_1^* C_0 = C_0 + 2E_1 + 3E_2 + 3E_3,$$

hence the transform $\tilde C_0$ of $C_0 \subseteq S$ in $\tilde\Lambda$ is

$$\tilde C_0 = \varepsilon^* C_0 - E_1 - 2E_2 - 3E_3 = C_0 + E_1 + E_2.$$

It is reduced with 3 nodes, hence counts with multiplicity 3 as a singular member of $\Lambda^\perp$. This shows that $\Lambda^\perp$ is an asymptotic tangent to $S^\vee$.

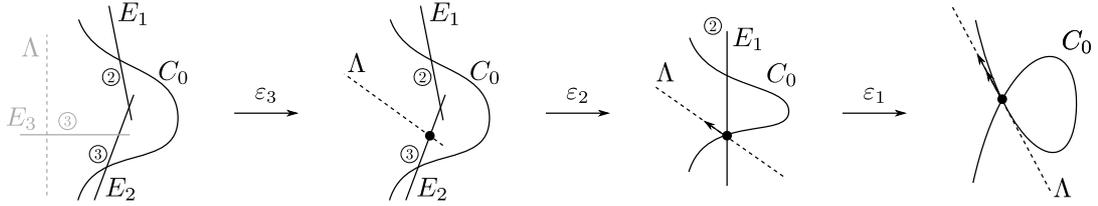

Figure 4: Pencil with base locus an asymptotic tangent line $\Lambda$

(c) Let $\Lambda$ be an asymptotic tangent at a general parabolic point $p$. The situation is the same as in case (b), except that the tangent section $C_0 = S \cap \mathbf{T}_p S$ now has a cusp at $p$, with tangent cone supported on $\Lambda$. We have

$$\varepsilon_1^* C_0 = C_0 + 2E_1$$
$$\varepsilon_2^* \varepsilon_1^* C_0 = C_0 + 2E_1 + 3E_2$$
$$\varepsilon_3^* \varepsilon_2^* \varepsilon_1^* C_0 = C_0 + 2E_1 + 3E_2 + 3E_3,$$

and the transform of $C_0$ in the base-point-free pencil $\tilde\Lambda$ is

$$\tilde C_0 = \varepsilon^* C_0 - E_1 - 2E_2 - 3E_3 = C_0 + E_1 + E_2.$$

It is reduced with a triple point, hence it counts with multiplicity 4 by [XI, Lemma **].

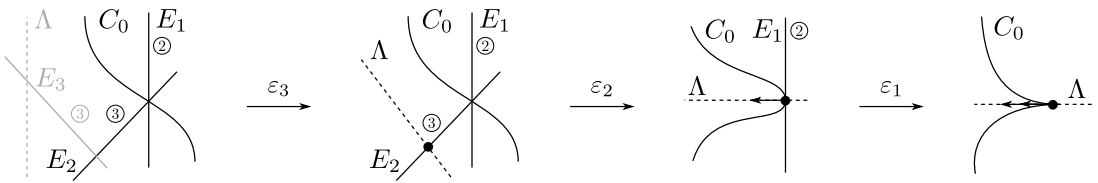

Figure 5: Pencil with base locus the double asymptotic tangent $\Lambda$ at a parabolic point

The point $(\mathbf{T}_p S)^\perp$ lies on the cuspidal double curve of $S^\vee$. A general line through this point intersects $S^\vee$ with multiplicity 2, a general line in the tangent cone of $S^\vee$ with multiplicity 3, and the tangent line of the cuspidal double curve with multiplicity 4. Hence $\Lambda^\perp$ is the latter line.

(d) Lastly, let $p$ be a point such that $C_0 = S \cap \mathbf{T}_p S$ has a tacnode at $p$, and $\Lambda$ be the asymptotic tangent of $S$ at $p$, which is also the support of the tangent cone of $C_0$ at $p$. The base-point-free pencil is defined as in cases (b) and (c). The computations now give

$$\varepsilon_1^* C_0 = C_0 + 2E_1$$
$$\varepsilon_2^* \varepsilon_1^* C_0 = C_0 + 2E_1 + 4E_2$$
$$\varepsilon_3^* \varepsilon_2^* \varepsilon_1^* C_0 = C_0 + 2E_1 + 4E_2 + 4E_3$$

so that

$$\tilde{C}_0 = \varepsilon^* C_0 - E_1 - 2E_2 - 3E_3 = C_0 + E_1 + 2E_2 + E_3.$$

This time the curve is not reduced. However this is not a reason for the topological formula behind Lemma (C.1) to fail, and the latter is still valid. The only caution is that the computation of $e(\tilde{C}_0)$ is no longer liable to [XI, Lemma **]; we compute it directly as the Euler number of $\tilde{C}_{0,\mathrm{red}}$, as indeed the topological formula doesn't see multiplicities, i.e., $e(\tilde{C}_0) = e(\tilde{C}_{0,\mathrm{red}})$.

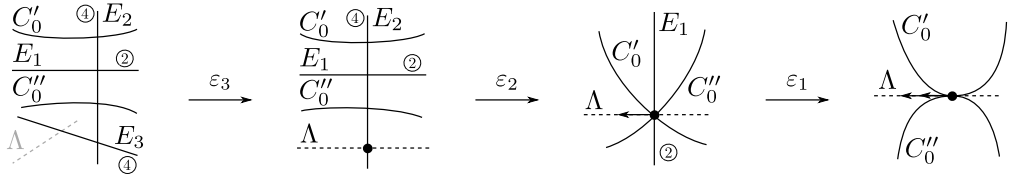

Figure 6: Pencil with base locus the double asymptotic tangent $\Lambda$ at a swallowtail point

Since the curve $C_0 \subseteq S$ has a tacnode, its proper transform in $\tilde{S}$ has genus $g - 2$, where $g$ stands for the genus of a general member of $\Lambda^\perp$. It then follows from the additivity of the Euler number as in the proof of [XI, Lemma **], that

$$e(C_0 + E_1 + E_2 + E_3) = e(C_0) + e(E_1) + e(E_2) + e(E_3) - 4$$
$$= \bigl(2 - 2(g-2)\bigr) + 2 + 2 + 2 - 4 = e(C_\mathrm{gen}) + 6.$$

One then argues as in case (c), remembering the local swallowtail model of $S^\vee$ at $(\mathbf{T}_p S)^\perp$ given in (3.1). The point $(\mathbf{T}_p S)^\perp$ is triple for $S^\vee$, and the only line through this point intersecting $S^\vee$ with multiplicity strictly larger than 4 is the common support of the respective tangent cones of both double curves of $S^\vee$ at $(\mathbf{T}_p S)^\perp$. □

*(C.2.1) Remark.* In the swallowtail model of (3.1), the support of the tangent cones of the double curves of $S^\vee$ is the line $b = c = 0$, which is downright contained in the discriminant hypersurface (3.1.2). Of course it is not the case in general that, in the notation of (d), the line $\Lambda^\perp$ is contained $S^\vee$.

## References

**Chapters in this Volume**[a]

[II]   Th. Dedieu and E. Sernesi, *Deformations of curves on surfaces.*

---